\documentclass{amsart}
\usepackage{amsmath}
\usepackage{amssymb}
\usepackage{amsbsy}
\usepackage{amsfonts}
\usepackage{pict2e}
\usepackage{graphicx}
\usepackage{epstopdf}
\usepackage{verbatim}
\usepackage{subfig}
\usepackage{algorithm}
\usepackage{algpseudocode}
\usepackage{amsmath}
\usepackage{amssymb}

\usepackage{epsfig}
\usepackage{color}
\usepackage[table]{xcolor}

\newcommand{\argmin}[1]{\underset{#1}{\mathrm{argmin}}}

\usepackage{color}
\usepackage{varwidth}
\usepackage[normalem]{ulem}

\usepackage[table]{xcolor}

\usepackage{graphicx}
\usepackage{epstopdf}
\usepackage{subfig}
\usepackage{epsfig}
\usepackage{fullpage}
\newcommand{\supp}{\operatorname{supp}} 

\newtheorem{theorem}{Theorem}[section]

\theoremstyle{definition}

\theoremstyle{remark}
\newtheorem{remark}[theorem]{Remark}

\numberwithin{equation}{section}

\begin{document}

\title{PDEs with Compressed Solutions}

\author{Russel E. Caflisch}

\author{Stanley J. Osher}

\author{Hayden Schaeffer}

\author{Giang Tran}

\begin{abstract}
Sparsity  plays a central role in recent developments in signal processing, linear algebra, statistics, optimization, and other fields. In these developments, sparsity is promoted through the addition of an $L^1$ norm (or related quantity) as a constraint or penalty in a variational principle. We apply this approach to partial differential equations that come from a variational quantity, either by minimization (to obtain an elliptic PDE) or by gradient flow (to obtain a parabolic PDE).  Also, we show that some PDEs can be rewritten in an $L^1$ form, such as the divisible sandpile problem and signum-Gordon. Addition of an $L^1$ term in the variational principle leads to a modified PDE where a subgradient term appears. It is known that modified PDEs of this form will often have solutions with compact support, which corresponds to the discrete solution being sparse.   We show that this is advantageous numerically through the use of efficient algorithms for solving $L^1$ based problems.  
\end{abstract}

\maketitle

\section{Introduction}\label{intro_sec}

Sparsity has played a central role in recent developments in fields such as signal processing, linear algebra, statistics and optimization. 
Examples include compressed sensing \cite{Candes:2006robust,donoho2006compressed}, matrix rank minimization \cite{recht2010guaranteed}, phase retrieval \cite{candes2013phase} and robust  principal component analysis \cite{candes2011robust,AspremontSPCA2007,Qi:2013sparse},
as well as many others. A key step in these examples is the use of an $L^1$ norm (or related quantity) as a constraint or penalty term in a variational formulation. 
In all of these examples, sparsity is for the coefficients (\textit{i.e.}, only a small set of coefficients are nonzero) in a well-chosen set of modes  for representation of the corresponding vectors or functions.

The use of sparse techniques in physical sciences and partial differential equations (PDEs) has been limited, but recent results have included  numerical solutions of PDEs with multiscale oscillatory solutions \cite{Schaeffer:2013sparse}, efficient material models derived from quantum mechanics calculations \cite{Ozolins2013compressive},  ``compressed modes" for variational problems in mathematics and physics \cite{Ozolins:2013CMs}, and ``compressed plane waves" \cite{Ozolins:2013CPWs}. In the latter two examples, sparsity is used in a new way, in that the solutions are sparse and localized in space (as opposed to sparsity of the coefficients in some modal representation). Sparse solutions with respect to low-rank libraries are used in modeling and approximating dynamical systems, see for example \cite{brunton2013compressive}.

Motivated by these works and by the early theoretical framework established in \cite{BrezisBook,brezis1974monotone,brezis1974solutions,brezis1976estimates}, we investigate PDEs with $L^1$ subdifferential terms. The PDE is either an elliptic PDE  coming from a variational principle or a parabolic PDE coming from a gradient flow of a convex functional. In either case, the $L^1$ term in the convex functional leads to a subgradient term in the PDE. Fortunately, the subgradient term has a simple explicit form, so that the PDEs are amenable to analysis and computation.

The goal of this work is to present fast computational schemes for these modified PDEs, provide some additional theoretical insights, and show some connections to known physical equations. Our starting point is the convex functional:
\begin{equation}
E(u) = \int \frac{1}{2} (\nabla u) \cdot M (\nabla u)  - uf + \gamma |u| dx, 
\label{variational_function}
\end{equation}
where $\gamma\geq 0$, $M=M(x)$ is a symmetric, positive definite matrix as a function of $x$,  and $f=f(x)$ or $f=f(x,t)$ will be a specified function depending on $x$ or $(x,t)$. Define the partial differential operator $Au= - \nabla \cdot ( M \nabla u)$.  Minimization of $E(u)$ for $f=f(x)$ leads to the following elliptic PDE
\begin{equation}
\begin{aligned}
  Au= f - \gamma p(u),
\label{general_elliptic_pde}
\end{aligned}
\end{equation}
and gradient descent $\partial _t u = - \partial_u E(u)$, starting from initial data $g(x)$, leads to the following parabolic PDE
\begin{equation}
\begin{aligned}
  u_t + Au= f- \gamma p(u)\\
u(x,0)  =g(x),
\label{general_parabolic_pde}
\end{aligned}
\end{equation}
in which $p(u)$ is a subgradient of $\|u\|_{L^1}$, \textit{i.e.}, $\|v\|_1\geq \|u\|_1 + \langle {v-u,p(u)}\rangle$, for any $u$ and $v$, where $\langle{\,,\,}\rangle$ denotes the $L^2$ inner product. 

The paper is divided as follows: in Section~\ref{formulation_sec},  we provide the general formulation of the problem. In Section~\ref{various_sec}, we review known results and present various properties of solutions to the modified PDEs. The numerical implementation and simulations are presented in Sections~\ref{numerics_sec} and ~\ref{comp_sec}, and we conclude in Section~\ref{conc_sec}.  


\section{Problem Formulation}\label{formulation_sec}
The problem we consider in this work is to numerically solve the following PDE 
\begin{equation}
\begin{aligned}
  u_t + Au= f- \gamma p(u)\\
u(x,0)  =g(x),
\end{aligned}
\end{equation}
and to verify theoretical results. The difficulty with such equations is the multivalued nature of the subgradient term. Fortunately for this type of equation, we can explicitly identified the subgradient as 
\begin{equation}
p(u) = \begin{cases}
  \textrm{sign} (u) & \mbox{if $|u|>0$ } \\
\argmin{|q| \leq 1} |f-\gamma q|&  \mbox{if $u=0$. }
\end{cases} 
\end{equation}
Note that if $u=0$ and $|f(x)| \leq \gamma$, then $p=f(x)/\gamma$.
This specification for $u$ was proved in general in \cite{AubinBook,BrezisBook}. It can be shown directly from Equations \eqref{general_elliptic_pde} and \eqref{general_parabolic_pde}, as follows. For $u=0$ in an open set, the left side of the equations is $0$ so that $f(x)-\gamma p(u)=0$, which is only possible if $f(x)\leq \gamma$ and $p(u)=f(x)/\gamma$. The value of $p(u)$ on a lower dimensional set does not matter, since the value of the forcing terms on a lower dimensional set does not affect the solution $u$ of the differential equations. For the elliptic equation \eqref{general_elliptic_pde} one can also show directly that this identification of $p(u)$ gives $u=0$ as the unique minimizer of $E(u)$ (see Appendix).

\section{Various Properties} \label{various_sec}
In this section we recall the established existence theory for the elliptic equation \eqref{general_elliptic_pde} and the parabolic equation \eqref{general_parabolic_pde}, and provide some further insights to the behavior of solutions.

\subsection{Review of Theoretical Results}
Equation~\eqref{general_elliptic_pde} is related to the general class of elliptic equation:
$$ -\Delta u = F(u),$$
where $F$ contains a discontinuous component. The existence and uniqueness of the solution $u$ are studied in \cite{komura1969differentiability, kato1970accretive, crandall1969semi}. Solutions also satisfy the standard maximum and comparison principles given the correct sign of $F$. The solutions are compactly supported in both the elliptic and parabolic case, under  some additional conditions \cite{brezis1974solutions, brezis1976estimates}. For the parabolic equation, the solutions are Lipschitz continuous  and right differentiable in time. Furthermore, solutions exhibit finite speed of propagation \cite{brezis1976estimates}. More precisely,  let $S(t)$ be the support set of $u(x,t)$, then for small times $t$:
\begin{itemize}
\item if $u(x,0)$ does not vanish on $\partial S(0)$ , then $$S(t) \subset S(0) + B(c \sqrt{t \log(t) }),$$
\item if $u(x,0)$ and $\nabla u(x,0) $ vanishes on $\partial S(0)$, then $$S(t) \subset S(0) + B(c \sqrt{t }),$$
\end{itemize}
where $B(r)$ is the ball of radius $r$ centered at the origin. In a simple case, we can construct the exact bounds in order to verify the convergence of the method to a known solution.

At a number of places in the manuscript, we will simplify the presentation by assuming that $x\in\mathbb{R}^1$ and that $M=1$, so that the elliptic PDE \eqref{general_elliptic_pde} becomes Laplace's equation with nonlinear forcing:
\begin{equation}
u_{xx} = -f +\gamma p(u),
\label{elliptic_pde}
\end{equation}
and the parabolic PDE \eqref{general_parabolic_pde} becomes the heat equation with nonlinear forcing:
\begin{equation}
\begin{aligned}
 u_t - u_{xx} = f - \gamma p(u),  
\label{parabolic_pde}
\end{aligned}
\end{equation}
\subsection{A Free Boundary Formula}\label{freebndysection}
 In 1D, consider the following equation
\begin{equation}
\begin{aligned}
	u_t - u_{xx} &=
	\begin{cases}
       & f(x) -\gamma,\quad |x|<a(t)\\ 
       & 0,\hspace{1.5 cm}|x|>a(t) 
	\end{cases}\\
u(x,0) &= 0.
\end{aligned}
\label{eq:freeBoundary}
\end{equation}
For simplicity assume that $f(x) = f(|x|)$ and $f$ is a decreasing function with $f(|x|)\rightarrow 0$ as $|x|\rightarrow\infty$. Denote $a_0\geq 0$ such that $f(a_0) =\gamma$ and assume that $f_x(a_0) \neq 0$. Then, the free boundary's endpoint is governed by (for small time $t$): 
\begin{equation}
a(t) = a_0+a_1\sqrt {t},\label{a(t)} 
\end{equation}
for some $a_1\geq 0$ (for the proof, see Appendix). A similar result holds for zero force and non-zero (finitely supported) initial data.

\subsection{Support Size}\label{sec_compsupp}
Since it is known that the support is compact, we would like to estimate its size. In fact, by integrating Equation~\eqref{general_elliptic_pde} (see Appendix), the support of $u$ satisfies
\begin{equation}
 |\text{supp}(u)| \leq \gamma^{-1} {\int_{\supp(u)} |f| dx}.
 \label{supportBound1}
\end{equation}
A slight modification of (\ref{supportBound1}) shows that for any nonnegative $\alpha$ and $\beta$ with $\alpha+\beta=1$, we have
\begin{equation}
 |\text{supp}(u)| \leq (\alpha \gamma)^{-1} {\int (|f|-\beta\gamma)^+ dx}.\label{supportBound2}
\end{equation}
In this inequality, the superscript $+$ denotes the positive part; \textit{i.e.}, $(x)^+=\max(x,0)$. For the parabolic case, a similar bound on the support size holds:
\begin{align}
| \supp_{(x,t)} u(x,t)|  \leq (\alpha \gamma)^{-1} \left( \int  |g| dx
+ \iint  (|f|-\beta\gamma)^+ dx \, dt \right),
\label{para_bound2}
\end{align}
 for any nonnegative $\alpha$ and $\beta$ with $\alpha+\beta=1$. 
 
 \subsection{$L^1$ Contraction and Total Variation Diminishing}

Let $u$ and $v$ be solutions of Equation~\eqref{parabolic_pde} with initial data $g(x)$ and $h(x)$, respectively. First, note that for any subgradient $p$ of a convex functional, we have 
\begin{equation}
\text{sign}(u-v) (p(u)-p(v)) \geq 0.
\label{sign_sub}
\end{equation}
We wish to show that the solutions are $L^1$ contractive and TVD by computing the following:
\begin{align*}
\frac{d}{dt} ||u-v||_{L^1}& = \frac{d}{dt}\int\limits_{|u-v|>0} |u-v| dx \\
&= \int\limits_{|u-v|>0} \text{sign}(u-v) (u_t-v_t) dx \\
&= \int\limits_{|u-v|>0} \text{sign}(u-v) (u-v)_{xx} - \gamma \ \text{sign}(u-v) (p(u)-p(v)) dx .
\end{align*}
The first term is zero by the divergence theorem and the second term is negative by Equation~\eqref{sign_sub},  so we have $\frac{d}{dt} ||u-v||_{L^1}\leq 0$, and thus the modified PDE is an $L^1$ contraction. Moreover, if we take $h(x)=g(x+\delta)$ for any $\delta>0$ we have
\begin{align*}
\frac{d}{dt} \|u(x,t)-u(x+\delta,t)\|_{L^1} \leq 0. 
\end{align*}
Dividing the equation above by $\delta$ and taking the supremum over all $\delta$, the following inequality holds:
\begin{align*}
\frac{d}{dt} ||u||_{TV}  \leq 0.
\end{align*}
Therefore, Equation~\eqref{parabolic_pde} is TVD.

\subsection{Entropy Condition}

The $L^1$ contraction and TVD results are directly analogous to those that are obtained by solving the viscosity regularized nonlinear conservation laws:
\begin{equation*}
w^{\epsilon}_t= \epsilon \,w^{\epsilon}_{xx} -f(w^{\epsilon})_x,
\end{equation*}
for $\epsilon>0$. Then by letting $\epsilon \rightarrow 0$, one recovers the unique inviscid limit, see \cite{lax1973hyperbolic}.  

We can also easily obtain an ``entropy inequality'' in the same spirit. Consider the scaled modified heat equation:
\begin{equation}
u_t= \epsilon \,u_{xx} -\gamma p(u).
\label{eps_eqn}
\end{equation}
We deliberately put an $\epsilon$ in front of the diffusion term to emphasize the similarities to the theory of scalar conservation laws. The following argument holds in more general cases.

Let $ K(u)$ be a convex function of $u$ with subgradient $q(u)$. Multiplying Equation~\eqref{eps_eqn} by the subgradient (as in \cite{lax1973hyperbolic}) yields:
\begin{equation}
  \frac{d}{dt}K(u) \leq \epsilon \, \frac{d^2}{dx^2}K(u) -\gamma q(u)p(u).
\label{entropy_ineq}
\end{equation}
For example, if $K(u)=|u|$, then whenever $u \neq 0$, we have
\begin{equation}
 |u|_t \leq \epsilon \,|u|_{xx} -\gamma. \label{absolutevalue}
\end{equation}

We integrate Equation~\eqref{entropy_ineq} over the region $\mathcal{S}(t)$, the support set of $u(x,t)$ defined in Section~\ref{sec_compsupp}, to get
\begin{equation}
 \frac{d}{dt} \int_{\mathcal{S}(t)} K(u) \, dx \leq -\gamma \int_{\mathcal{S}(t)} q(u)p(u) \, dx ,
\label{entropy_ineq2}
\end{equation}
since the spatial gradient is zero along the boundary. By choosing $K(u)=\frac{1}{a} |u|^a$ for $a\geq1$, Equation~\eqref{entropy_ineq2} provides $L^a$ estimates of the solutions. Furthermore, if $K(u)=(u-c)^+$ for $c>0$, then 
\begin{equation}
 \frac{d}{dt} \int_{\mathcal{S}^+_c(t)} (u-c)^+ \, dx \leq -\gamma |\mathcal{S}^+_c(t)|,
\end{equation}
where $\mathcal{S}^+_c(t)$ is the set of $x$ for which $u(x)>c$. 

\subsection{Regularity}

We can show that the solutions of the Laplace's equation \eqref{elliptic_pde}  and of the heat equation \eqref{parabolic_pde} are smooth. Let $\Omega_+$, $\Omega_-$, and $\Omega_0$ denote the sets $\{u>0\}$, $\{u<0\}$ and $\{u=0\}$, respectively. Then the solution $u$ of the Laplace's equation \eqref{elliptic_pde} can be represented by 
\begin{equation}
  u(x) =   \int_{\Omega_+} G(x-y)(f(y)-\gamma) dy  +  \int_{\Omega_-} G(x-y)(f(y)+\gamma) dy, 
\label{elliptic_solution}
\end{equation}
and the solution of the heat equation \eqref{parabolic_pde} can be written as 
\begin{equation}
\begin{aligned}
  u(x,t) = \int G(x-y,t) g(y)dy + \int_0^t \int_{\Omega_+(s)} G(x-y,t-s)(f(y)-\gamma) dy ds \\
\quad
+ \int_0^t \int_{\Omega_-(s)} G(x-y,t-s)(f(y)+\gamma) dy ds,
\label{heat_solution}
\end{aligned}
\end{equation}
in which  the Green's function $G(x,t)$ for the heat equation
and the Green's function $G(x)$ for the Laplace's equation are given by
\begin{equation}
\begin{aligned}
& G(x)=|x|/2,\\
 & G(x,t) = (4\pi t)^{-1/2} \exp(-x^2/4t).
\label{GreensFunctions}
\end{aligned}
\end{equation}
 From these formulas, if $f$ is continuous, then one can see that $u$ is $C^2(x)$ and $C^1(t)$ away from $u=0$ and that  $u$ is $C^1(x)$ everywhere.

\subsection{Traveling Wave}
To demonstrate finite speed of propagation,  consider the 1D-traveling wave solution $u(x,t) = v(s)$ for 
$s=x-\sigma t$, of the Equation~\eqref{parabolic_pde} with no forcing term.
To be specific, we will assume that $v(s)\geq 0$ for $s\geq 0$ and $v(s)=0$ for $s\leq 0$.
We see that $v$ must satisfy the ODE
\begin{equation}
v_{ss} +\sigma v_s -\gamma =0, \label{travelingwave}
\end{equation}
subject to the conditions
\begin{equation}
 v(0) = v'(0) =0.
\end{equation}
The general solution of Equation~\eqref{travelingwave} is 
\begin{eqnarray}
v(s)= \begin{cases}
 & \frac{\gamma}{\sigma}s + c_1 e^{-\sigma s} +c_2, \quad s\geq 0 \\
& 0,\quad\text{otherwise}.
\end{cases}
\label{tw_profile}
\end{eqnarray}
The boundary conditions imply
\begin{equation*}
c_1 = -c_2= \frac{\gamma}{\sigma ^2},
 \end{equation*}
so that the traveling wave solution of Equation~\eqref{parabolic_pde} is
\begin{eqnarray*}
 u(x,t) = \begin{cases}
           &  \frac{\gamma}{\sigma}(x-\sigma t) + \frac{\gamma}{\sigma ^2}\left(e^{-\sigma (x-\sigma t)} -1\right),\quad x\geq \sigma t\\
           &  0,\quad\text{otherwise}.
          \end{cases}
\end{eqnarray*}
We see that in this case we have one sided support.

\begin{remark}
This traveling wave solution is used as a reference solution to compute the error for our numerical scheme (see Section~\ref{numTW}). Also, the simple analytic form shows that solutions with non-trivial support sets are easy to find in the modified PDE. 
\end{remark}
\subsection{An Exact Solution}
We construct the exact solution of Equation~\eqref{elliptic_pde} with nonnegative force $f=(1+x^2)^{-3/2}$ and $\gamma\in [0,1]$.   The exact solution is given explicitly by:
\begin{equation*}
u = \begin{cases}
&-(1+x^2)^{1/2} +\dfrac{1}{2} \gamma x^2 +c,\quad |x| \leq a\\
&0,\quad |x|>a.
\end{cases}
\end{equation*}
where,
\begin{equation*}
c = \dfrac{\gamma + \gamma^{-1}}{2},\quad a = \sqrt{\gamma^{-2} -1}.
\end{equation*}
The boundary value $a$ and constant $c$ are determined so that $u(\pm a)=u_x(\pm a)=0$. At the boundary of the support, $f(\pm a)=\gamma^3 < \gamma$. The results show that the solution is nonnegative for nonnegative $f$, and that having $|f(x)| \leq \gamma$ does not imply $p(u(x))=\frac{f(x)}{\gamma}$.


\section{Numerical Implementation}\label{numerics_sec}

Given an elliptic operator $A$, we would like to solve problems of the form:
\begin{align}
Au+\partial\|u\|_{L^1} \ni f
\end{align}
or
\begin{align}
u_t+Au+\partial\|u\|_{L^1} \ni f
\end{align}
which corresponds to the elliptic or parabolic equations, respectively. We will present two methods to do so.  The first scheme is semi-implicit (also known as implicit-explicit or proximal gradient method), where the subgradient term is discretized forward in time and the diffusion term is lagged. We apply this method to solve the time dependent equations. The second scheme is the Douglas-Rachford method, which we use to solve both the elliptic problem and the parabolic problem. Both methods can handle the multivalued nature of the subgradient $\partial\|u\|_{L^1}$. In this section, we denote $h$ and $\tau$ the space and time steps of the finite difference schemes.

\subsection{Implicit-Explicit Scheme (Proximal Gradient Method)}

From the numerical perspective, the multivalued term $\partial\|u\|_{L^1}$ is the main source of difficulties, since the value is ambiguous. However, an operator of the form $I +\sigma \, \partial F$ ( where $F$ is convex) has an easy-to-compute inverse. The inverse operator  $(I + \sigma \, \partial F)^{-1}$, also known as the resolvent or proximal operator, $\text{prox}_{\sigma F}(\cdot)$, can be found by solving the following optimization:
\begin{align}
(I + \sigma\,  \partial F)^{-1}(z)= \argmin{v}  \ \frac{1}{2}||v-z||_{L^2}^2+ \sigma F(v).
\label{proximal_operator}
\end{align}
For example, if  $F(u) = ||u||_{L^1}$ and thus $\partial F(u) = \partial\|u\|_{L^1}$, we have:
\begin{align*}
(I + \sigma\, \partial\|\cdot\|_{L^1} )^{-1}(z)&= \argmin{v}  \ \frac{1}{2}||v-z||_{L^2}^2+ \sigma ||v||_{L^1}\\
&=S(v,\sigma),
\end{align*}
where the shrink operator, $S$, is defined point-wise as $S(v,\sigma):= \max(|v|-\sigma,0) \frac{v}{|v|}$. 

Using the proximal operator, we will write the discretization of Equation~\eqref{general_parabolic_pde} in a semi-implicit form. We first discretize Equation~\eqref{general_parabolic_pde}  in time:
\begin{align}
u^{n+1}-u^{n}+\tau Au+\tau \partial\|u\|_{L^1} \ni \tau f.
\end{align}
 Then to apply the proximal gradient method, the last two terms on the left are evaluated at $n$ and $n+1$ as follows:
  \begin{align}
u^{n+1}-u^{n}+\tau Au^n+\tau \partial\|u^{n+1}\|_{L^1} \ni \tau f.
\end{align}
The resulting iterative scheme is:
  \begin{align}
u^{n+1}=S(u^n-\tau Au^n+\tau f, \tau).
\end{align}
For example,  for the heat equation, where $A=-\Delta$, the iterative scheme is:
 \begin{align}
u^{n+1}=S(u^n+\tau \Delta u^n+\tau f, \tau),
\label{num_scheme}
\end{align}
and is convergent given $\tau \leq \frac{h^2}{4}$. This scheme has the same complexity as the corresponding standard explicit methods for PDE.

\subsection{Alternating Direction Implicit (Douglas-Rachford) Method}
The Douglas-Rachford algorithm for nonlinear multivalued evolution equation was studied in \cite{lions1979splitting}. Denote $Bu: = \partial \|u\|_{L^1}$, the iterative scheme for Equation~\eqref{general_parabolic_pde} is
\begin{equation}
u^{n+1} = (I + \tau B )^{-1}\left [ (I+\tau A)^{-1}(I-\tau B) + \tau B \right ]u^n,
\end{equation}
which can be rewritten as:
\begin{equation}
\begin{aligned}
u^{n+1}&=(I + \tau B)^{-1}\tilde{u}^{n}\\
\tilde{u}^{n+1} &= \tilde{u}^n +(I+\tau A)^{-1}(2u^{n+1} - \tilde{u}^n) - u^{n+1}.
\label{Douglas_Rachford}
\end{aligned}
\end{equation}
It was shown that the method is unconditionally stable and convergent for all $\tau>0$ \cite{DR1,lions1979splitting,DR2}. Also, note that the iterates $u^n$ converges to a solution of the stationary equation \eqref{general_elliptic_pde}. For the sandpile problem \cite{LevineSandpile}, the operators $A$ and $B$ are chosen specifically as follows:
\begin{equation}
Au = -\Delta u - f  ,\quad Bu = \partial \|u\|_{L^1},
\end{equation} 
so that the operation for $u^{n+1}$ in the iterative process, Equation~\eqref{Douglas_Rachford}, is a shrink. The corresponding proximal operators are 
\begin{align*}
\text{prox}_{\tau F}(z) &=(I+\tau A)^{-1}(z) = (I-\tau\Delta)^{-1} (z+\tau f)\\
\text{prox}_{\tau G}(z) &= (I+\tau B)^{-1}(z) =  S(z,\tau),
\end{align*}
where $F(z) = \frac{1}{2}\|\nabla z\|_{L^2}^2 -\langle f,z \rangle $ and $G(z) = \|z\|_{L^1}$. To compute $(I-\tau \Delta)^{-1}$ numerically, we use the FFT, where the discrete Laplacian $\Delta_h u$ is viewed as the convolution of  $u$ with the finite difference stencil.
\begin{remark}
Since the shrink operator is the last step of the iterative process, this method provides a numerically well-defined support set for $u$, making it easier to locate the free boundary. 
\end{remark}

\section{Computational Simulations} \label{comp_sec}
In this section we show convergence of our numerical scheme to known solutions, approximations to the support set evolution, and numerical solutions for higher dimension.

\subsection{Numerical Convergence}\label{numTW}
\begin{figure}[t!]
\centering
\subfloat[$t=0$]{\includegraphics[width = 2in]{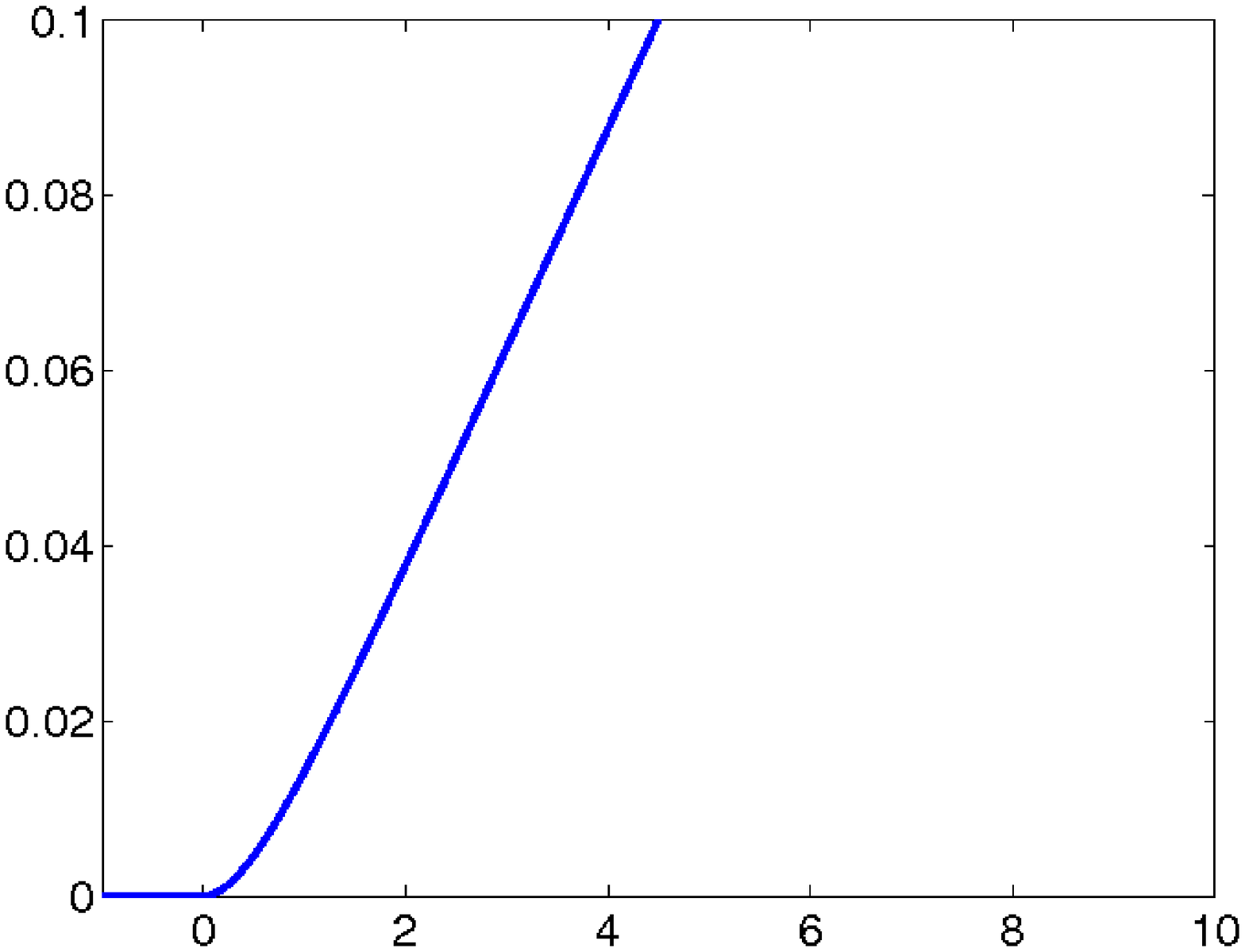}} \quad
\subfloat[$t=0.398$]{\includegraphics[width = 2in]{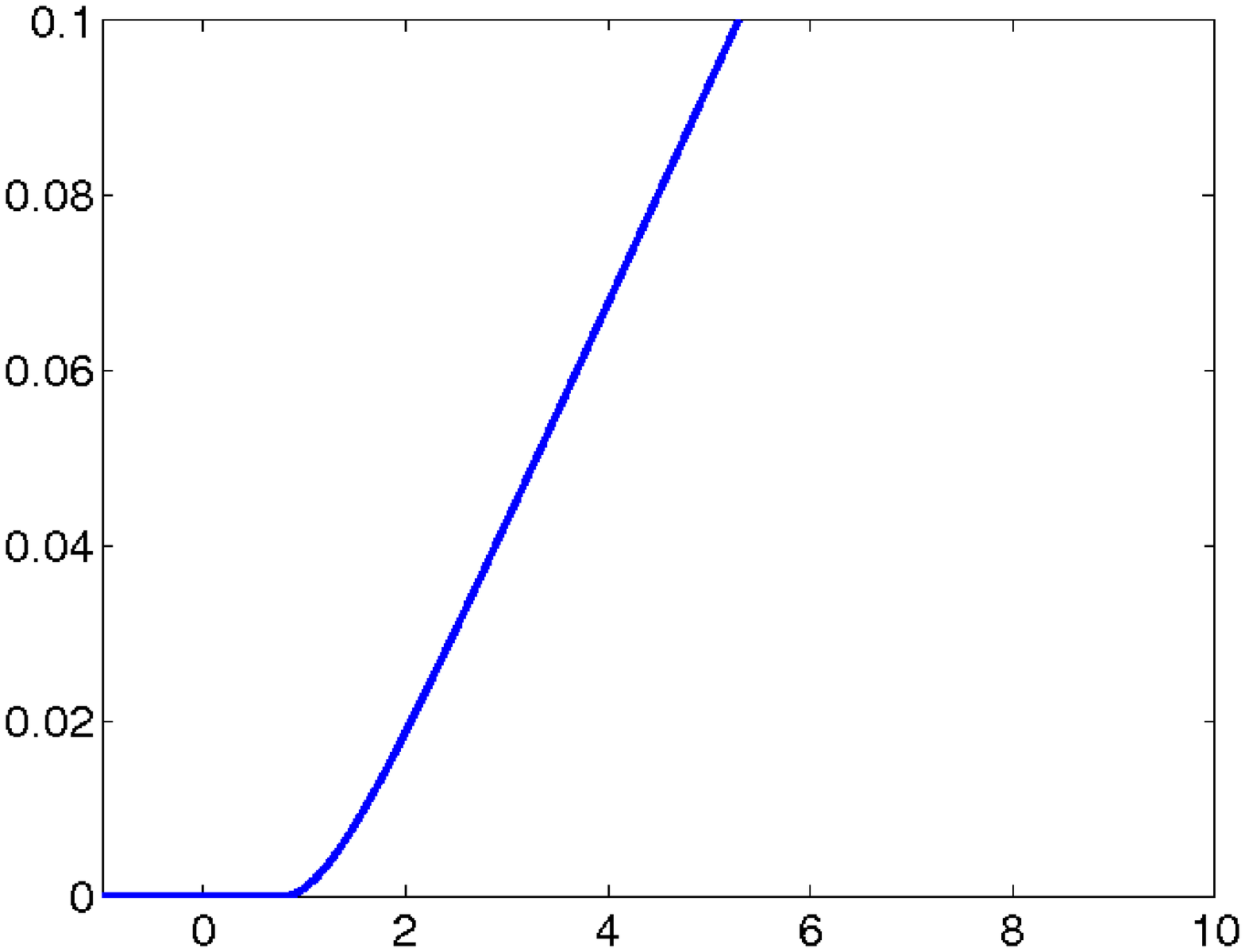}}\\
\subfloat[$t=1.594$]{\includegraphics[width = 2in]{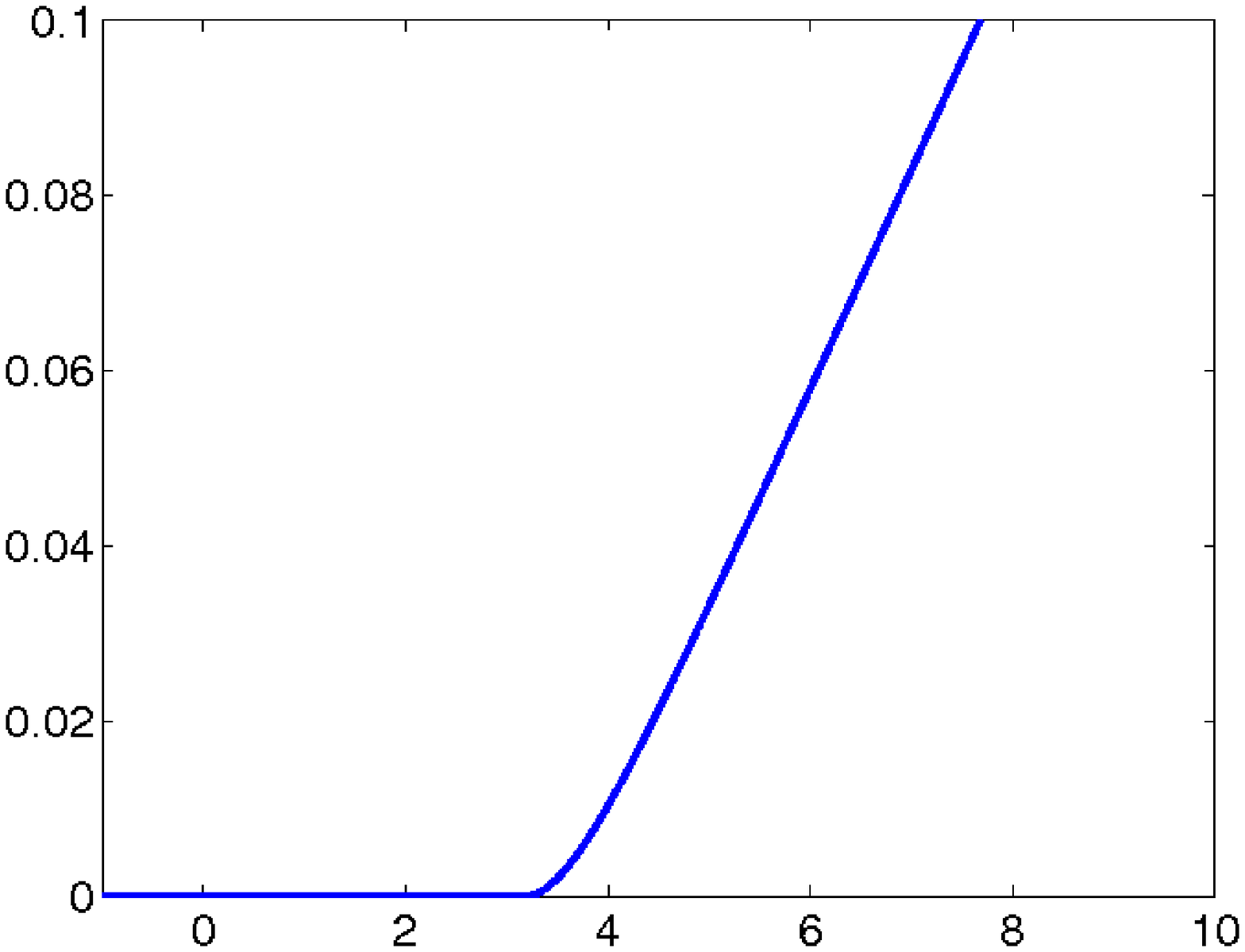}}\quad
\subfloat[$t=2.390$]{\includegraphics[width = 2in]{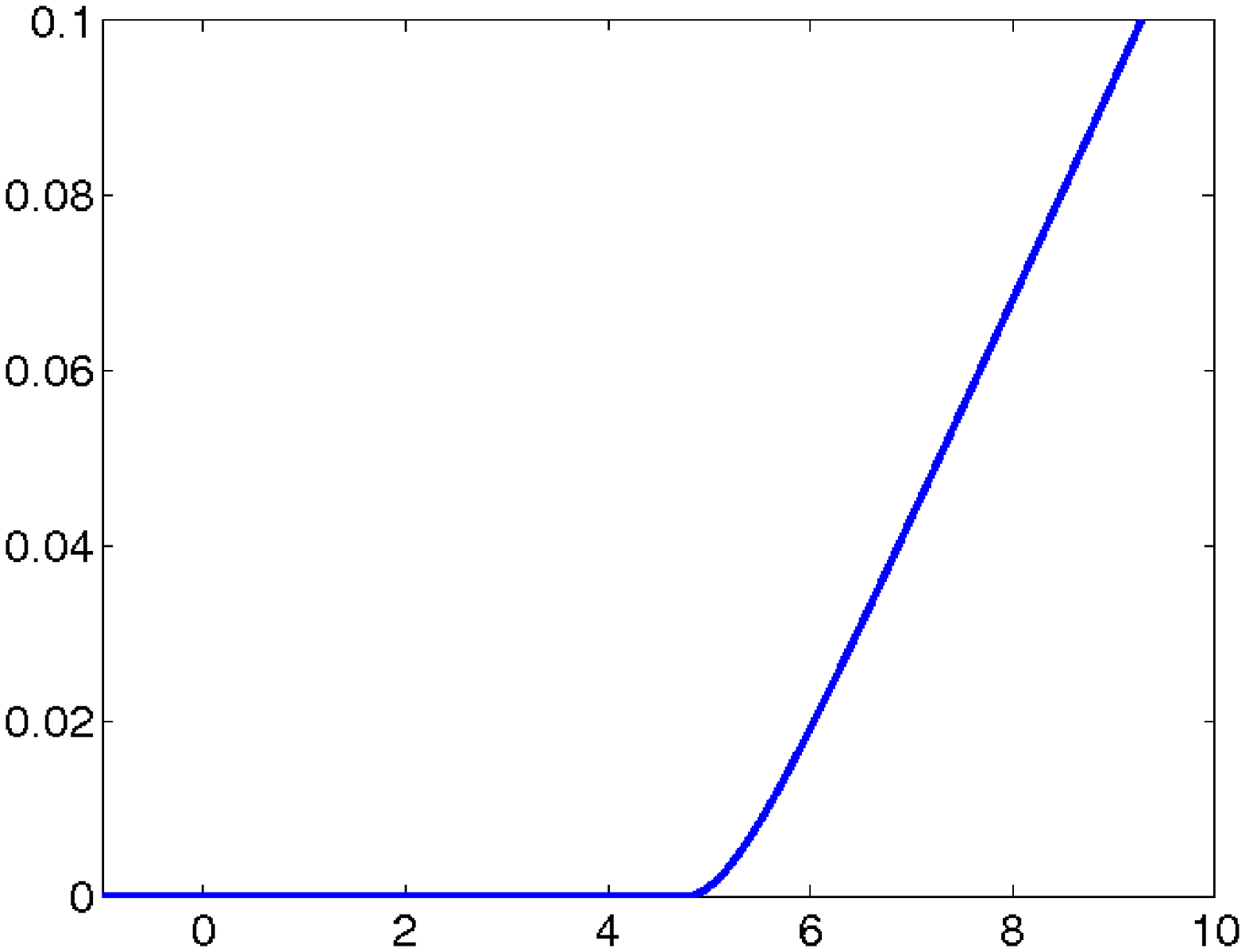}} 
\caption{Numerical solution starting with an initial traveling wave profile with $\sigma=2$ and $\gamma=0.05$ computed using  $500$ grid points.}
\label{fig:TW}
\end{figure}

In Figure~\ref{fig:TW}, we solve Equation~\eqref{parabolic_pde} (with $\gamma=0.05$) using the implicit-explicit scheme (Equation~\eqref{num_scheme}). The initial data is taken to be the traveling wave profile (Equation~\eqref{tw_profile}) with speed $\sigma=2$. The numerical solution has the correct support set and speed of propagation, validating the traveling wave solution as well as the numerical method. 

\begin{figure}[t!]
\centering
  \includegraphics[width=3in]{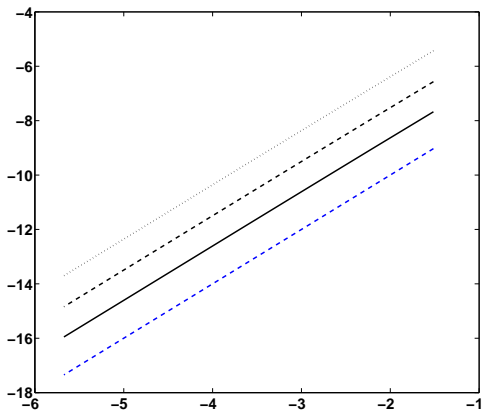}
\caption{Convergence analysis using the $L^1$ (dotted line), $L^2$ (dashed line), and $L^\infty$ (solid line) norms in space and $L^\infty$ norm in time. The $x$-axis is the $log$ of the grid resolution $h$ and the $y$-axis is the $log$ of the Error. The blue dashed line represents second order  convergence.}
  \label{fig:converge}
\end{figure}

This is further confirmed in Figure~\ref{fig:converge}, where the numerical solution is compared to the exact solution. To compute the error, we use the following norms: 
\begin{equation*}
\text{Error}_q(h) = \max_{n} ||u^n_h - u_{exact}||_{q},
\end{equation*}
where $q=1,2,\infty$ and $u^n_h$ is the solution at $t_n$ with space resolution $h$. The errors in these three norms are plotted along side the line representing the second order (dashed line) convergence. 

\begin{figure}[t!]
\centering
\subfloat[$t=0$]{\includegraphics[width = 2.1in]{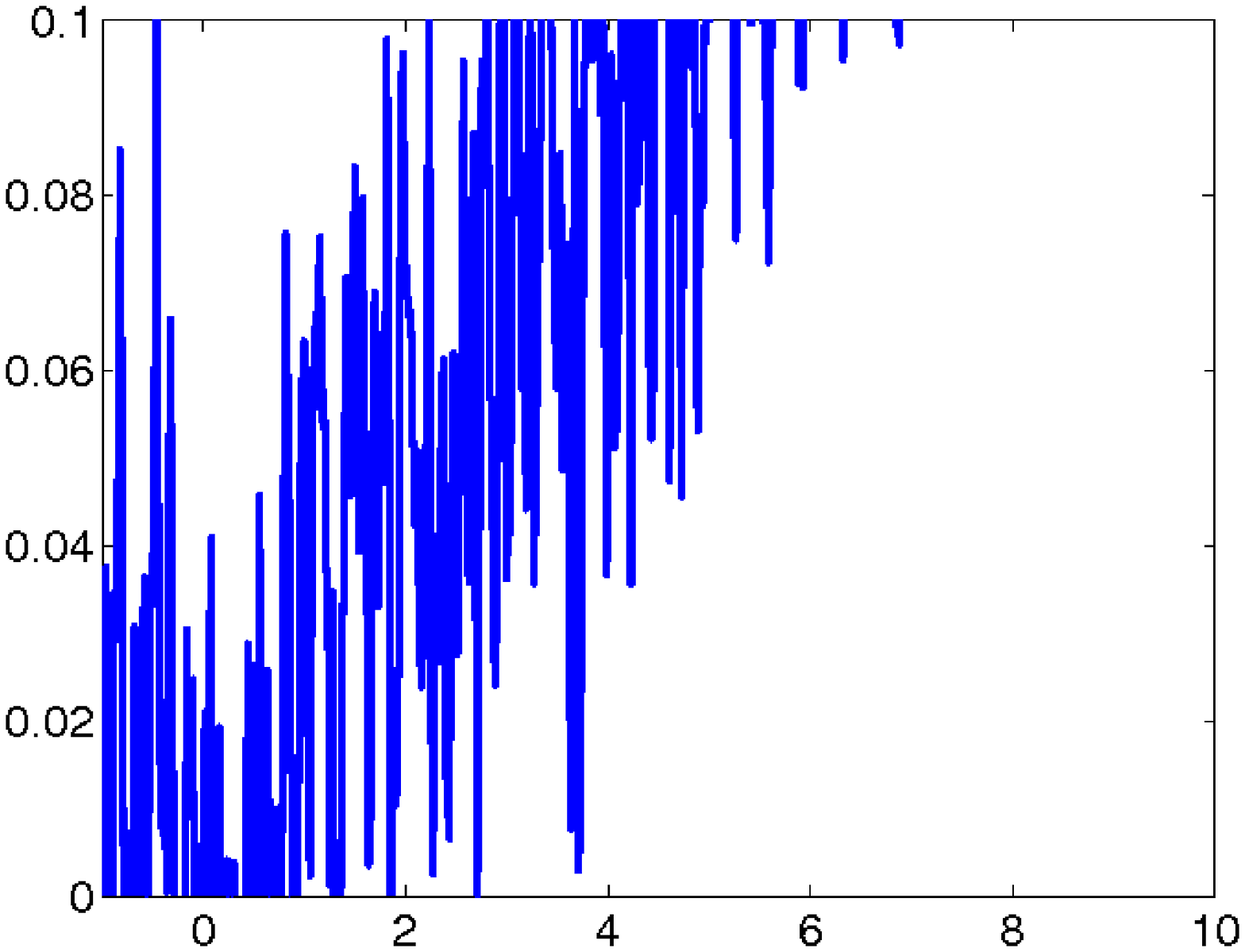}} \quad
\subfloat[$t=0.004$]{\includegraphics[width = 2.1in]{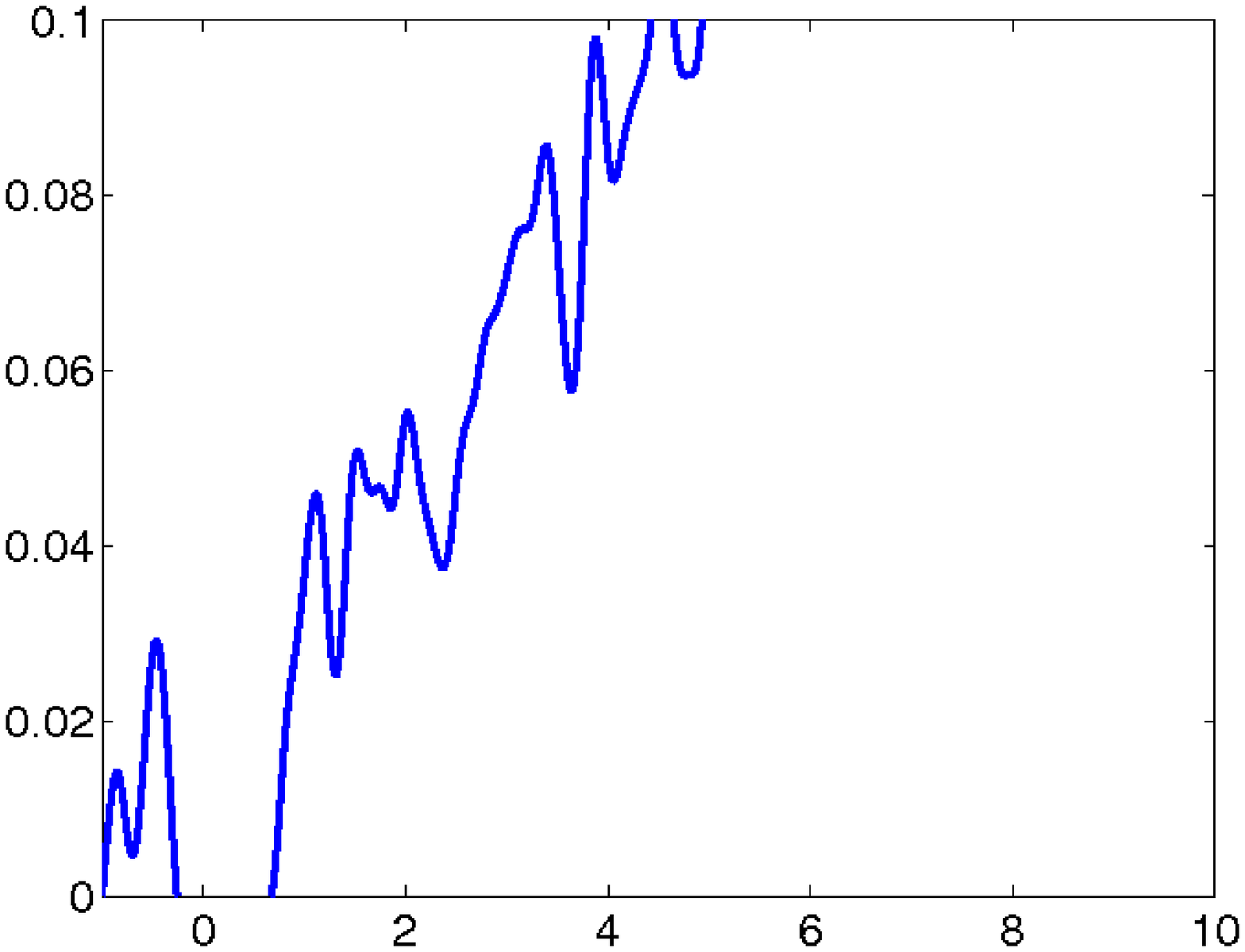}}\\
\subfloat[$t=0.040$]{\includegraphics[width = 2.1in]{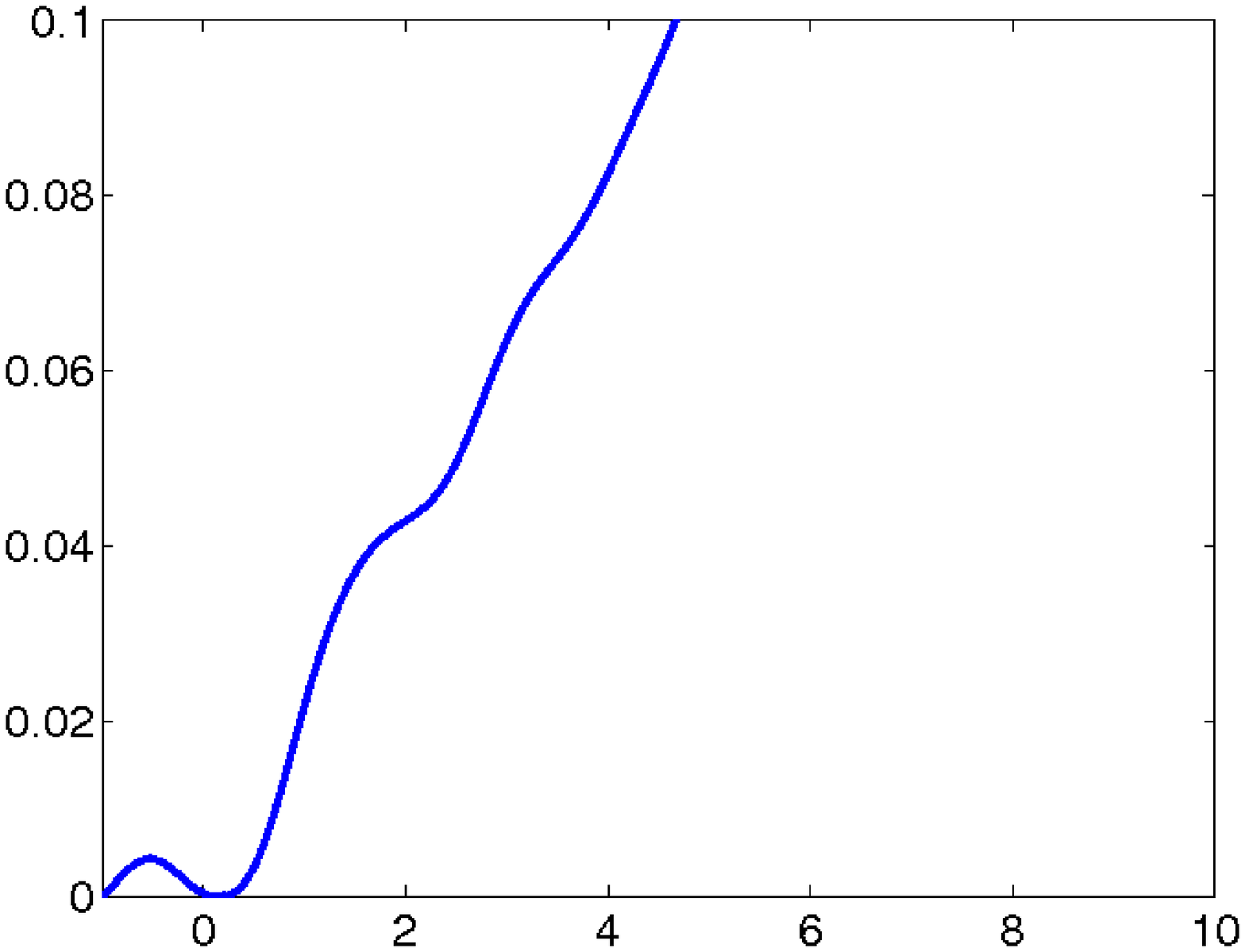}}\quad
\subfloat[$t=0.398$]{\includegraphics[width = 2.1in]{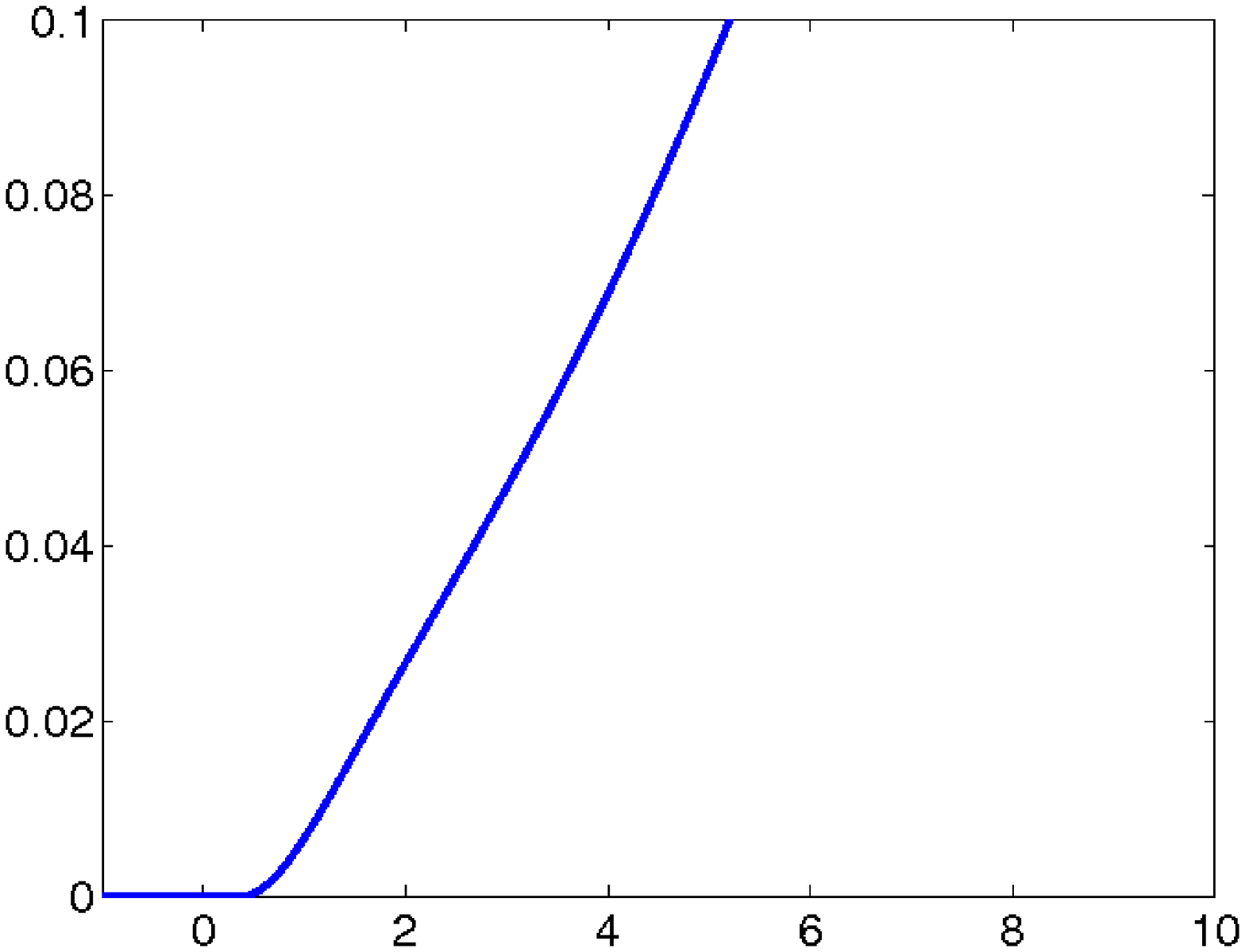}}\\ 
\subfloat[$t=1.594$]{\includegraphics[width = 2.1in]{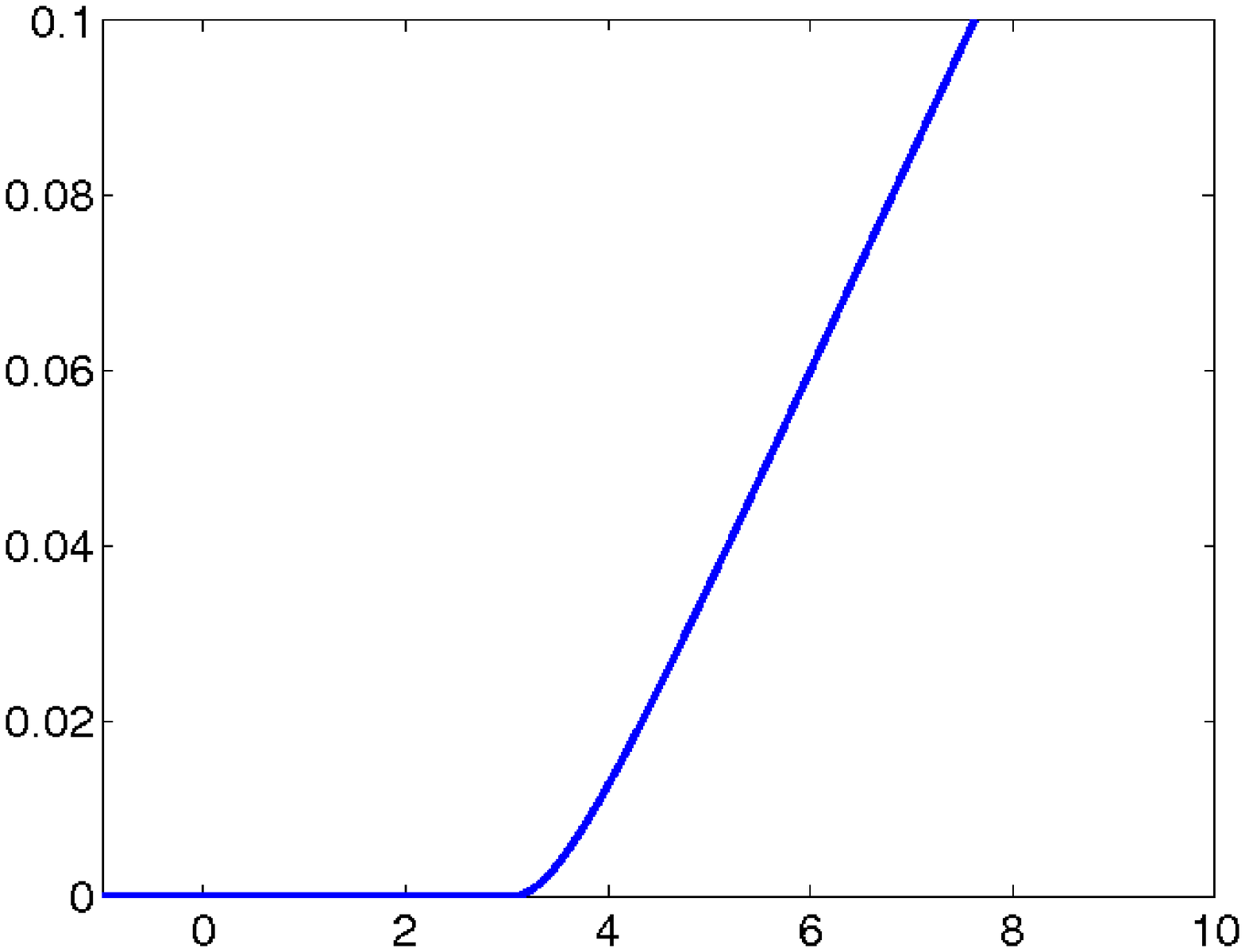}} 
\caption{Numerical solution starting with an initial traveling wave profile perturbed by uniformly random noise sampled from $[0,0.05]$ with $\sigma=2$ and $\gamma=0.05$. This is solved on a grid of $500$ points.}
\label{fig:TWn}
\end{figure}

To test the stability of these traveling wave solutions, we initialize our numerical scheme with the traveling wave profile perturbed by uniformly random noise sampled from $[0,0.05]$. The time evolution is shown in Figure~\ref{fig:TWn}. In a short time, the Laplacian term dominates the evolution, which is expected. The solution gradually smoothes down to a new traveling wave profile and begins to translate at the expected speed. This shows that the traveling wave solution is an attracting solution, at least locally.

\subsection{One Dimensional Heat Equation}
\begin{figure}[ht]
\begin{minipage}[b]{0.45\linewidth}
\centering
{\includegraphics[width=2.2in]{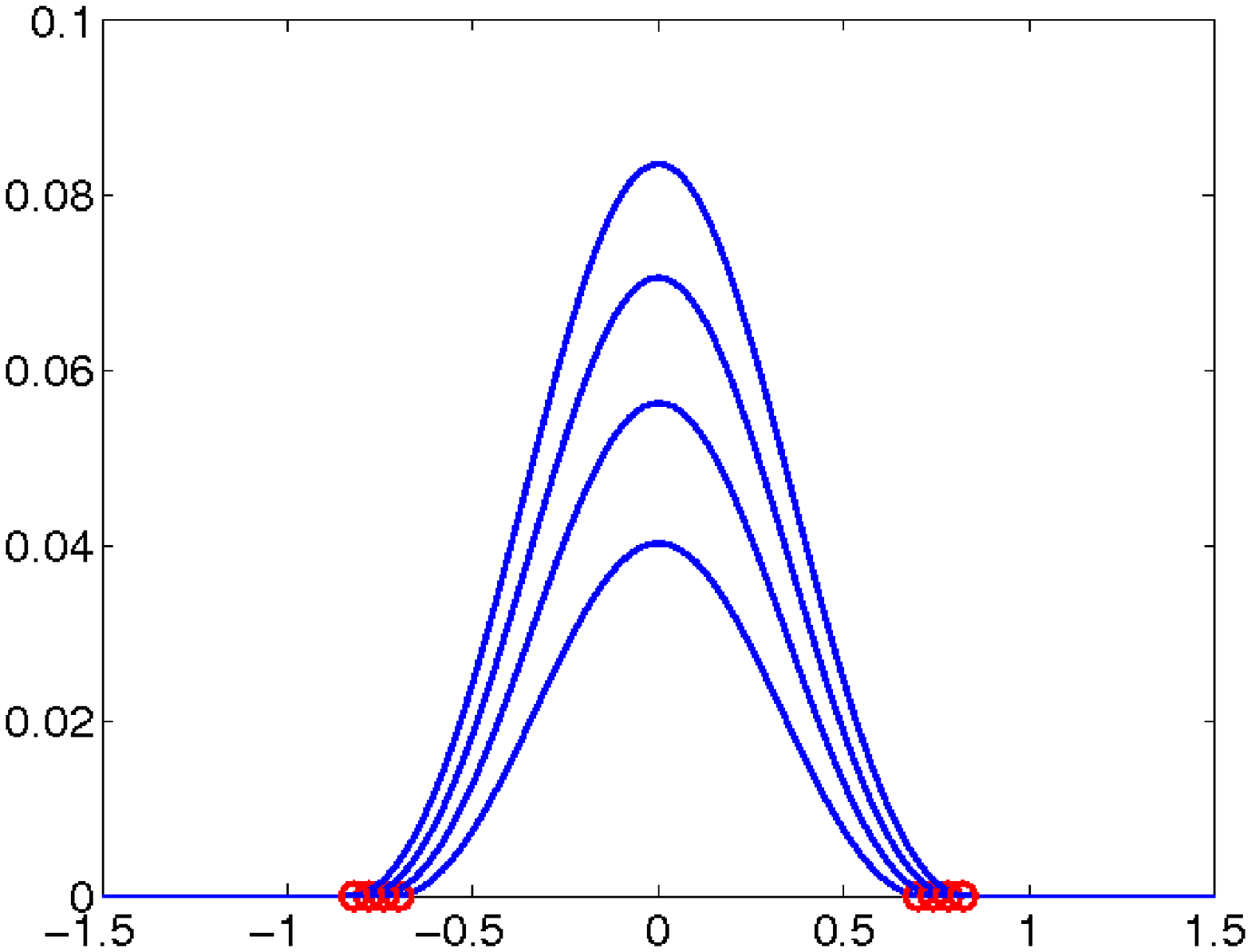}}
\caption*{ }
\end{minipage}
\begin{minipage}[b]{0.45\linewidth}
\centering
 \begin{tabular}{|c|c|} \hline
    {Number of grid points} & Estimate of $a_1$\\
    \hline
    256 & 0.948 \\
    512 & 0.979  \\
    1024 & 0.985 \\
    2048 & 0.991 \\
    4096 & 0.995 \\
    8192 & 0.997 \\
    16384 & 0.997 \\
    \hline  
    \end{tabular}
\caption*{ }
\end{minipage}
\caption{The graph is a 1D simulation of the heat equation with the subgradient term, zero initial data, and a Gaussian forcing function centered around zero, $f(x) = 2 e^{-5x^2}$. The solutions are growing upward in time and their support sets are marked by red circles. The table shows the estimate of the coefficient $a_1$ from Equation~\eqref{a(t)} under grid refinement.}
\label{fig:1dheat}
\end{figure}

In Figure~\ref{fig:1dheat}, the plot shows the modified heat equation (Equation~\eqref{parabolic_pde}) with zero initial data and force $f(x) = 2 e^{-5x^2}$. The solutions evolves upward in time with their support sets marked by red circles. We see that the computed solutions are indeed compactly supported in space, as the theory states.  The corresponding table provides a least squares fit to estimate the coefficient $a_1$ from Equation~\eqref{a(t)} under grid refinement. We see that the coefficient $a_1$ approaches the value 1 quickly within some small approximation error, which is used to verify that our numerical approximation is valid.

\subsection{Two Dimensional Heat Equation}
\begin{figure}[t]
\centering
\subfloat[$t=0$]{\includegraphics[width = 2in]{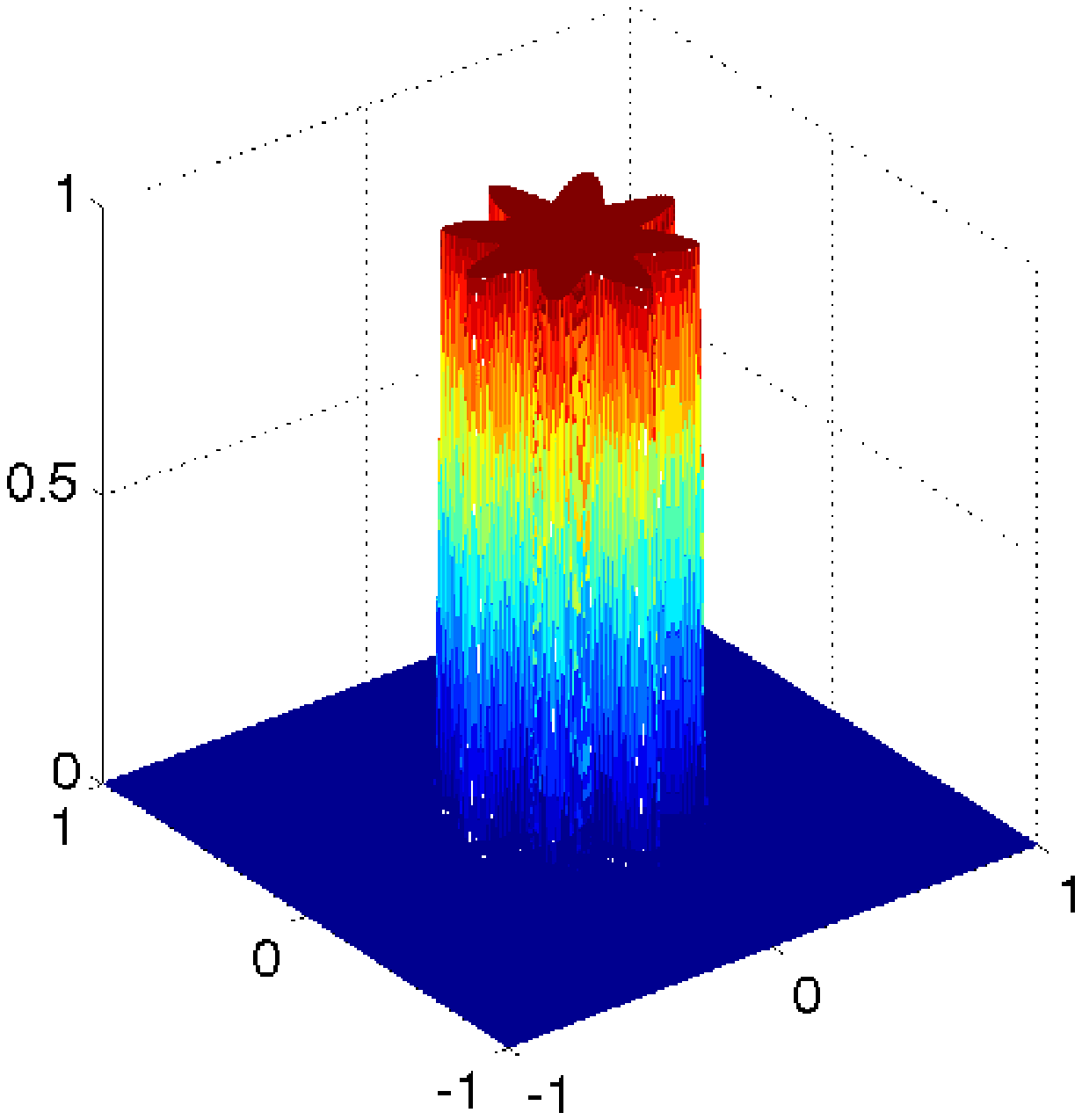}} 
\subfloat[$t=3.2\times 10^{-4}$]{\includegraphics[width = 2in]{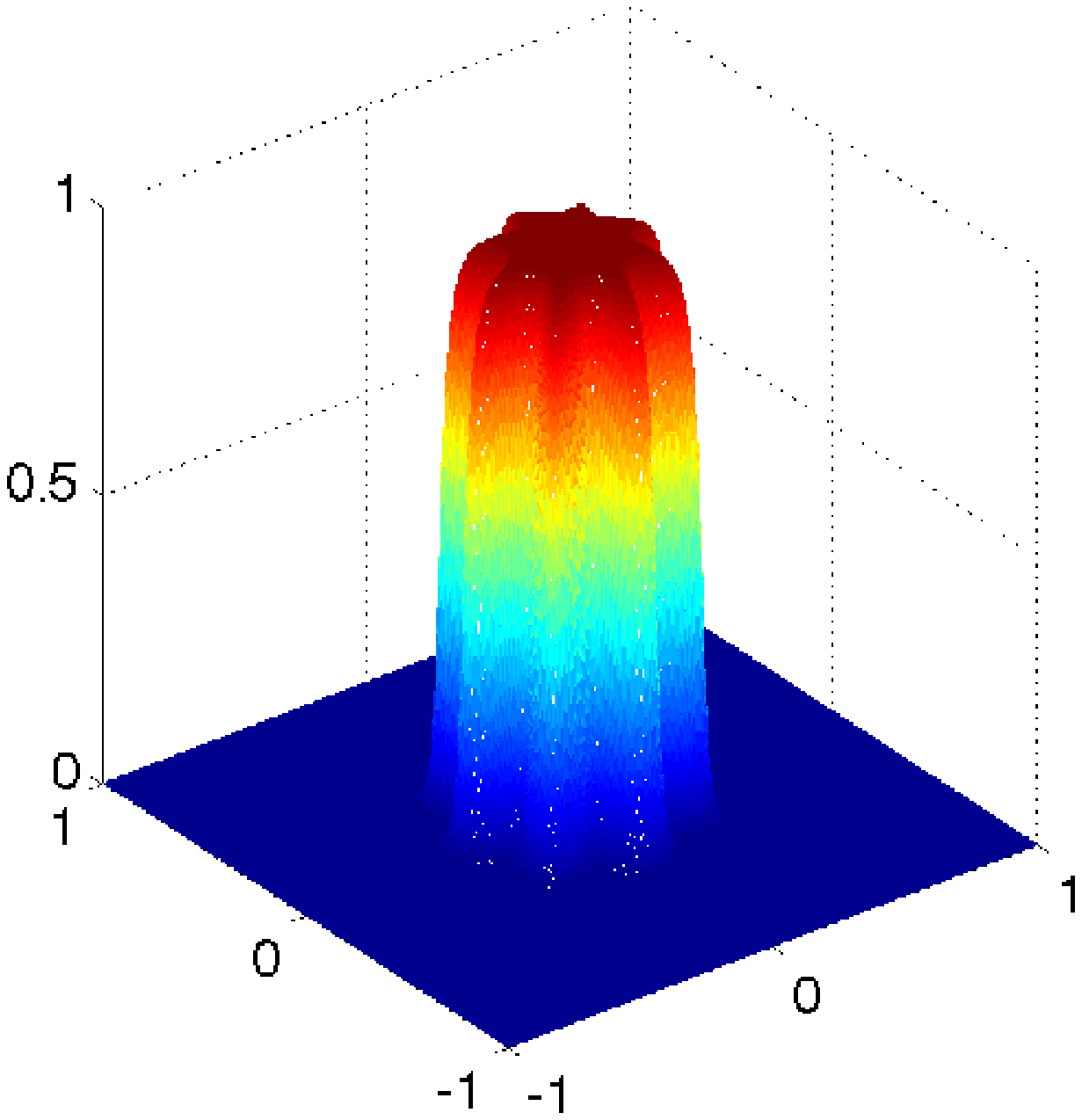}}\\
\subfloat[$t=3.2\times 10^{-3}$]{\includegraphics[width = 2in]{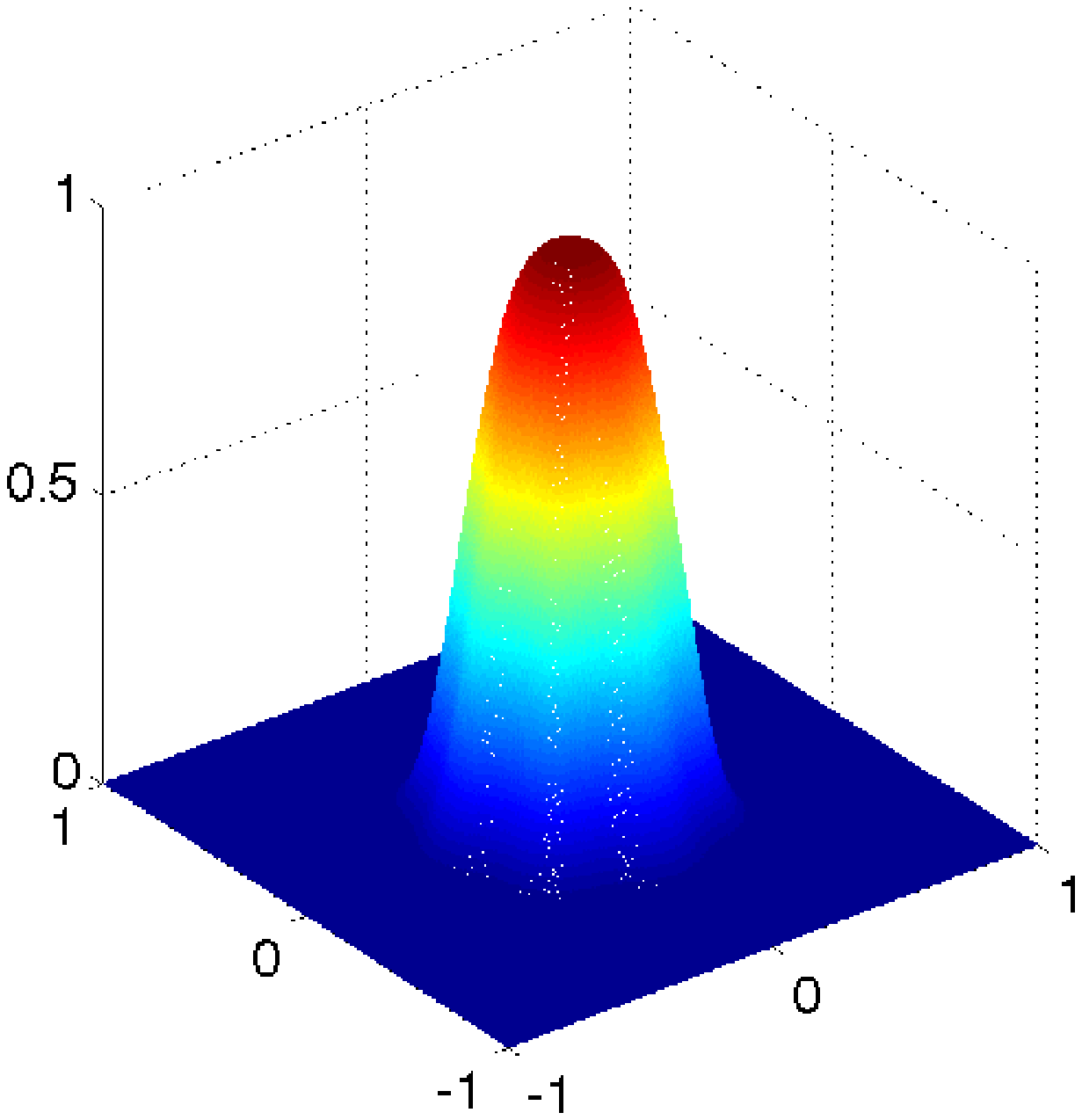}}
\subfloat[$t=1.76\times 10^{-2}$]{\includegraphics[width = 2in]{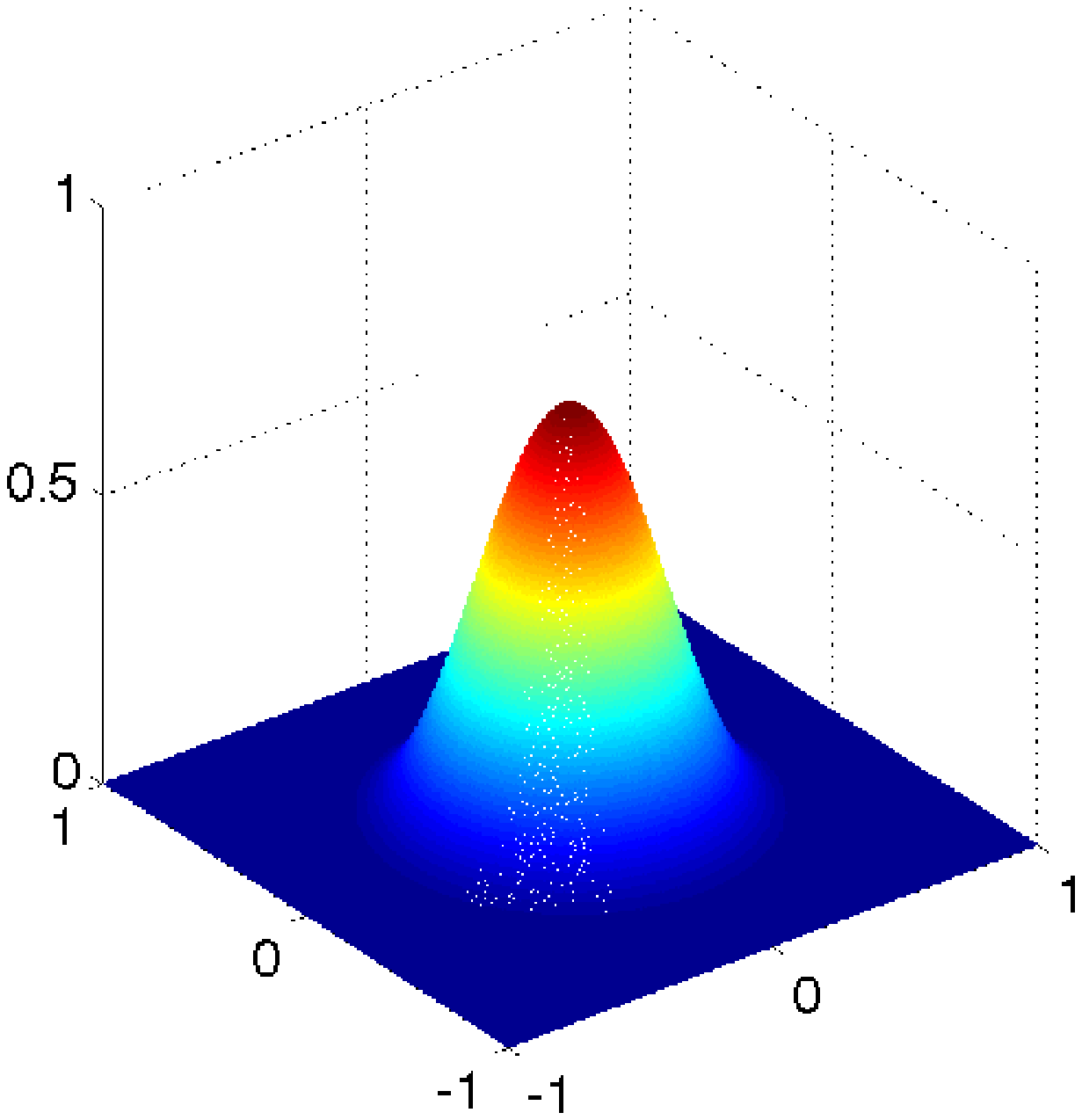}} \\
\subfloat[$t=4.8\times 10^{-2}$]{\includegraphics[width = 2in]{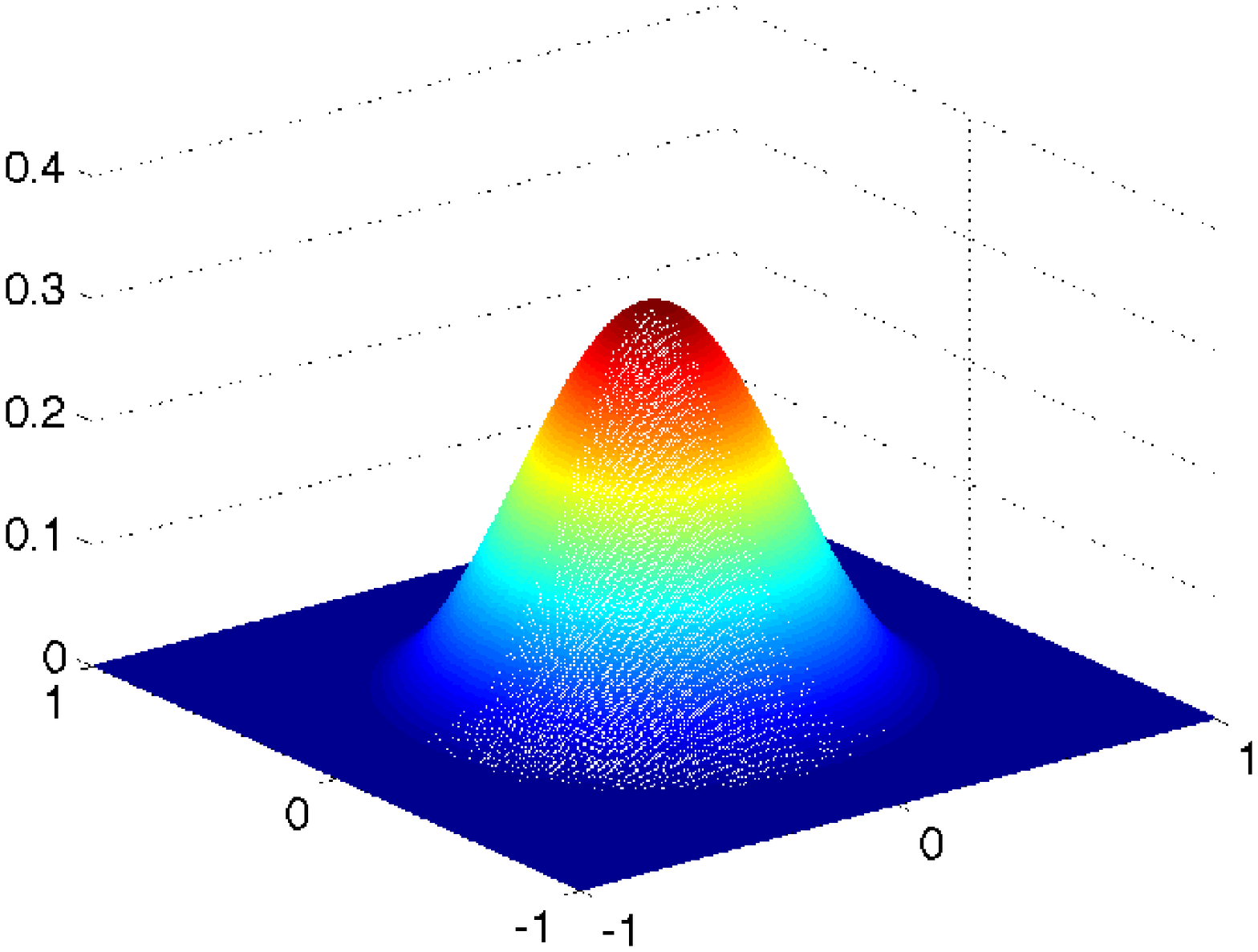}}
\subfloat[$t=0.112$]{\includegraphics[width = 2in]{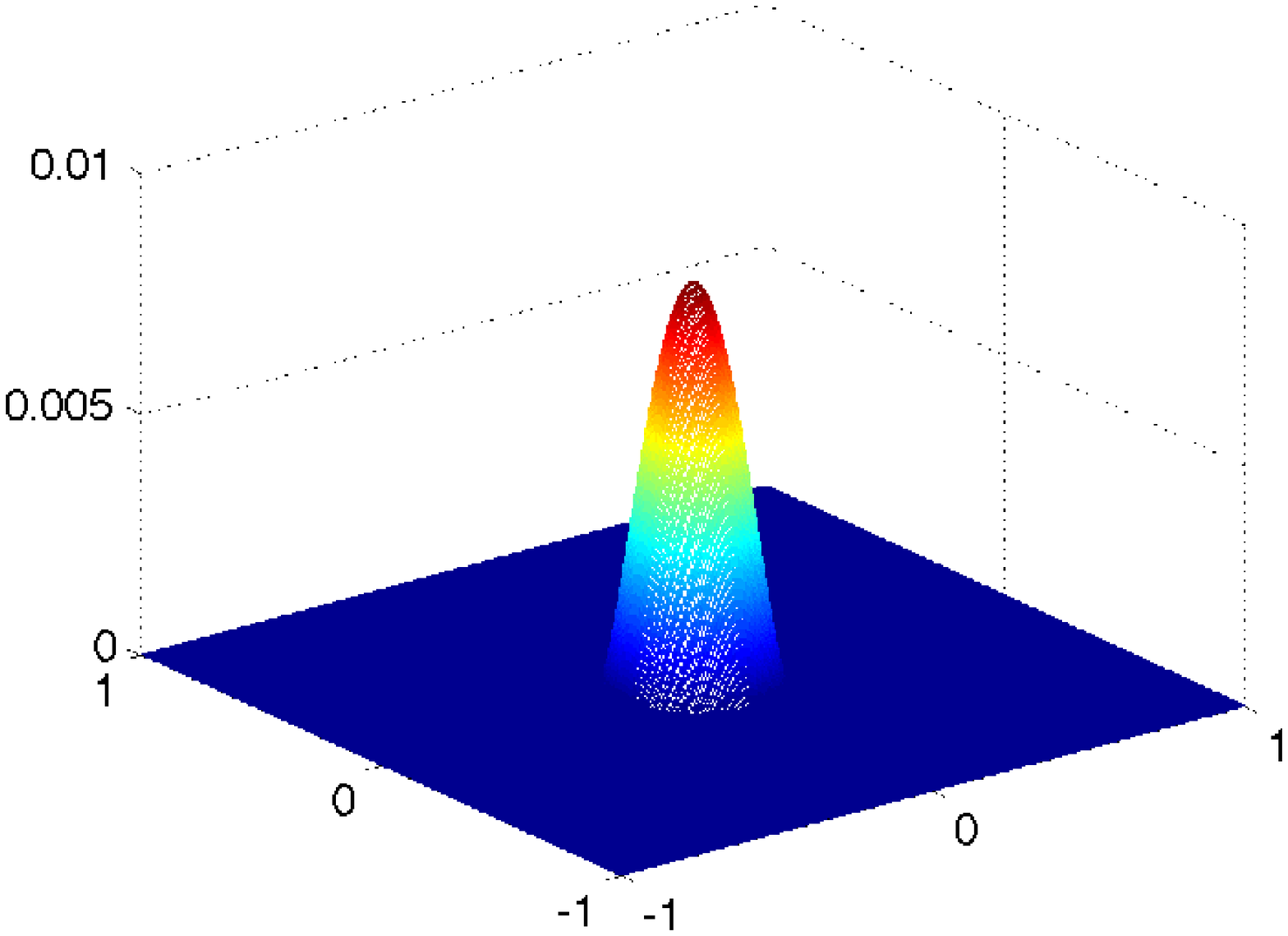}}

\caption{Solutions of the initial value problem (no forcing term) computed on a 500 by 500 grid with $\gamma=2$ at times indicated. The solution smoothes out and decays to zero.}
\label{fig:2dsoln}
\end{figure}

\begin{figure}[h]
\centering
\subfloat[$t=0$]{\includegraphics[width = 2in]{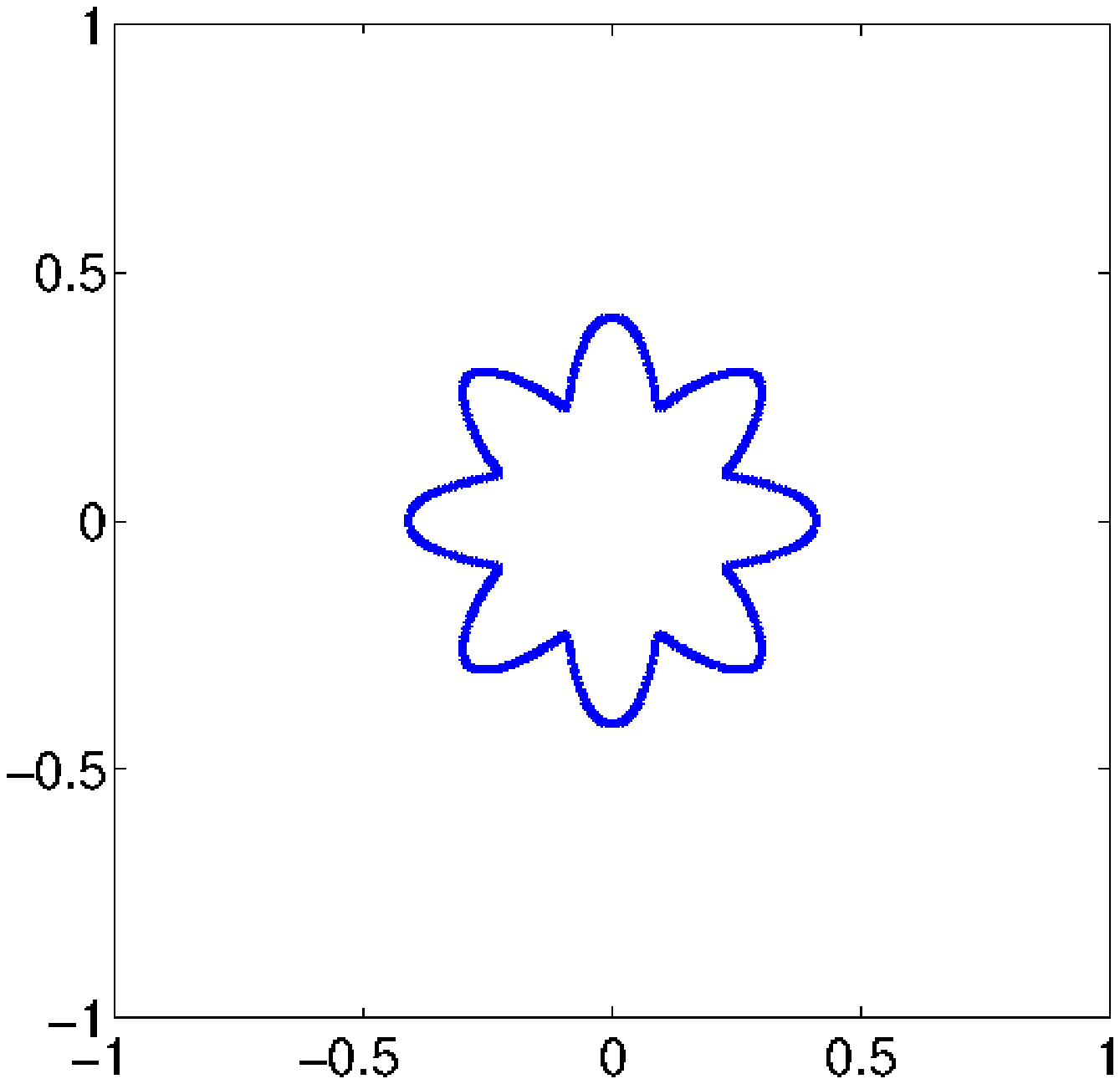}} 
\subfloat[$t=3.2\times 10^{-4}$]{\includegraphics[width = 2in]{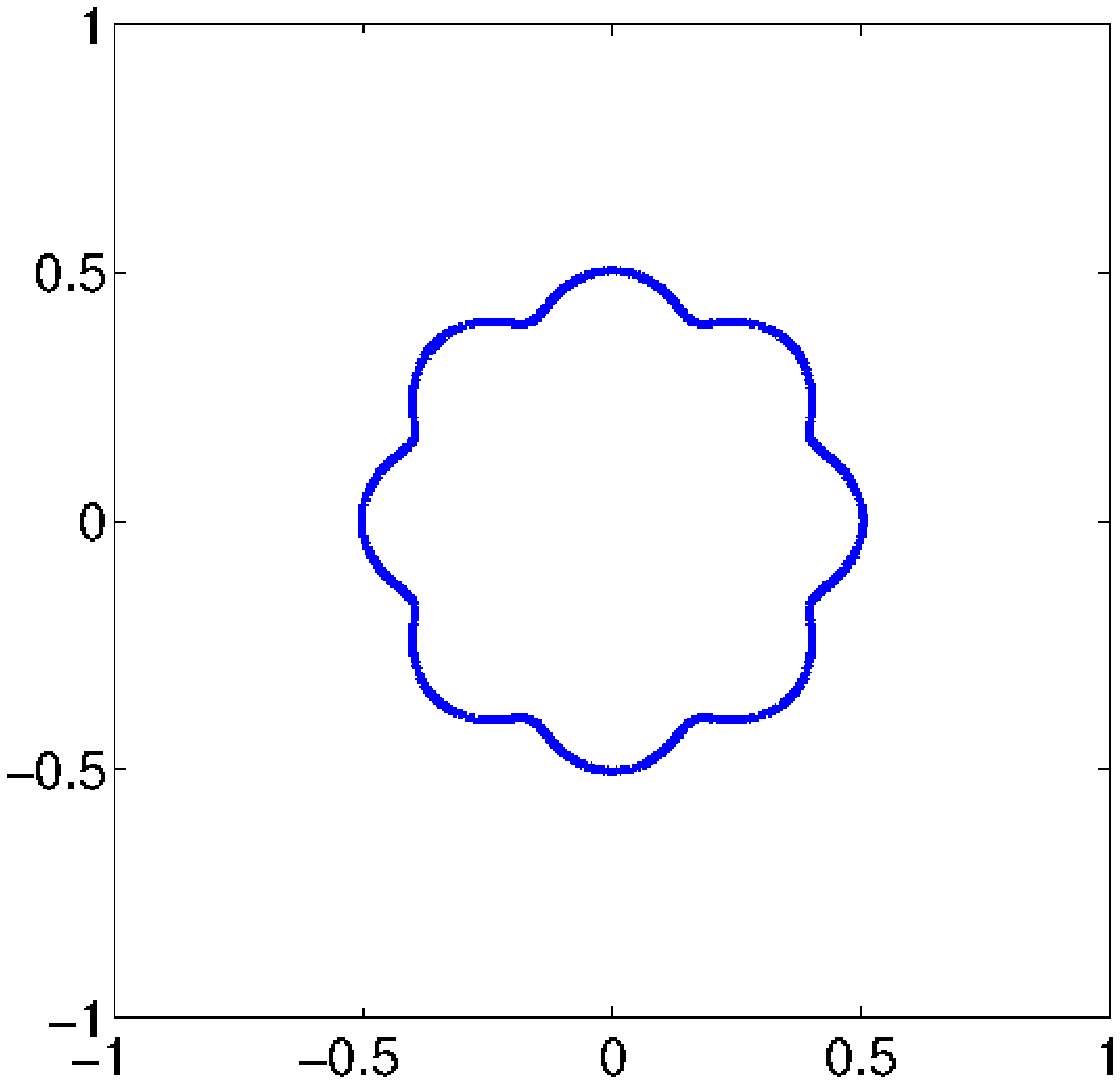}}\\
\subfloat[$t=3.2\times 10^{-3}$]{\includegraphics[width = 2in]{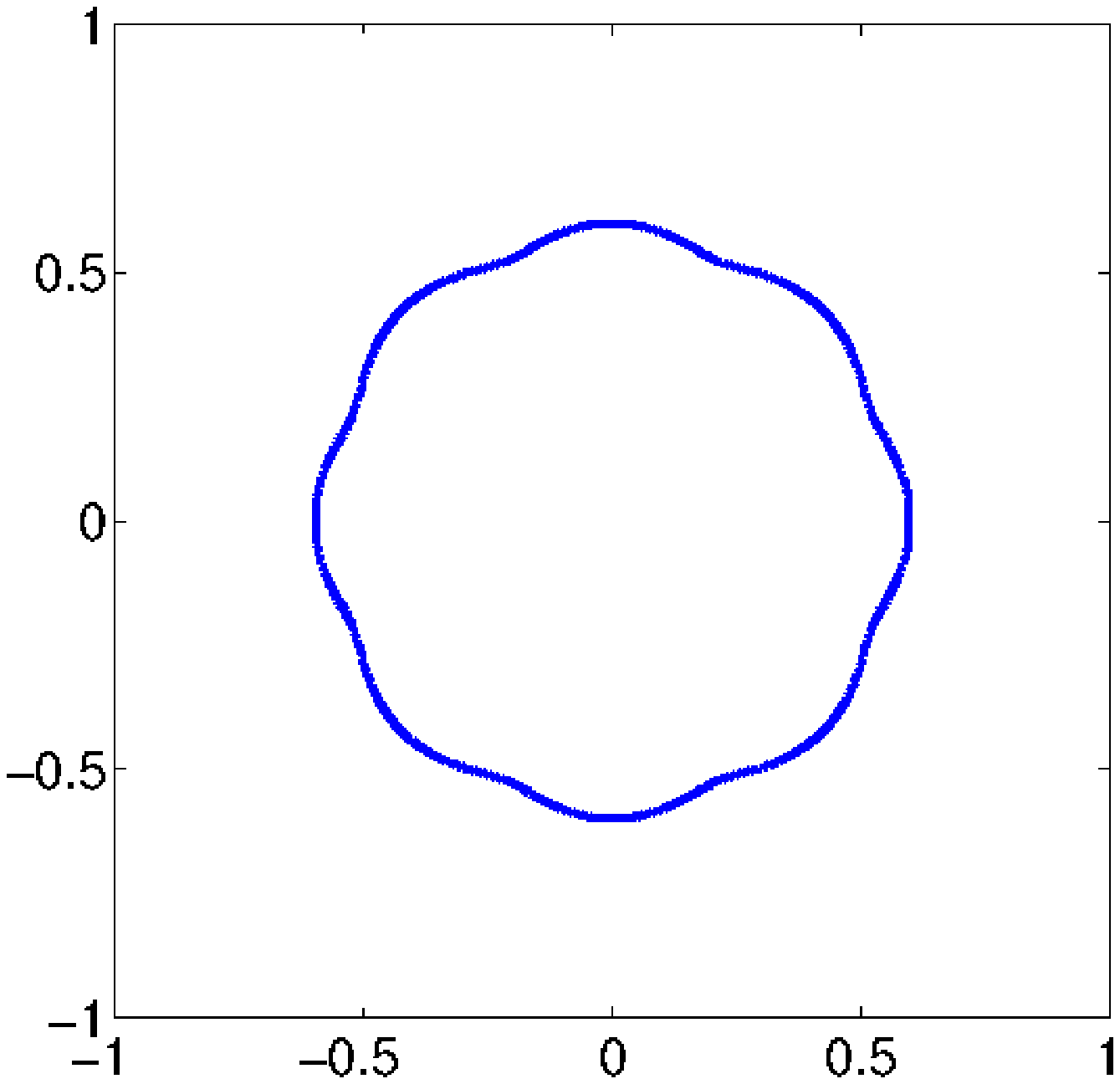}}
\subfloat[$t=1.76\times 10^{-2}$]{\includegraphics[width = 2in]{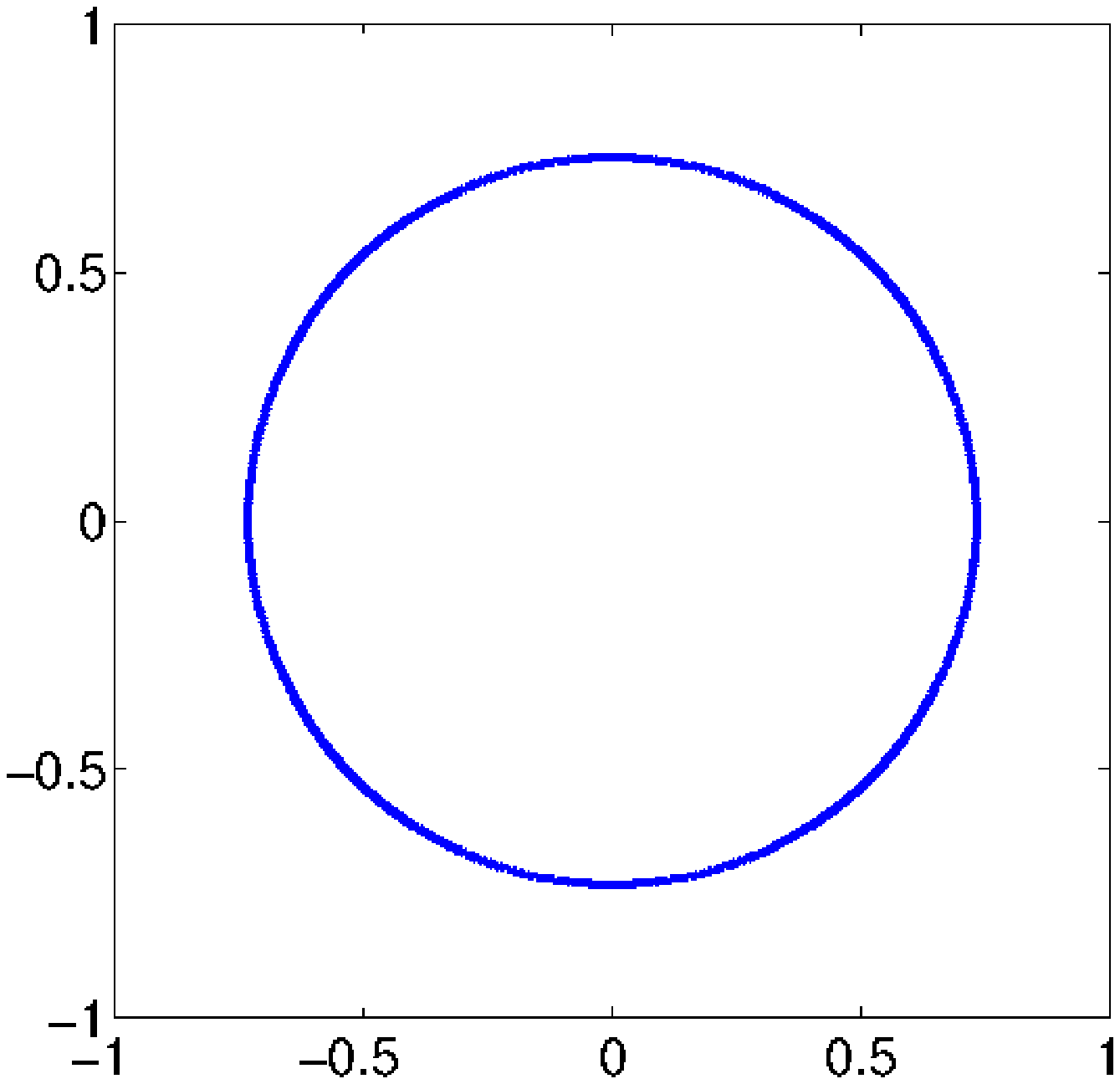}} \\
\subfloat[$t=4.8\times 10^{-2}$]{\includegraphics[width = 2in]{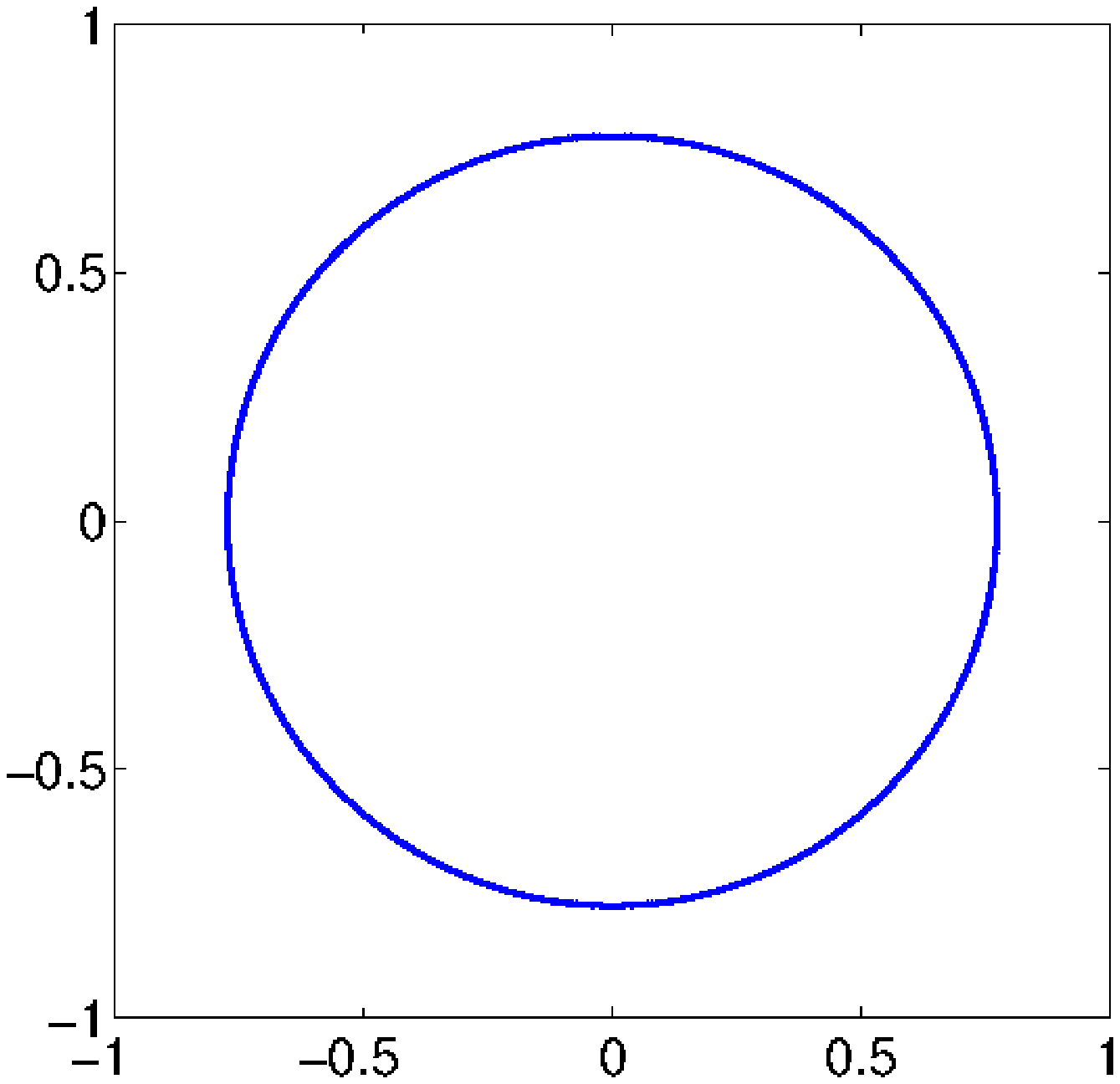}} 
\subfloat[$t=0.112$]{\includegraphics[width = 2in]{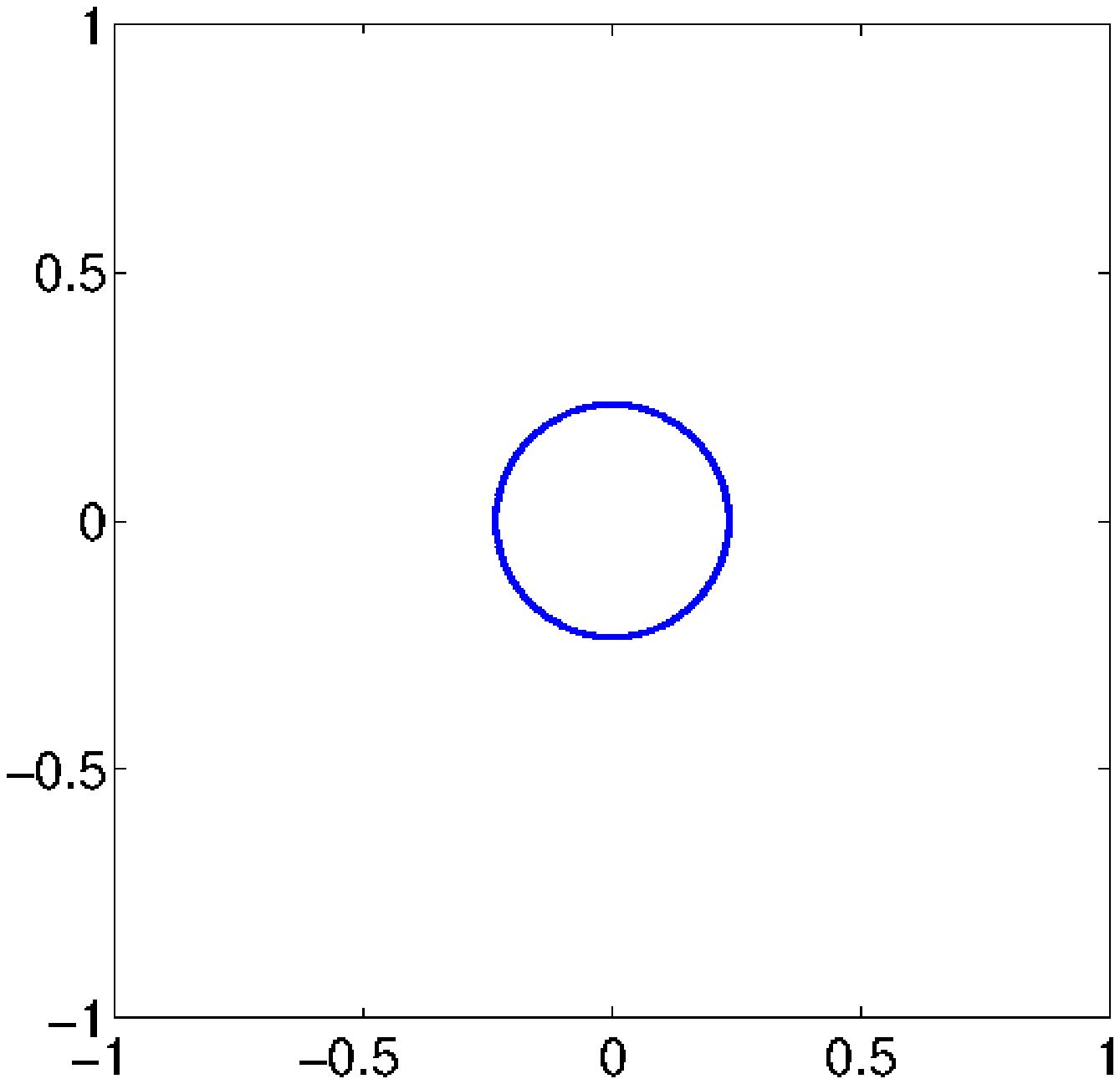}} 

\caption{Support set of the initial value problem in Figure~\ref{fig:2dsoln}. The support set grows outward to a maximum size and retracts inward as the solution decays to zero.}
\label{fig:2dsupp}
\end{figure}

In Figure~\ref{fig:2dsoln}, we compute the solution of  Equation~\eqref{parabolic_pde} with $\gamma=2$ and $f=0$. In this case, we apply the parabolic Douglas-Rachford algorithm, which allows for larger time-steps. The initial data is a smoothed indicator function on the star shaped domain. In Figure~\ref{fig:2dsupp}, the corresponding support set of Figure~\ref{fig:2dsoln} is shown. The support set grows outward to a maximum size and retracts inward as the solution decays to zero. The solution is identically zero at time $t=0.1152$.

\subsection{Graph Diffusion}
In higher dimensions, we can consider the standard normalized diffusion equation:
\begin{equation}
\begin{aligned}
u_t = L_g u& :=  -\left( \text{I} - D^{-1/2} A D^{-1/2}  \right) u\\
u(x,0)&=g(x),
\end{aligned}\label{graphLap}
\end{equation}
where $L_g$ is the graph Laplacian,  $A $ is the adjacency matrix,  and $D$ is the degree matrix. For more on the graph Laplacian, see \cite{GraphBook,SpectralClustering}.

\begin{figure}
\centering
\subfloat[$t=0$]{\includegraphics[width = 2in]{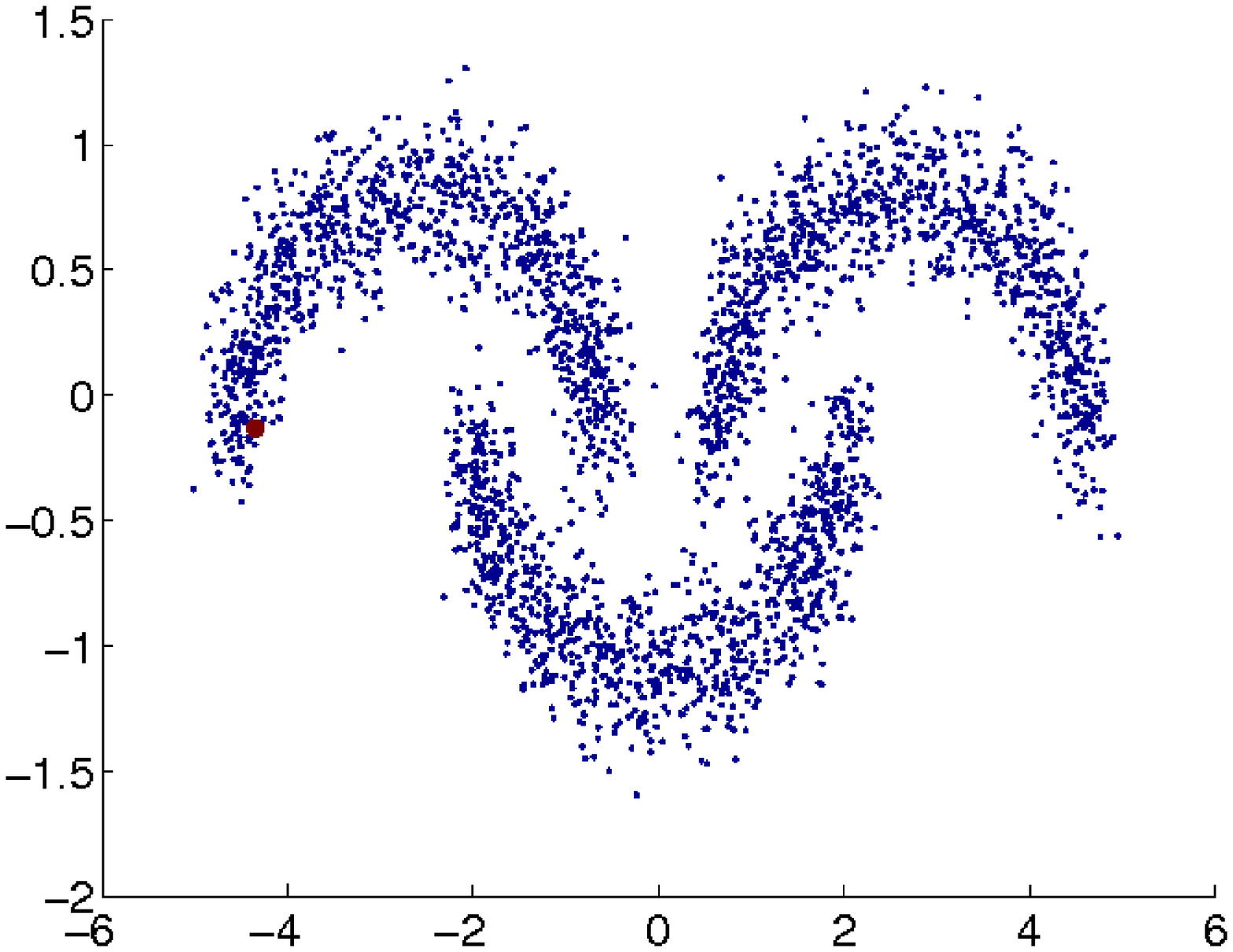}} 
\subfloat[$t=47.5$]{\includegraphics[width = 2in]{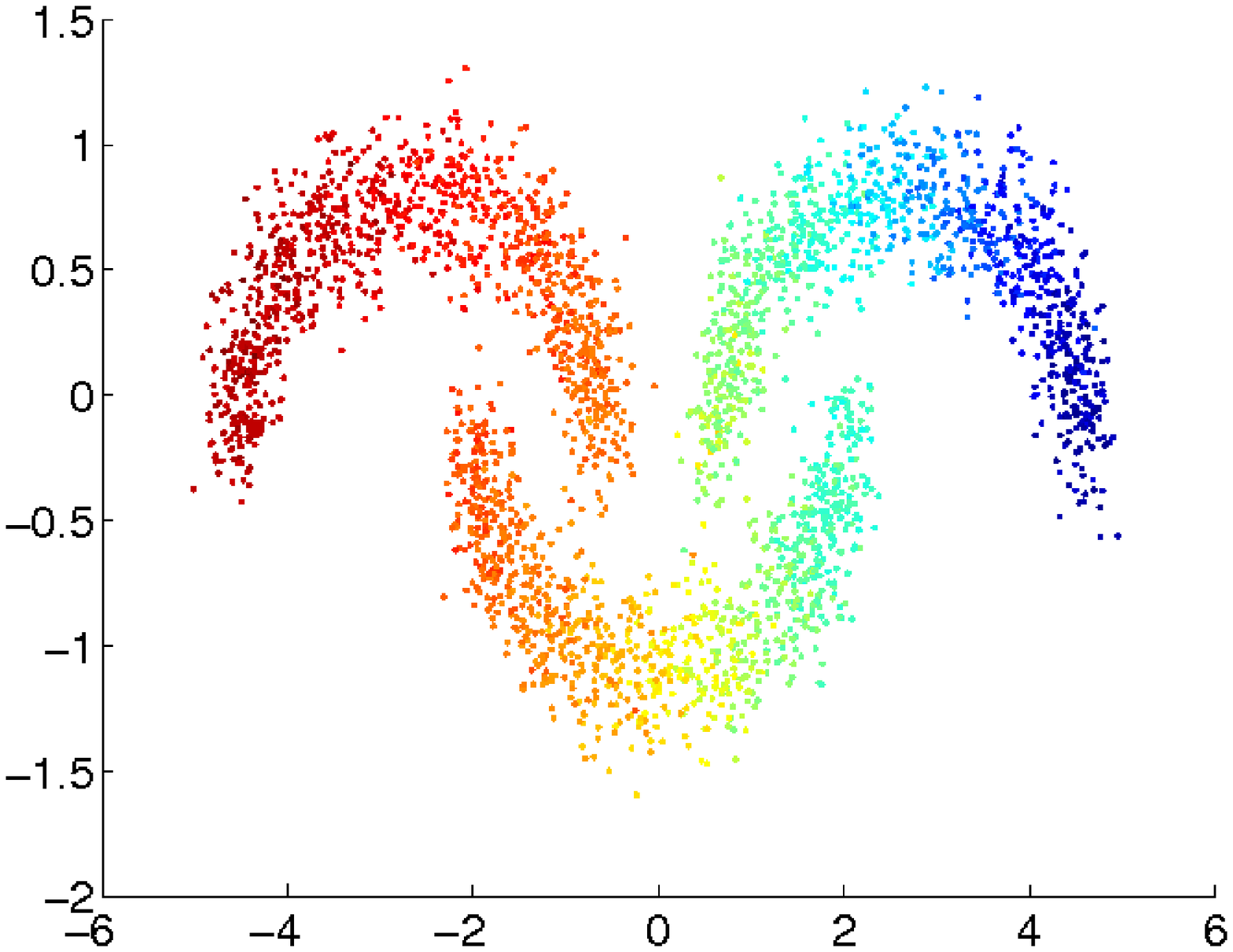}}\\
\subfloat[$t=475$]{\includegraphics[width = 2in]{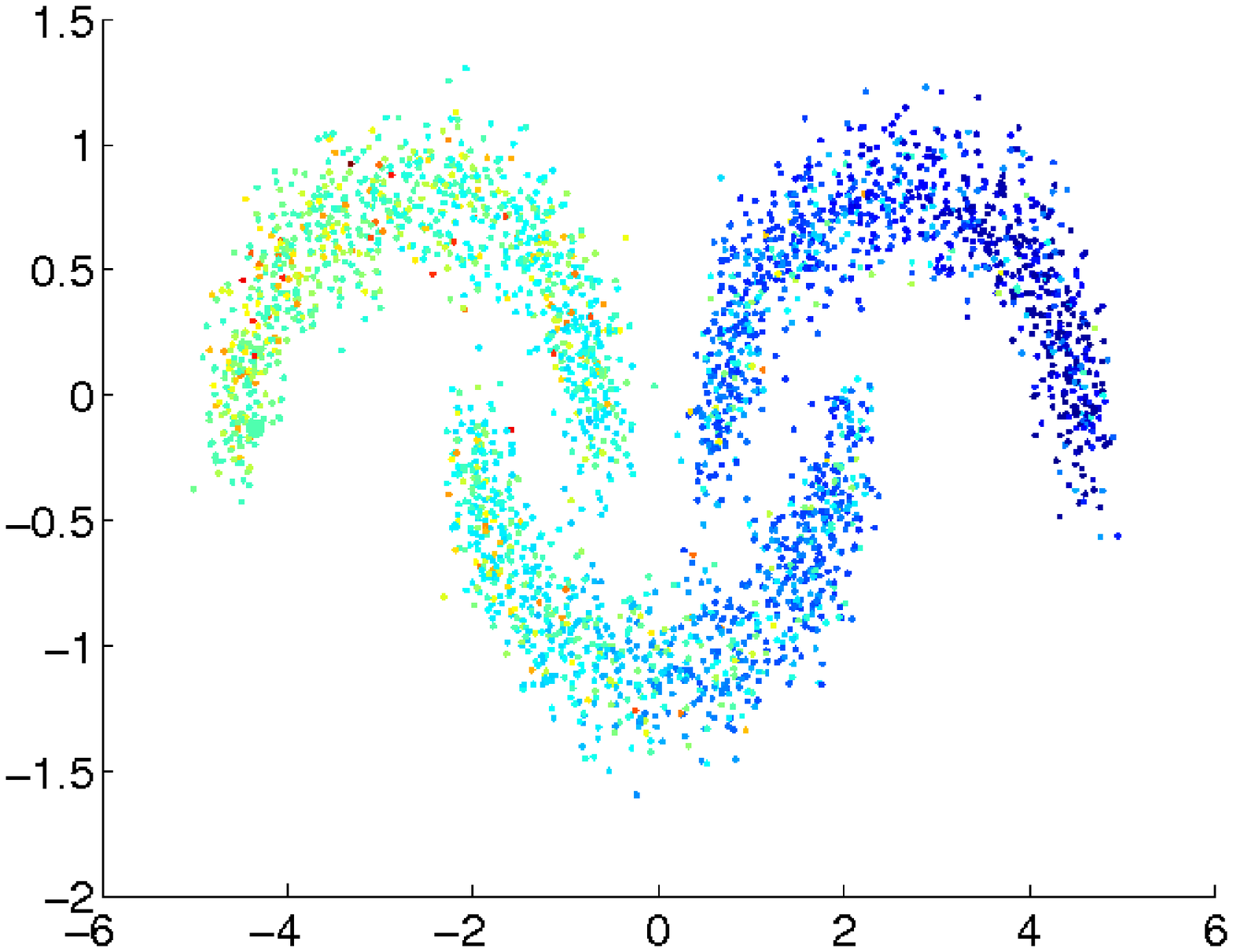}}
\subfloat[$t=1425$]{\includegraphics[width = 2in]{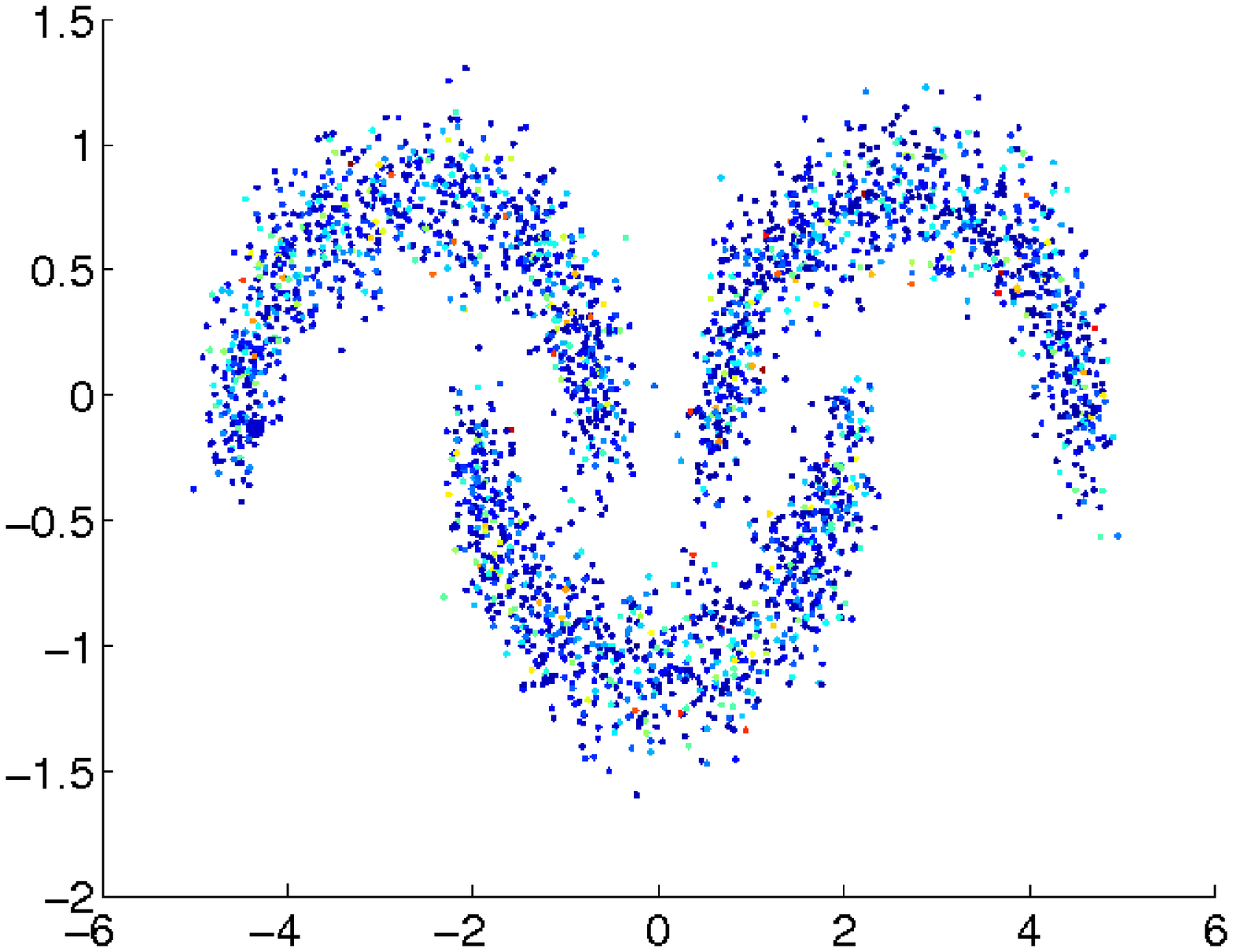}} 
\caption{Solution of the initial value problem diffusing standard normalized graph Laplacian. }
\label{fig:Graph1}
\end{figure}

In Figure~\ref{fig:Graph1}, the points represent the projection of vectors from $\mathbb{R}^{100}$ and each point is connected to many others in a non-local fashion. For the initial data, we concentrate the mass on one point in the far left, specifically, the $u(x_j,0)=\delta_{j,1000}$ where $\delta_{j,k}$ is the Kronecker delta function. As the system evolves governed by Equation~\eqref{graphLap}, the solution becomes strictly positive quickly. 

The modified equation is:
\begin{equation}
\begin{aligned}
&u_t = -\left( \text{I} - D^{-1/2} A  D^{-1/2}  \right) u - \gamma p(u),\\
&u(x,0)=g(x).
\end{aligned}\label{graphLap2}
\end{equation}

\begin{figure}
\centering
\subfloat[$t=0$]{\includegraphics[width = 2in]{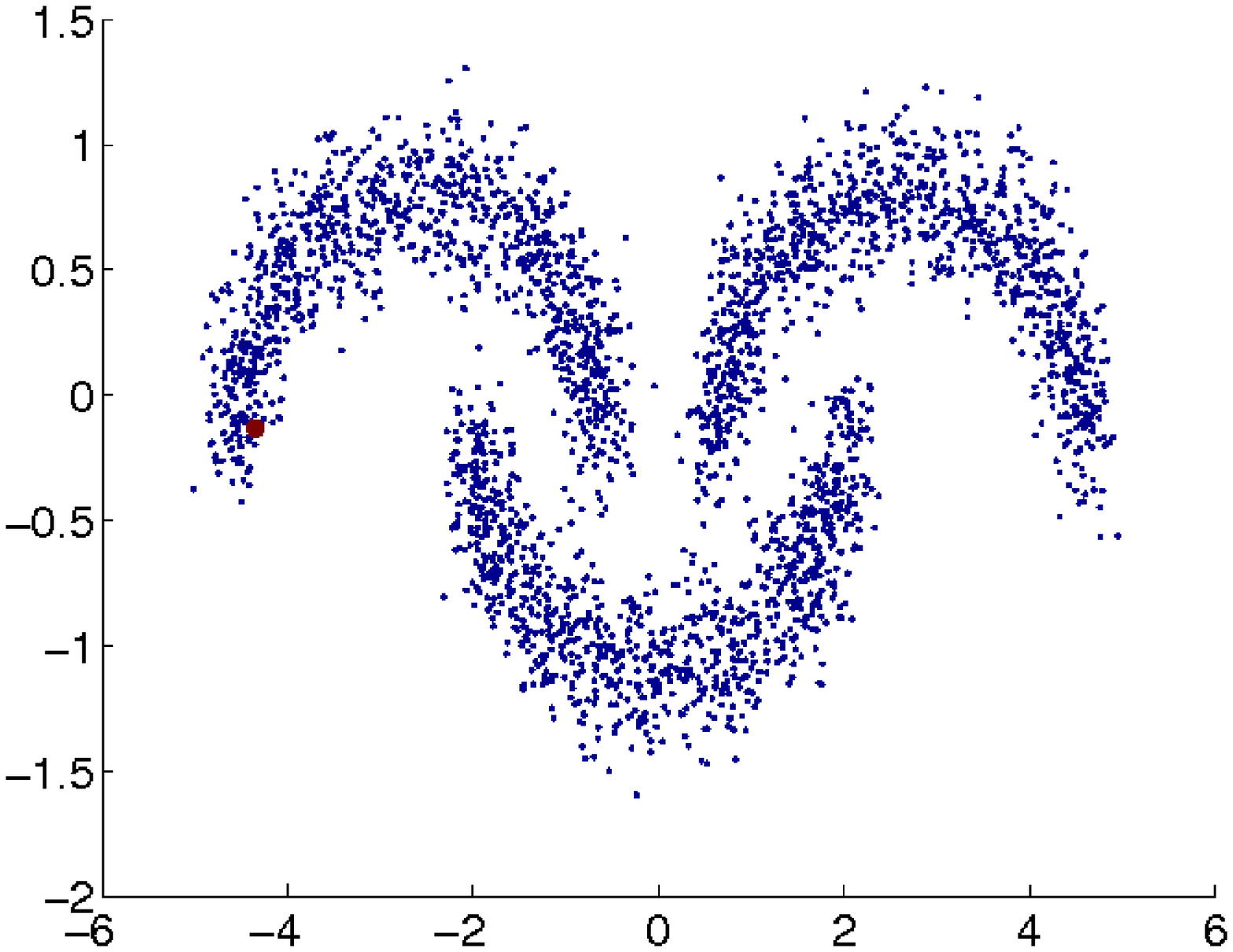}} 
\subfloat[$t=2.85$]{\includegraphics[width = 2in]{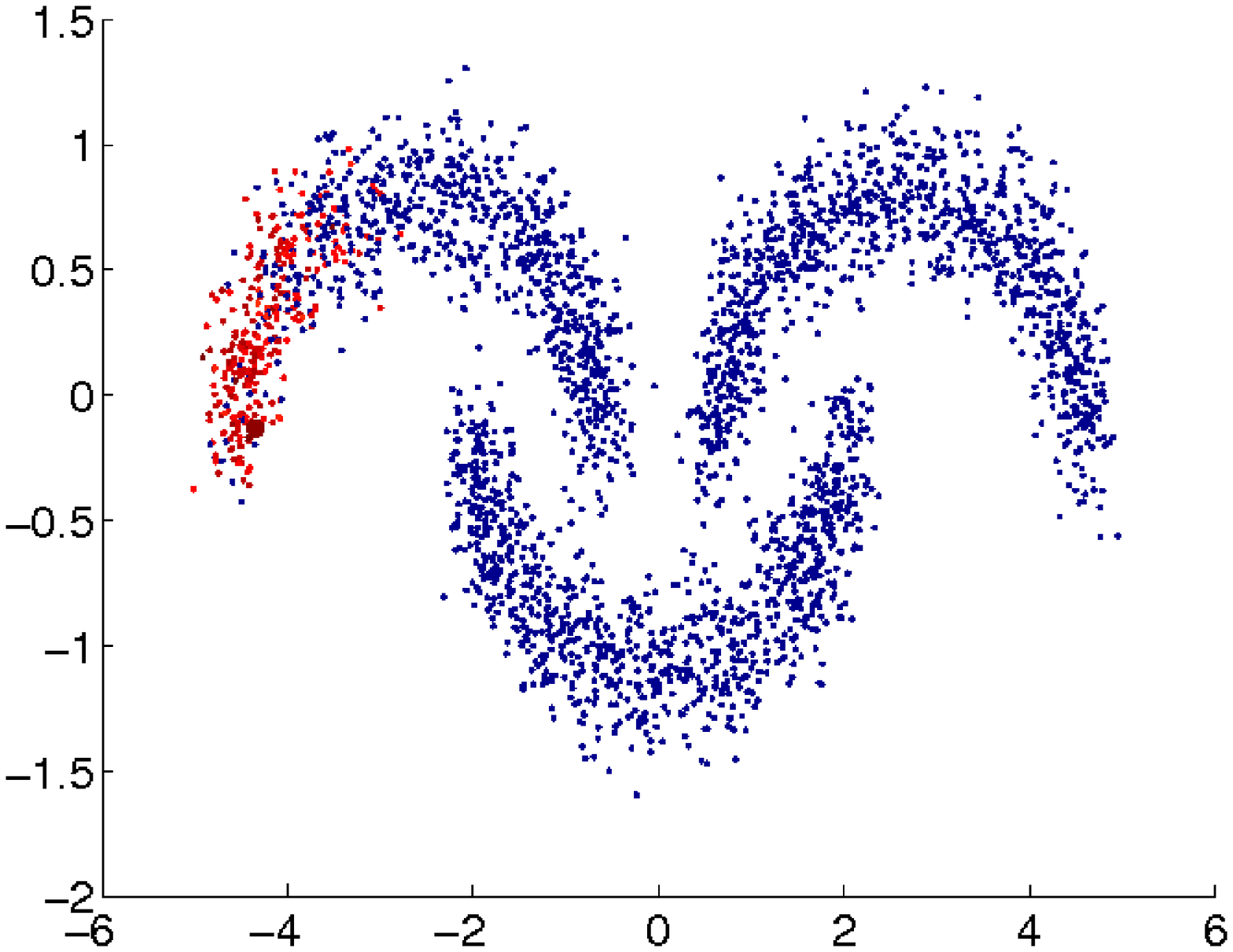}}\\
\subfloat[$t=6.65$]{\includegraphics[width = 2in]{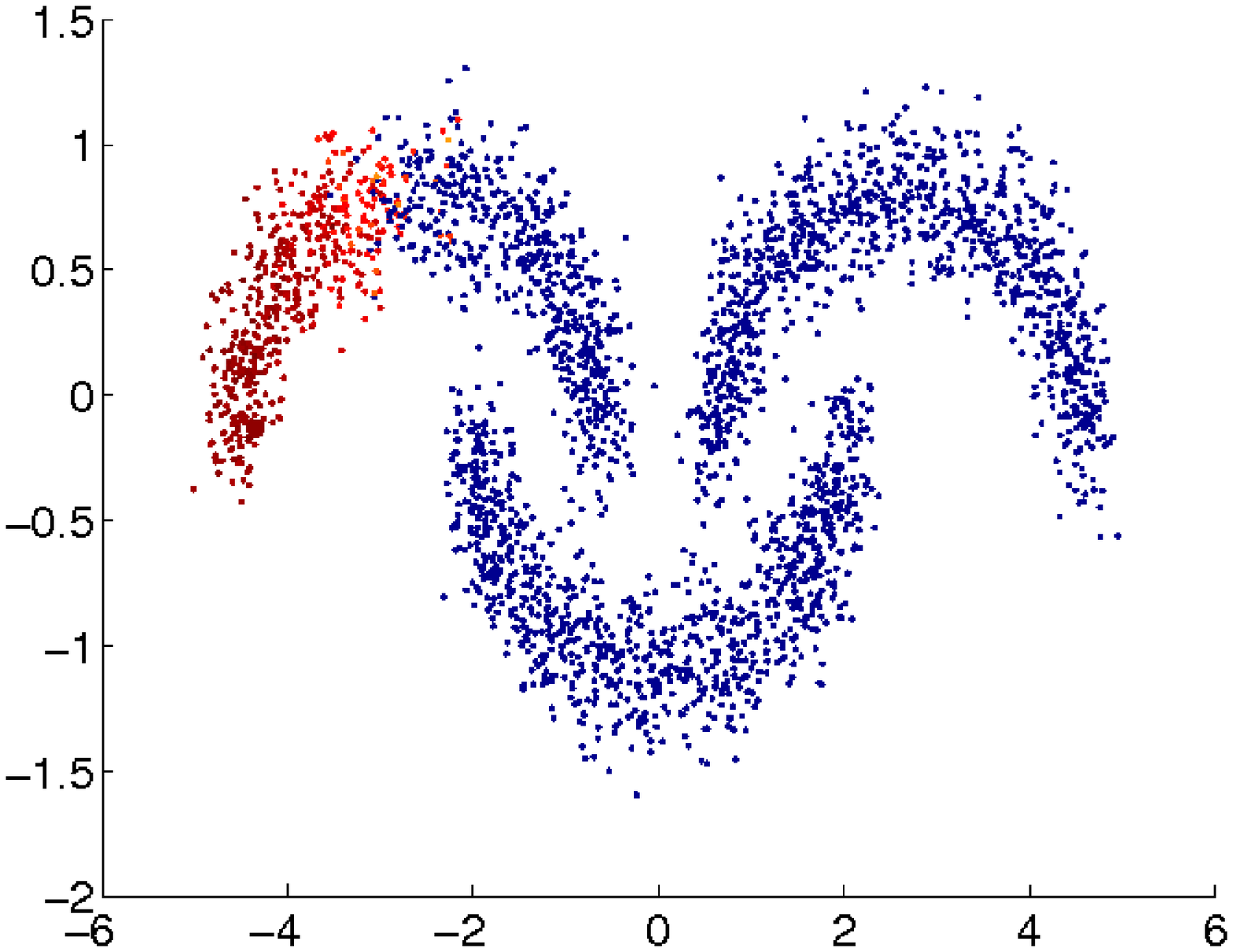}}
\subfloat[$t=11.4$]{\includegraphics[width = 2in]{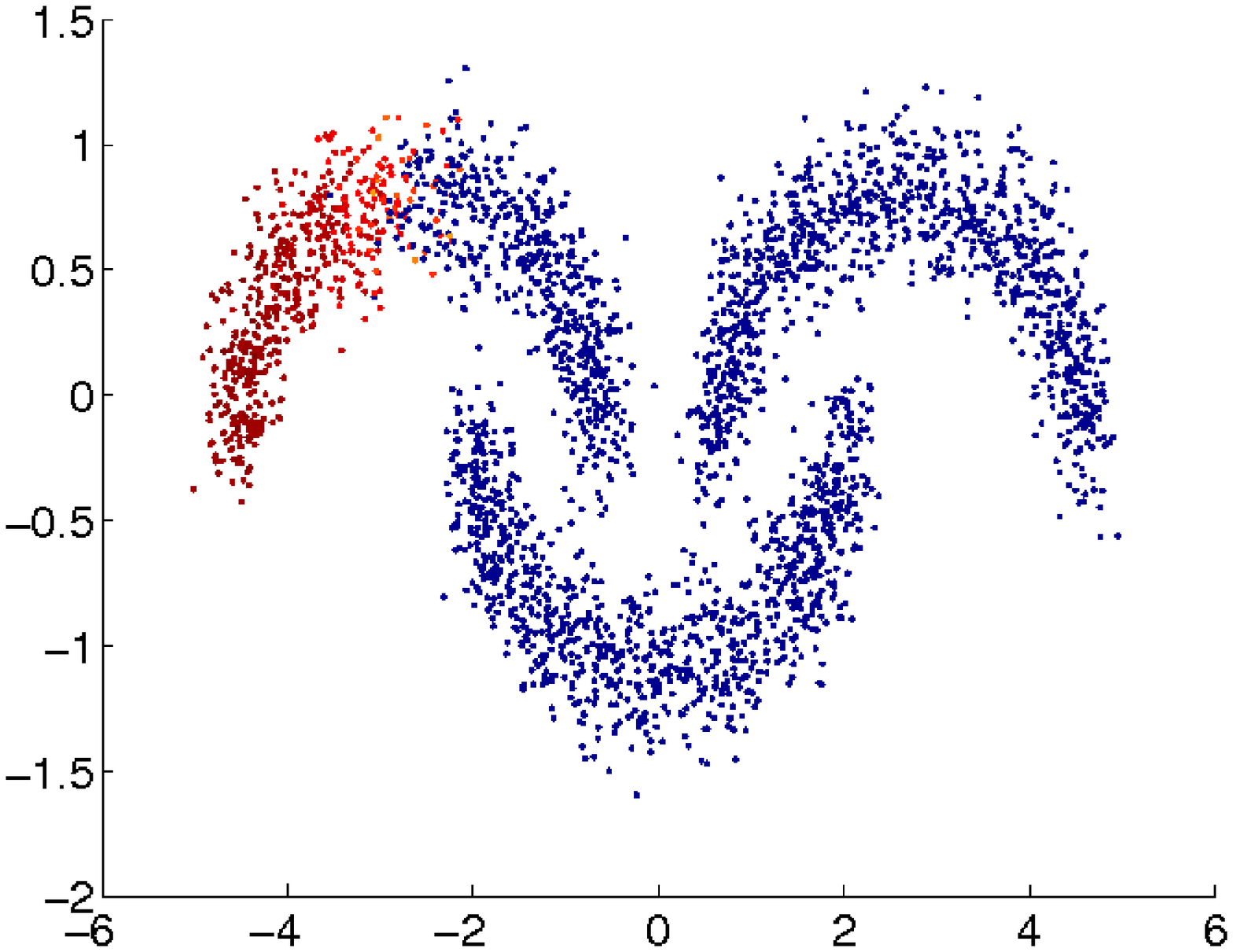}} 
\caption{Solution of the initial value problem with the subgradient term, $\gamma=5\times 10^{-5}$.}
\label{fig:Graph2}
\end{figure}

In Figure~\ref{fig:Graph2}, we begin with the same initial condition and see that over time the support set does not grow past a bounded region if $u$ evolves as in \eqref{graphLap2}. Therefore, numerically we show that the support is of finite size for the case of graph diffusion. In Figure~\ref{fig:Graph2}(d), the solution begins to decay to zero which causes its support set to retract towards the initial support before vanishing.

\subsection{Signum-Gordon Equation}

\begin{figure}
\includegraphics[width = 2 in]{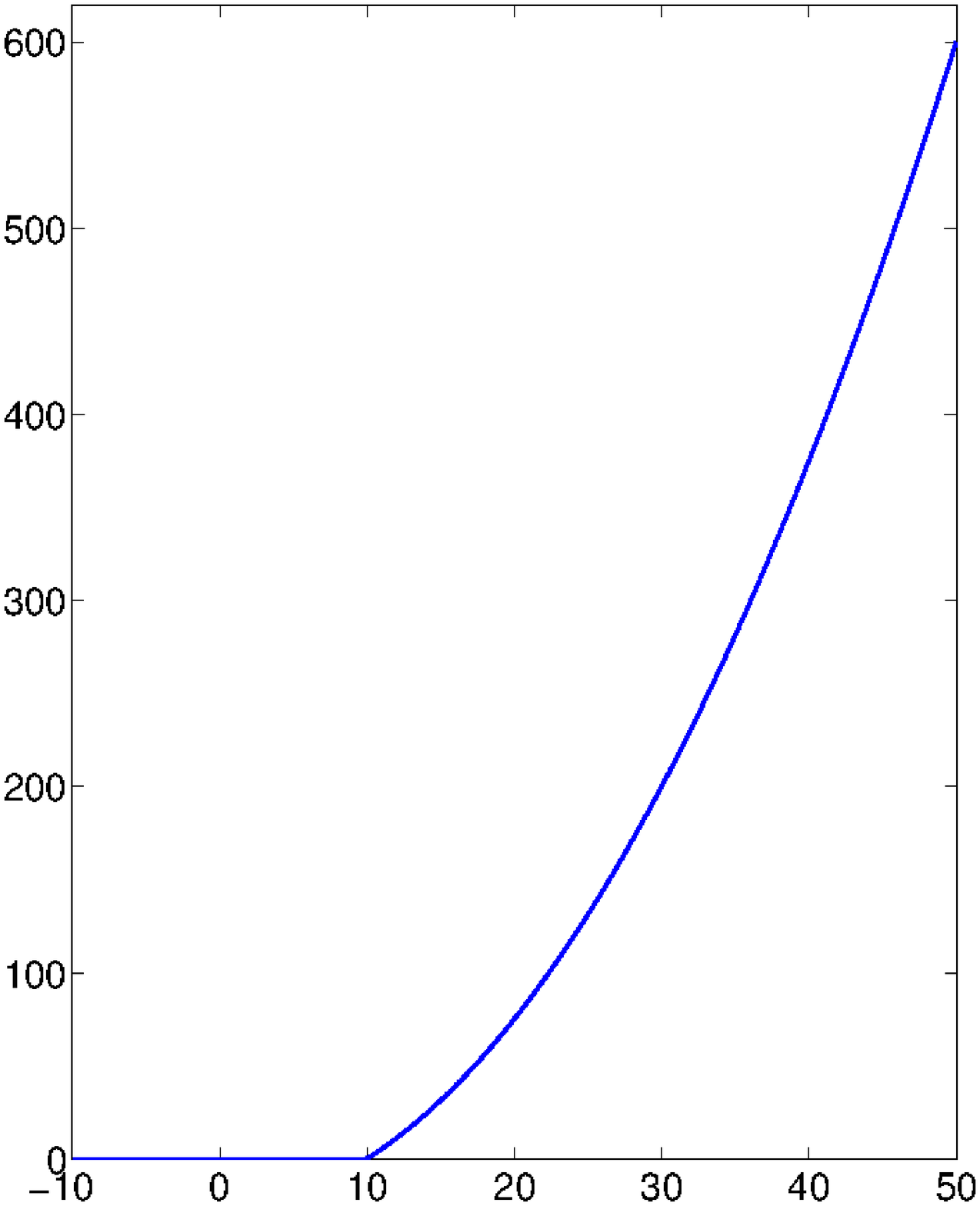}
\caption{Our numerical approximation to a compact traveling wave solution to the Signum-Gordon equation.}
\label{tw:sg}
\end{figure}

\begin{table}[t!]
\begin{center}
\rowcolors{1}{lightgray}{lightgray}
 \begin{tabular}{|c|c|c|c|c|c|c|c|}\hline
 Grid Size & 128  & 256  & 512  & 1024  & 2048 & 4096 & 8192\\ \hline
$L^2$-Error  & 0.4601 & 0.2319 & 0.1133 & 0.0569 & 0.0284 & 0.0143 & 0.0072\\
\hline

\end{tabular}
\label{table:SG}
\caption{Error between our numerical solution and the analytic solution of the Signum-Gordon Equation.}
\end{center}

\end{table}

The signum-Gordon equation has an interpretation as an approximation to certain physical models \cite{arodz2008compact, arodzWave, arodz2011swaying}. The equation takes the form of a second order nonlinear hyperbolic equation:
\begin{equation}
\begin{aligned}
u_{tt} - \Delta u &= -\text{sign}(u)\\
u(x,0)&=g_1(x)\\
u_t(x,0)&=g_2(x),
\label{signum_gordon}
\end{aligned}
\end{equation}
and exhibits both compactly supported traveling waves and oscillatory (stationary) soliton-like structures. This equation can be derived from the Lagrangian with the following $L^1$ potential:
\begin{align*}
L = \text{Kinetic}-\text{Potential} = \frac{1}{2}  |u_t|^2 - \frac{1}{2} |\nabla u|^2  - |u|.
\end{align*}
The equation of motion can be derived from the Lagrangian:
\begin{align*}
u_{tt} - \Delta u &= -p(u)\\
u(x,0)&=g_1(x)\\
u_t(x,0)&=g_2(x),
\end{align*}
which is the same as Equation~\eqref{signum_gordon} by replacing the $\text{sign}(u)$ term with the subgradient $p(u)$. 

To discretize the problem, we apply the ideas from the proximal gradient method, by placing $p(u)$ in the future: 
\begin{align*}
u^{n+1}-2u^{n}+u^{n-1} - \tau^2 \Delta u^n &= -\tau^2 p(u^{n+1}),
\end{align*}
and thus,
\begin{align*}
u^{n+1}= S( 2u^{n}-u^{n-1} +\tau^2 \Delta u^n, \tau^2).
\end{align*}

\begin{figure}
\includegraphics[width = 2 in]{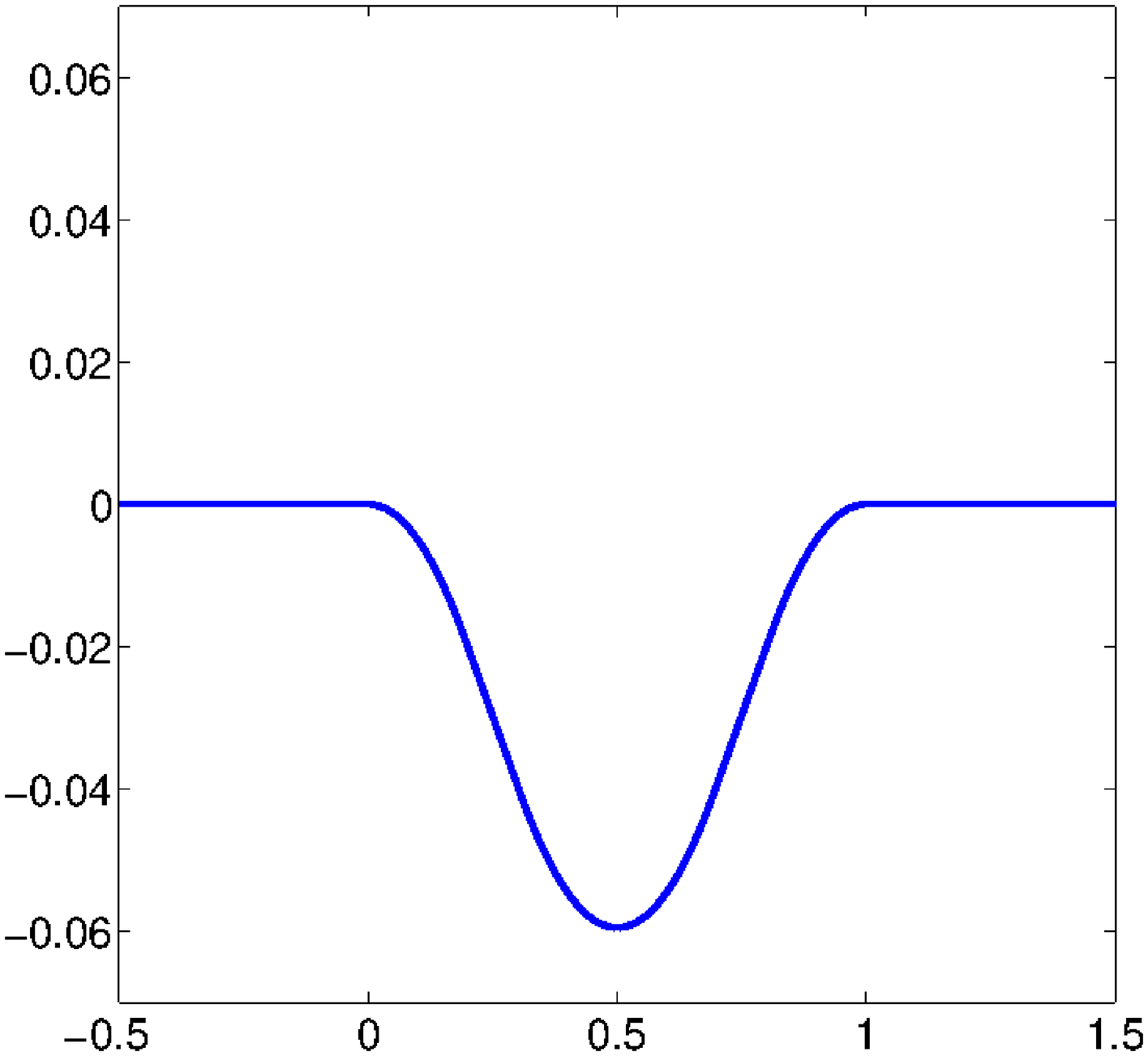}\quad
\includegraphics[width = 2 in]{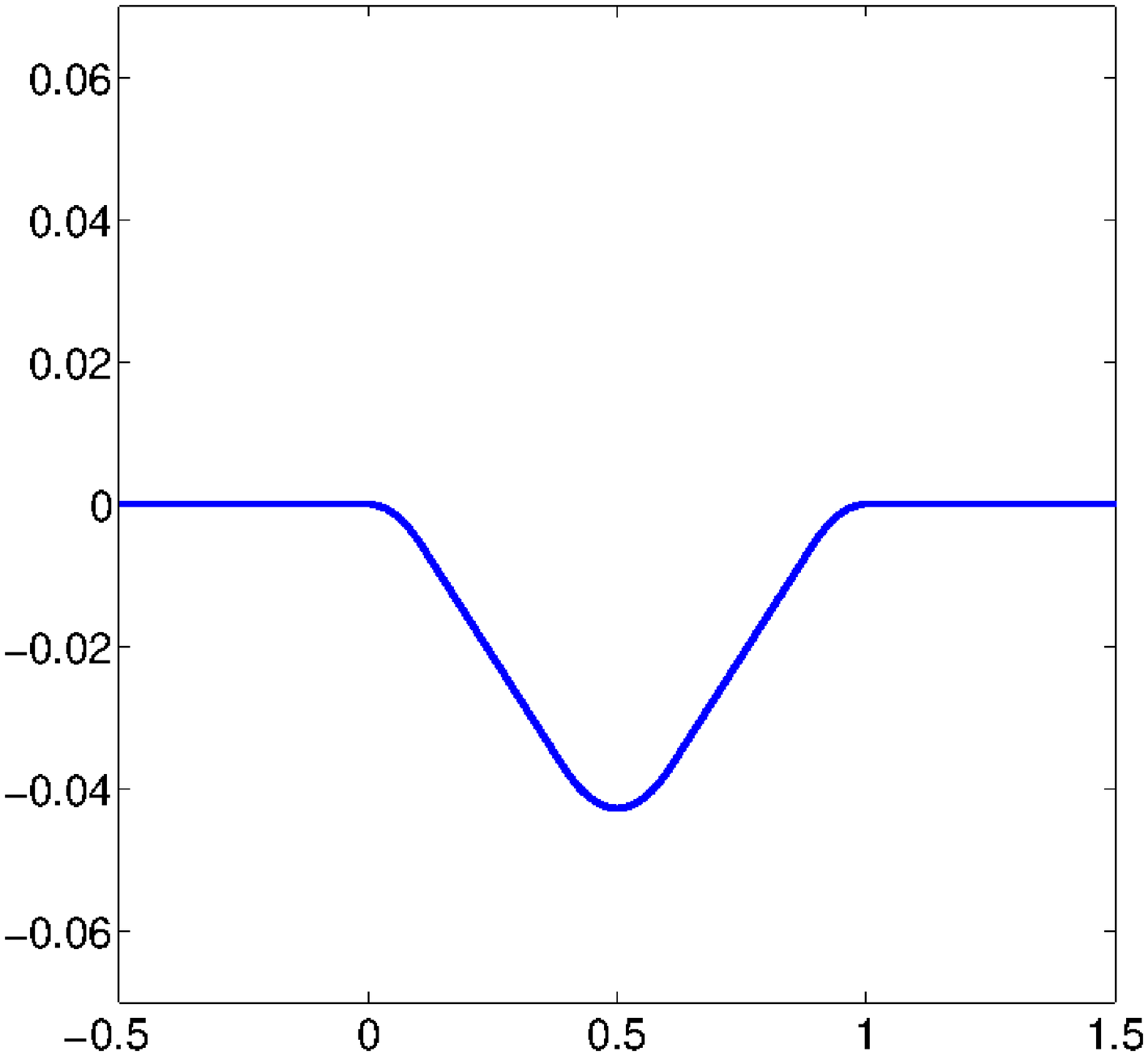}\\
\includegraphics[width = 2 in]{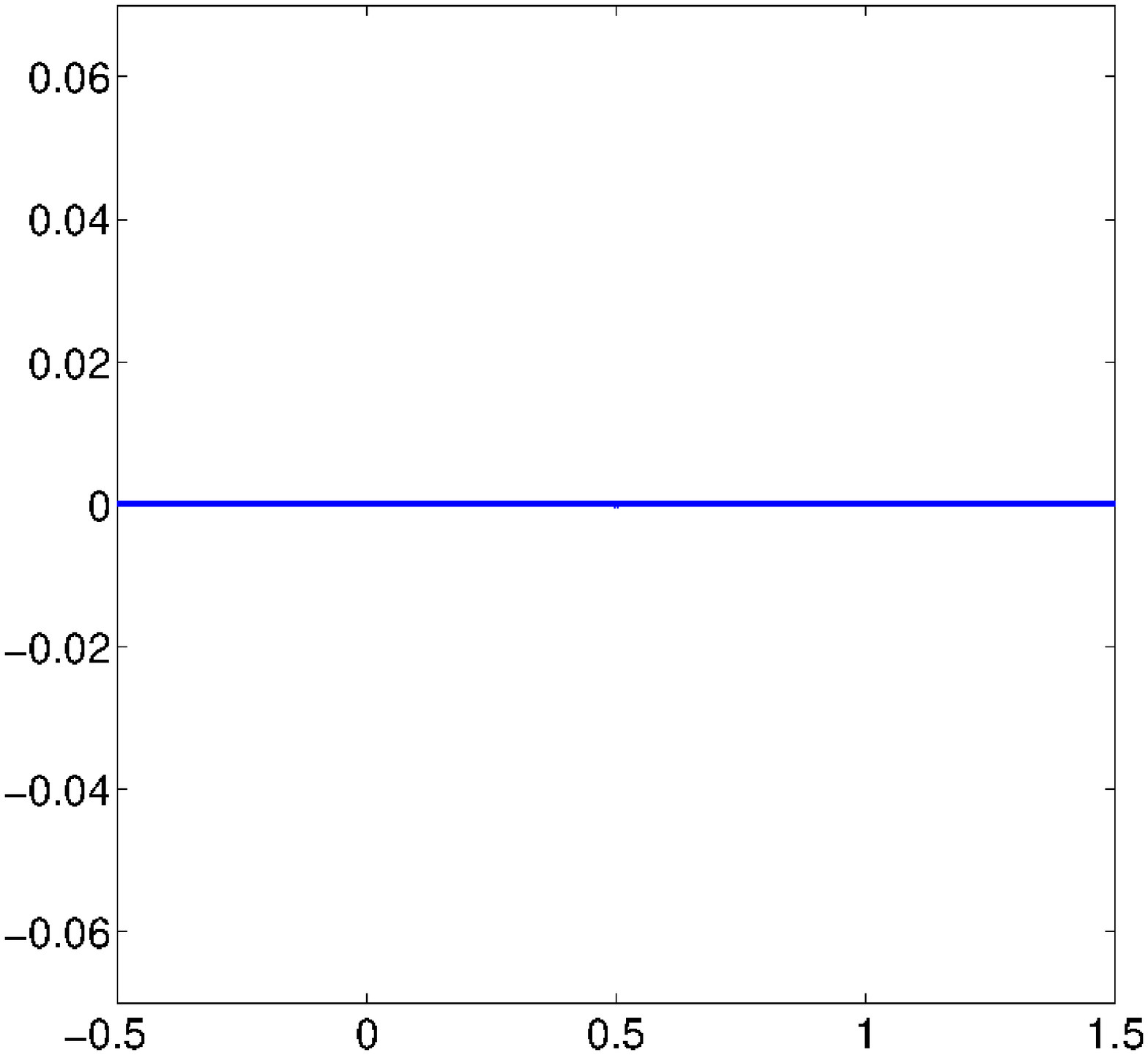}\quad
\includegraphics[width = 2 in]{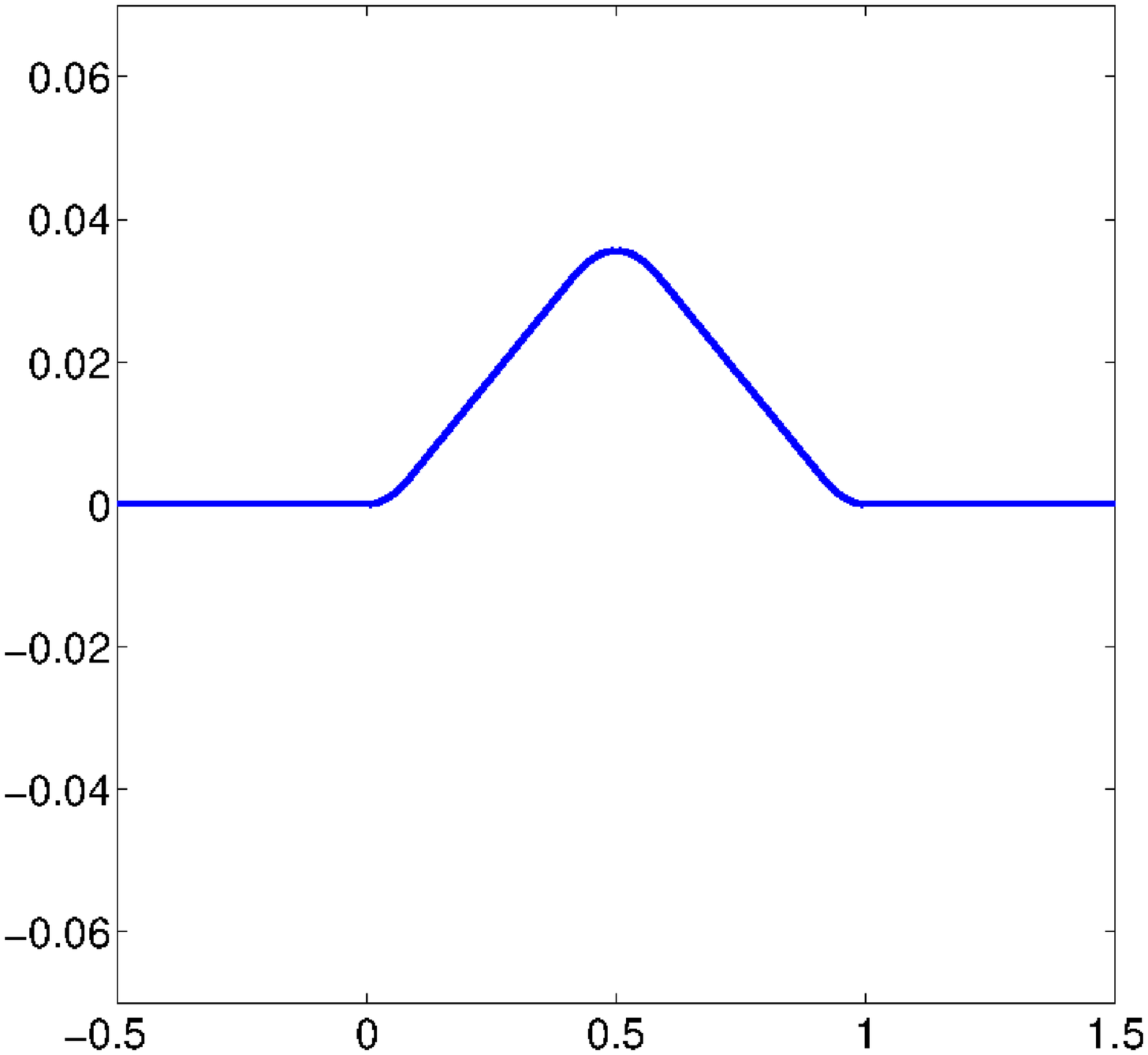}\\
\includegraphics[width = 2 in]{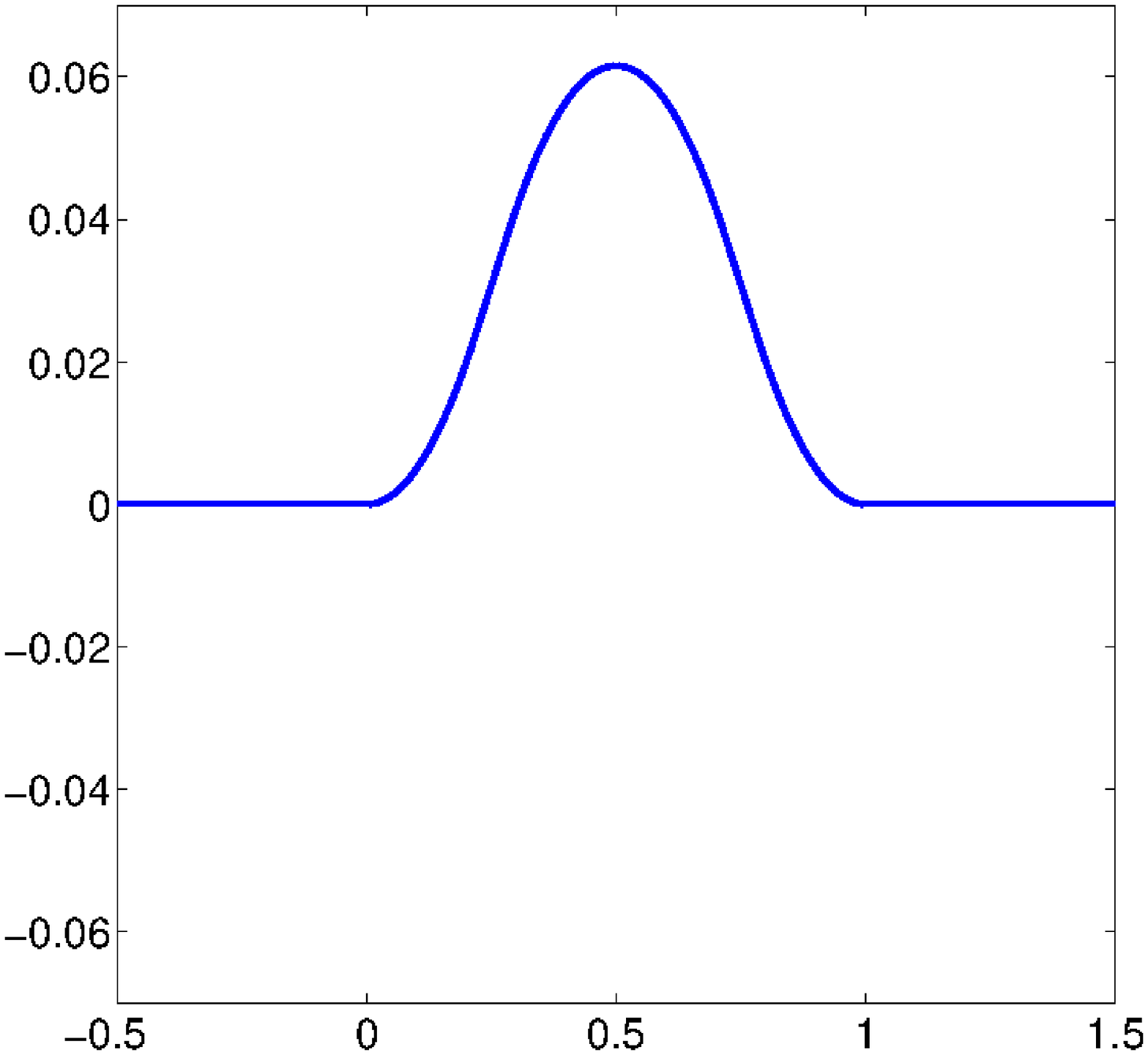}
\caption{The dynamics of an oscilatory compact solution of the Signum-Gordon equation. }\label{o:sg}
\end{figure}

In Figure~\ref{tw:sg}, we plot our numerical approximation to the traveling wave solution found in \cite{arodzWave}. Since the traveling wave profile is also known analytically, we show numerical convergence of our scheme as $h \rightarrow 0^+$ (see Table 1). Also, in Figure~\ref{o:sg}, we show the time evolution of an oscillatory compact soliton-like structure which appears in \cite{arodz2008compact, arodz2011swaying}. These examples show the range of behaviors that appear via the addition of an $L^1$ subgradient term.

\subsection{Divisible Sandpile}
\begin{figure}
\includegraphics[width = 2.25in]{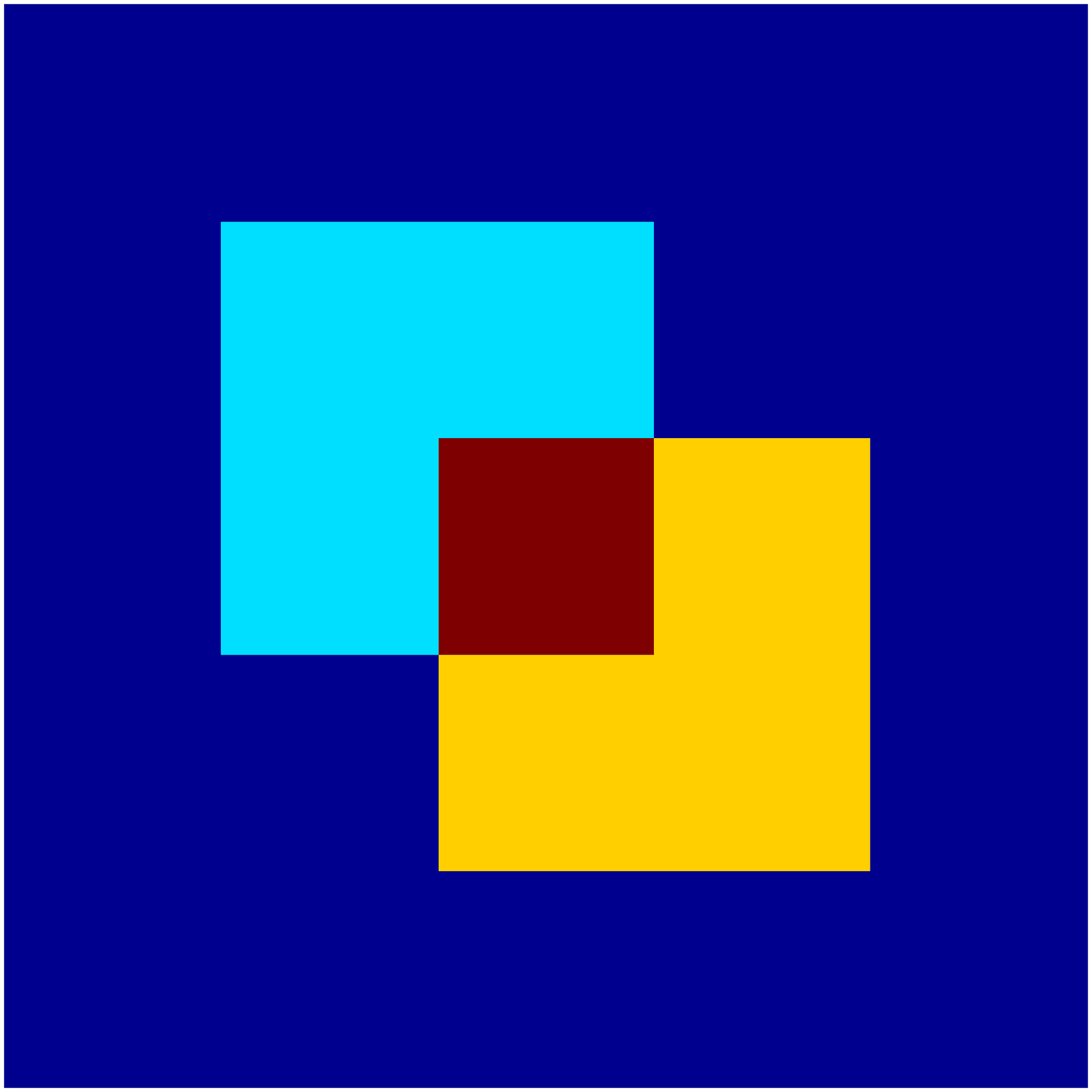}\quad
\includegraphics[width = 2.25in]{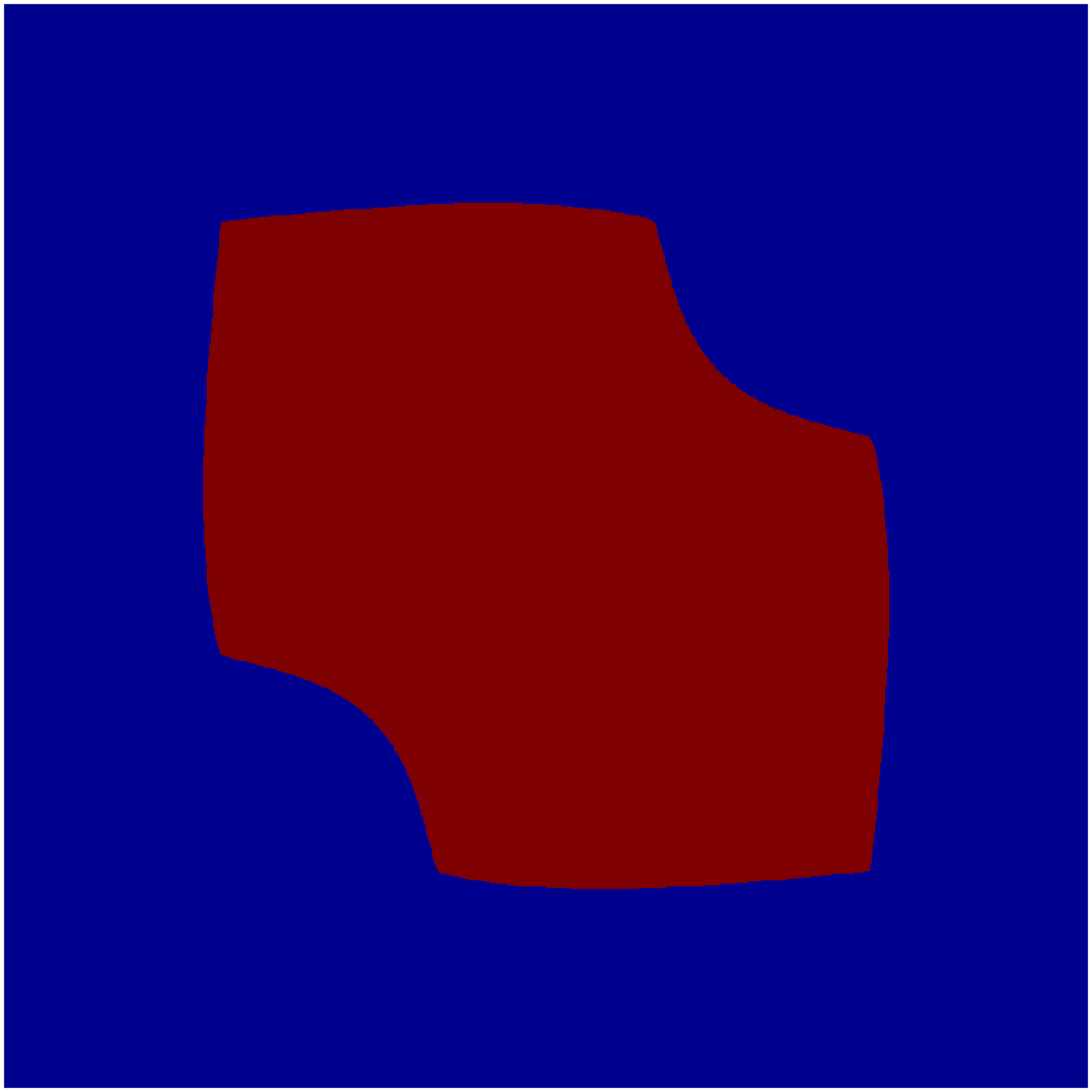}
\caption{A two-region sandpile problem, where each of the larger squares defines the set $S_j$, for $j=1,2$. The darker blue region has no mass, the lighter blue and yellow region has a mass density of 1,   and the red (overlap) region has a density of 2. On the left, the region of positive mass is displayed.}
\label{Sandpile}
\end{figure}

\begin{figure}
\includegraphics[width = 2.25in]{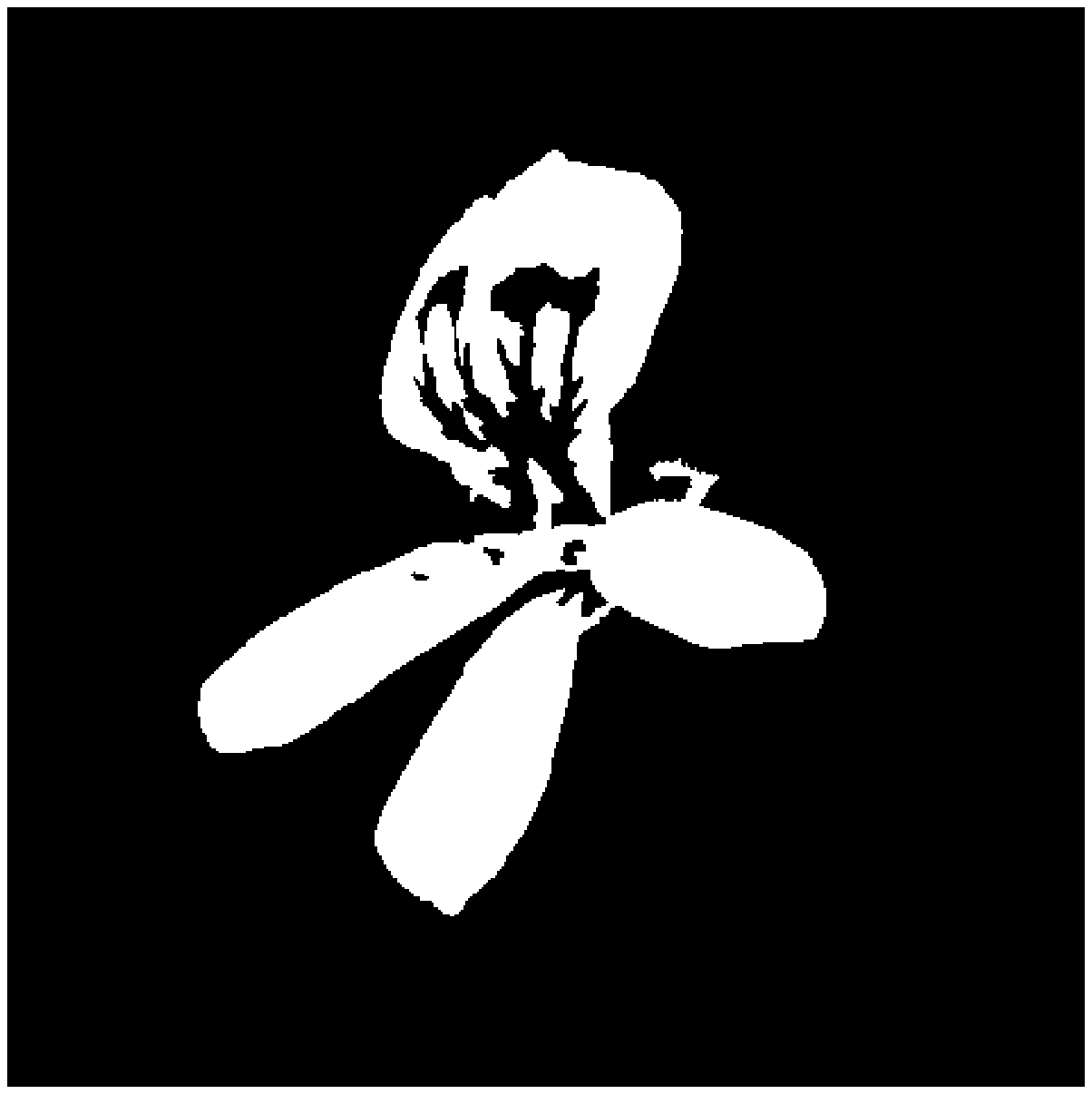}
\includegraphics[width = 2.25in]{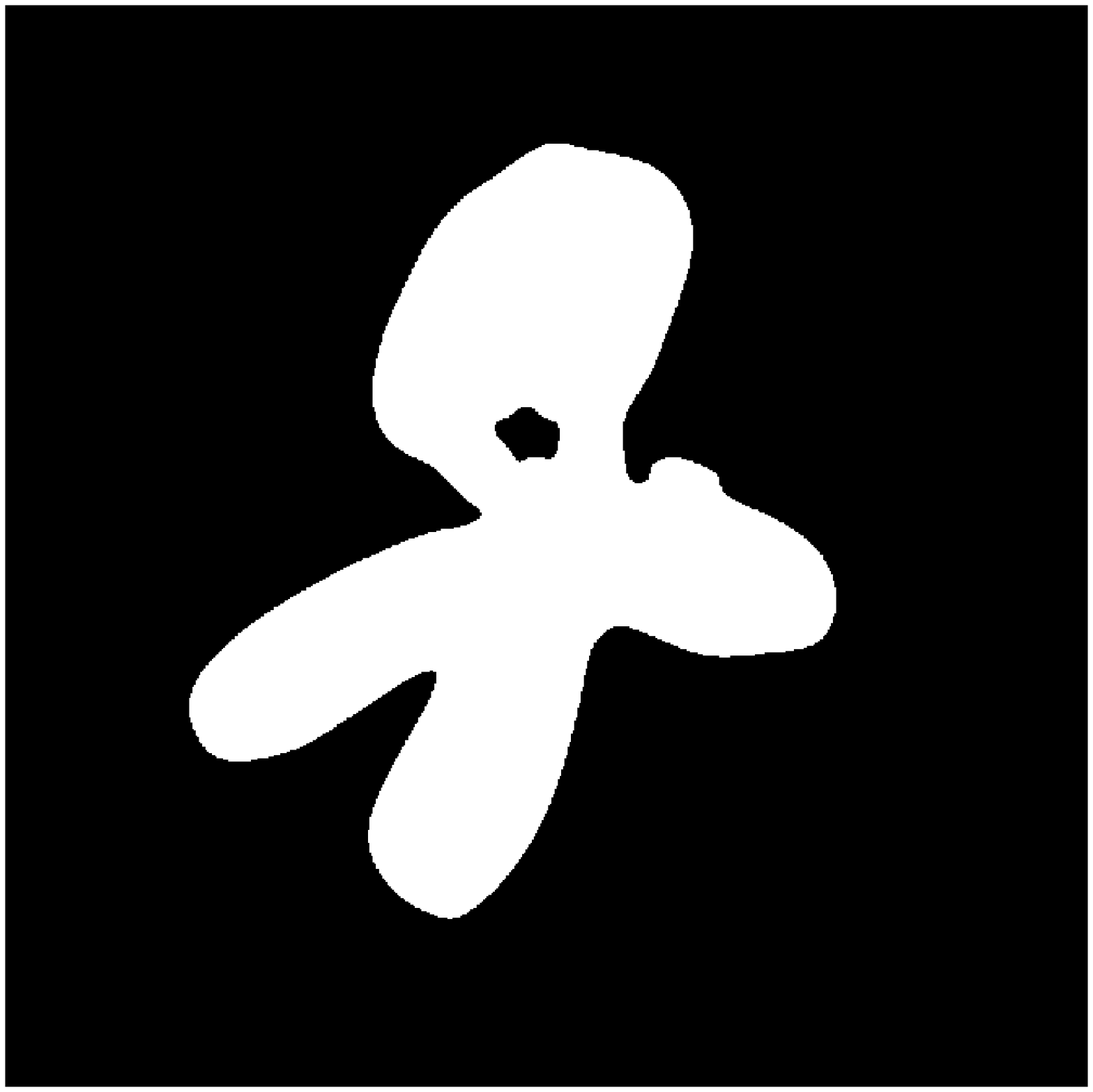}\\
\includegraphics[width = 2.25in]{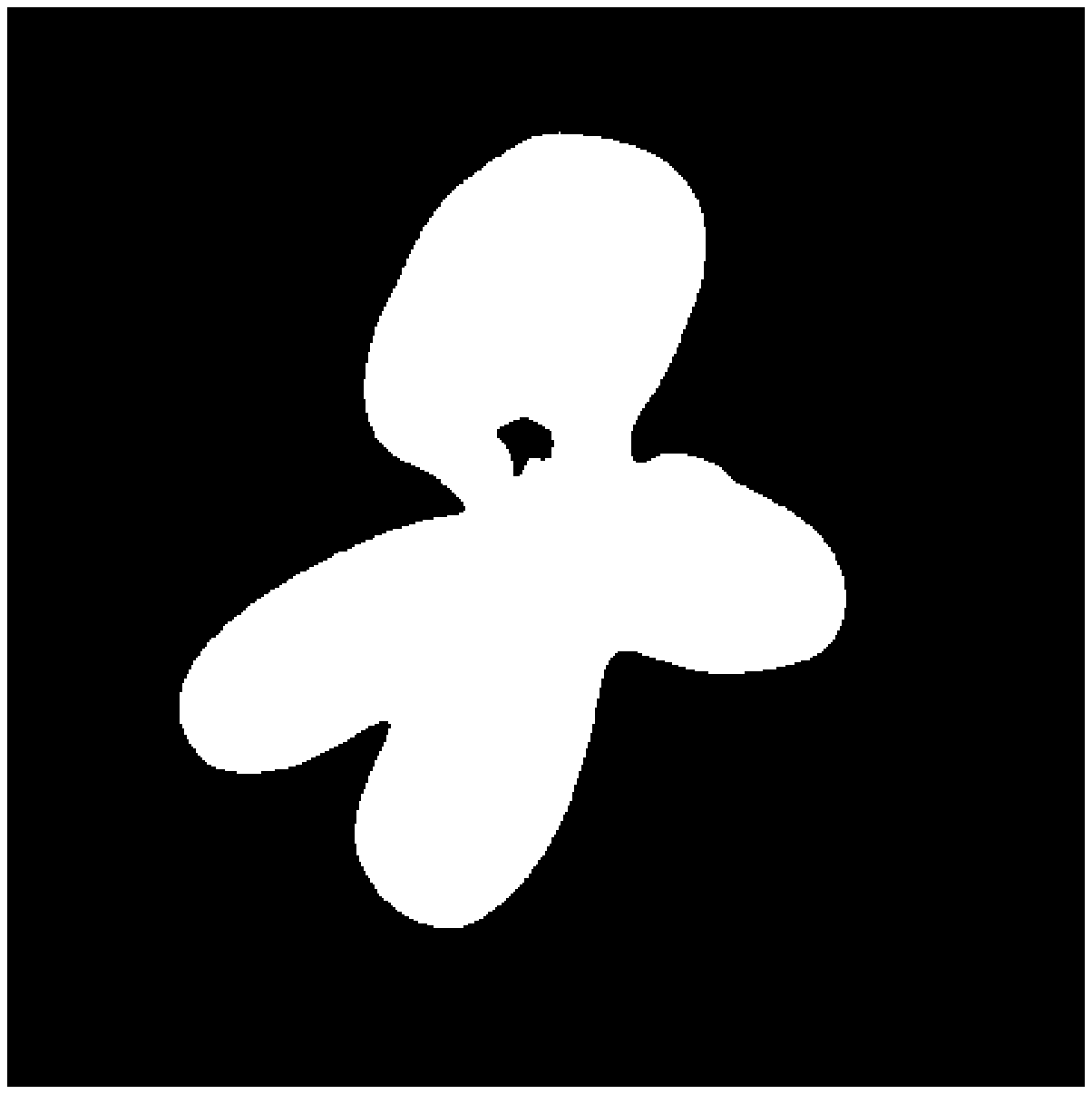}
\includegraphics[width = 2.25in]{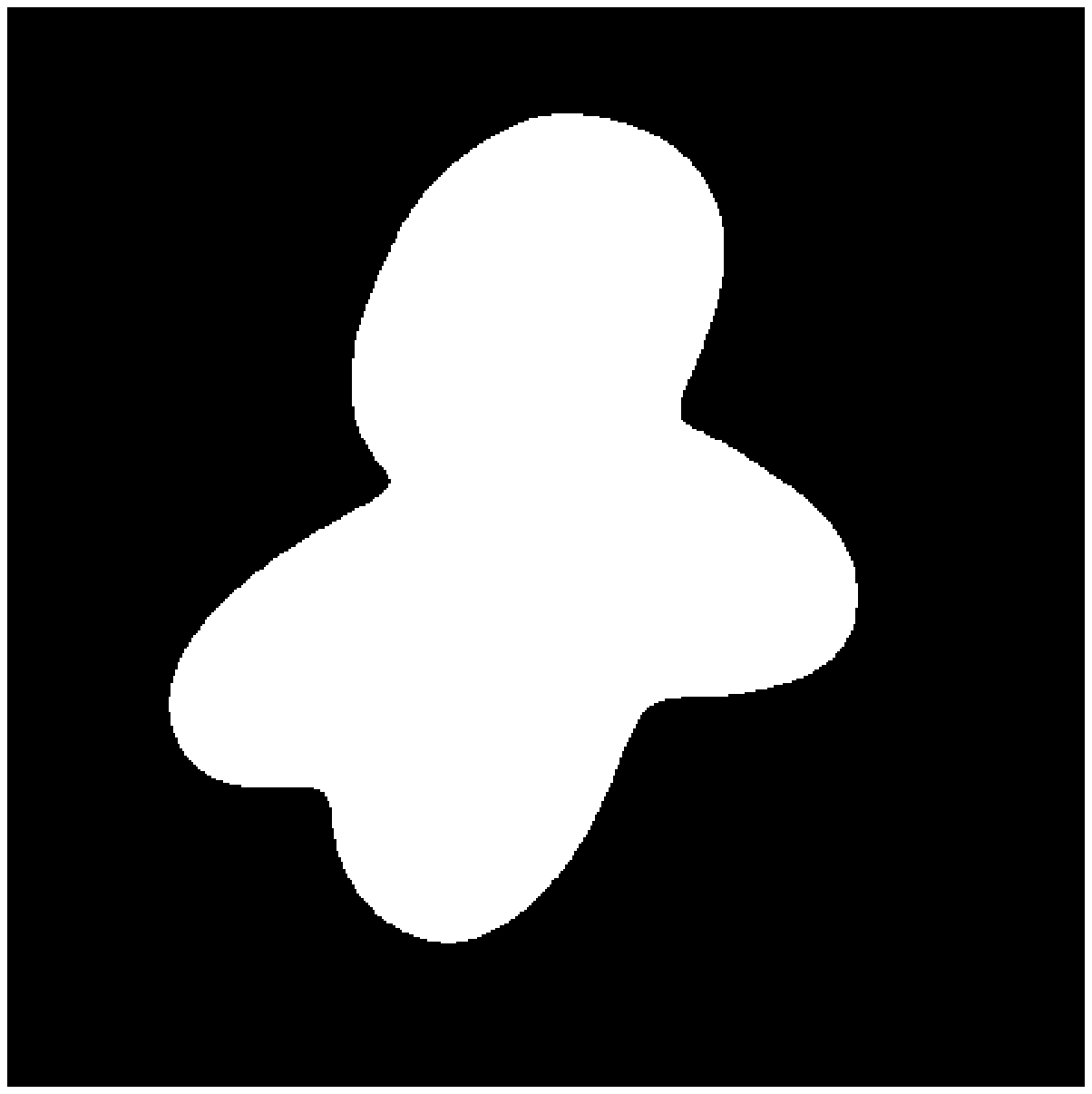}
\caption{The iterative evolution of our sandpile problem algorithm applied to a flower-shaped region $S$ on the top left with $f = 2\chi_S$. The final state appears in the bottom right corner. }
\label{Sandpile:flower}
\end{figure}

\begin{figure}
\includegraphics[width = 2.25in]{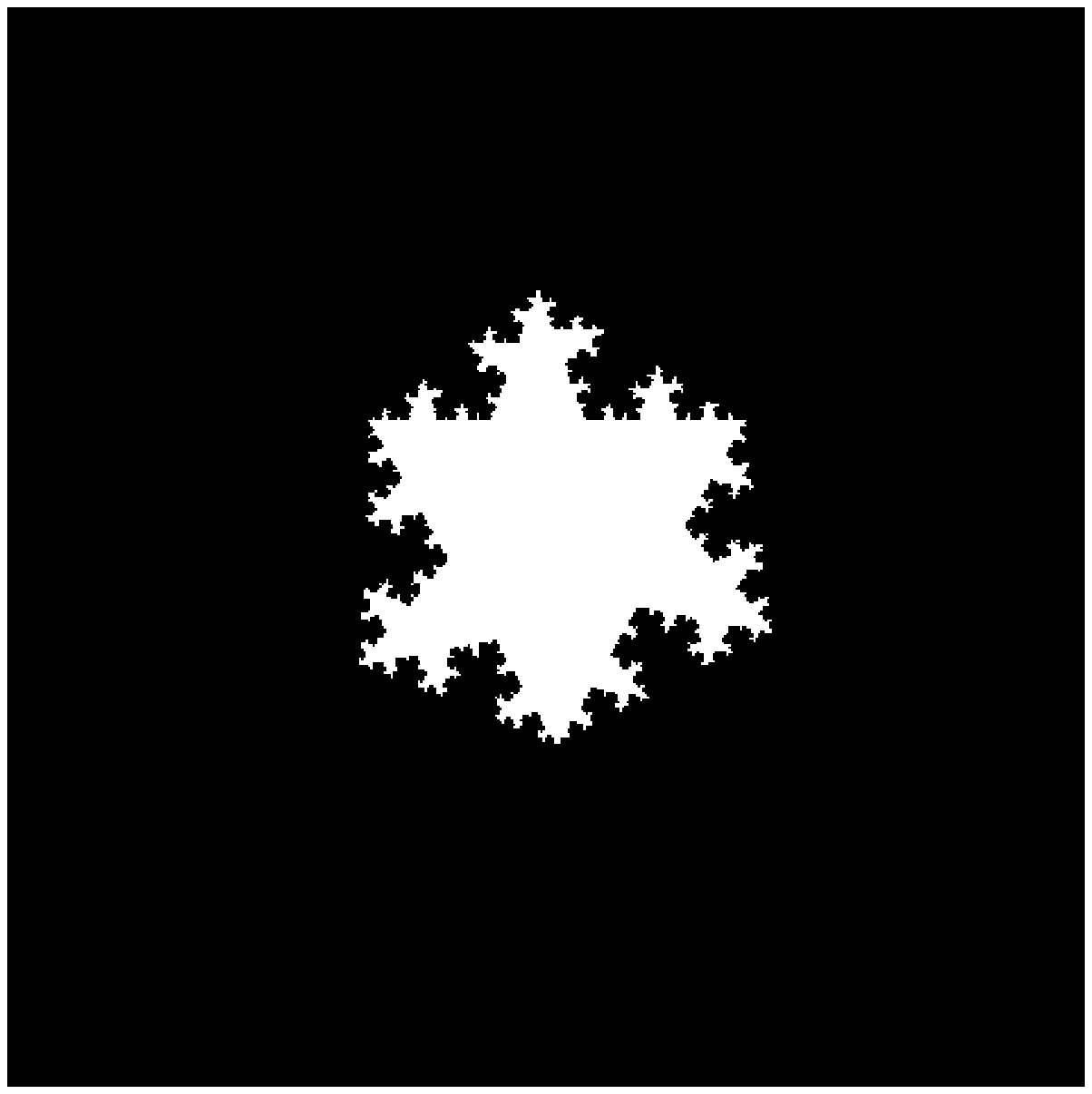}
\includegraphics[width = 2.25in]{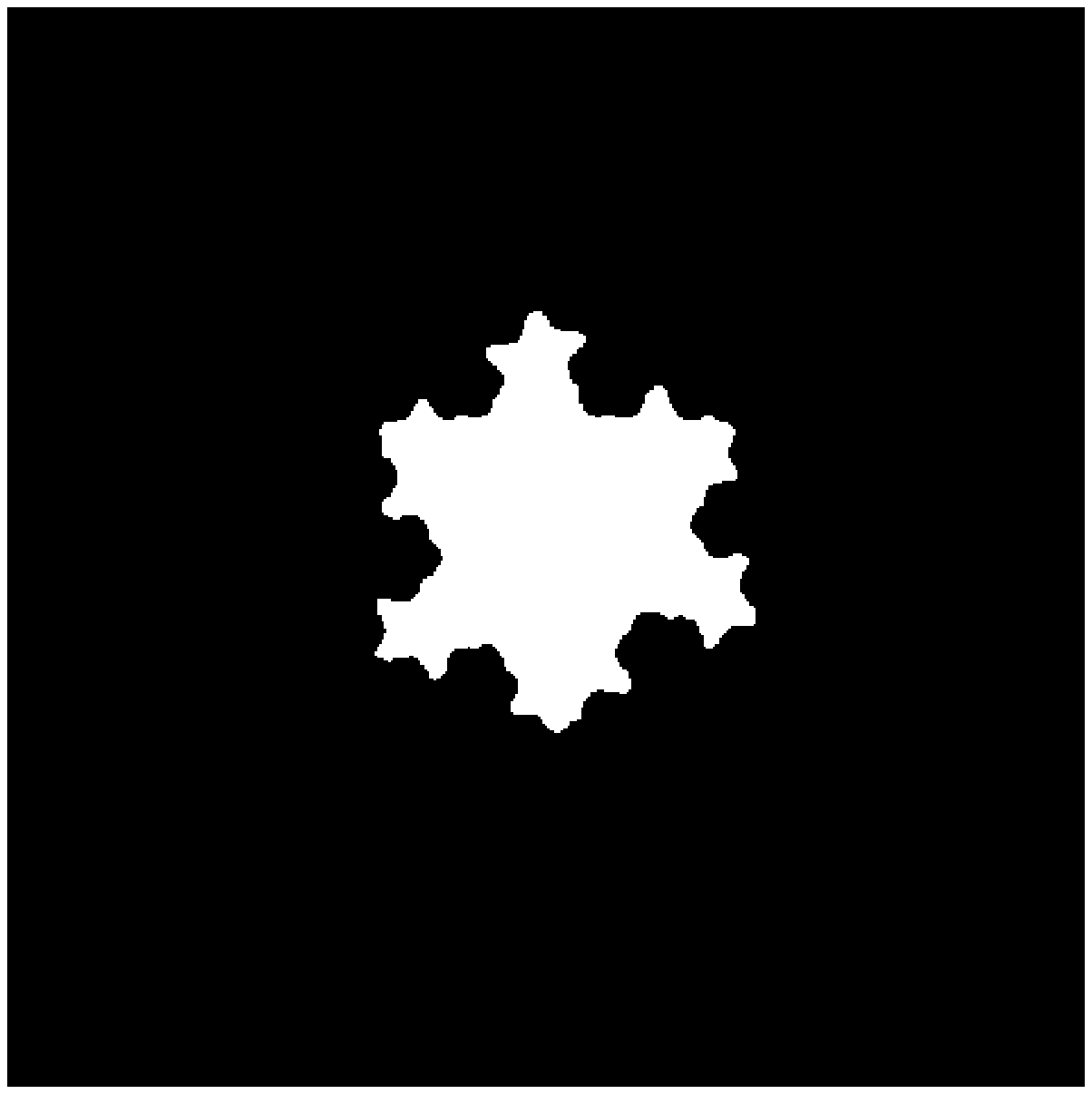}\\
\includegraphics[width = 2.25in]{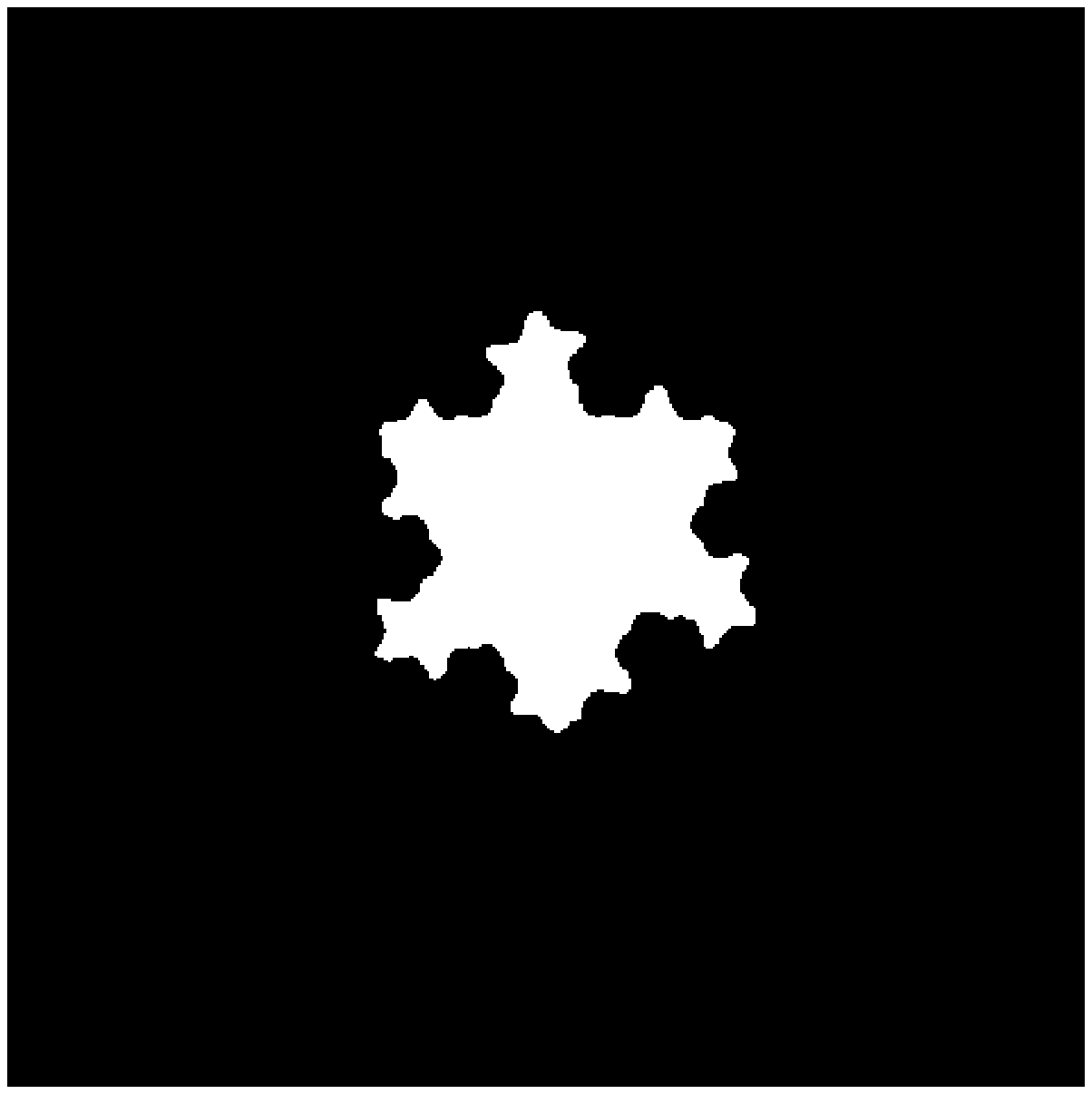}
\includegraphics[width = 2.25in]{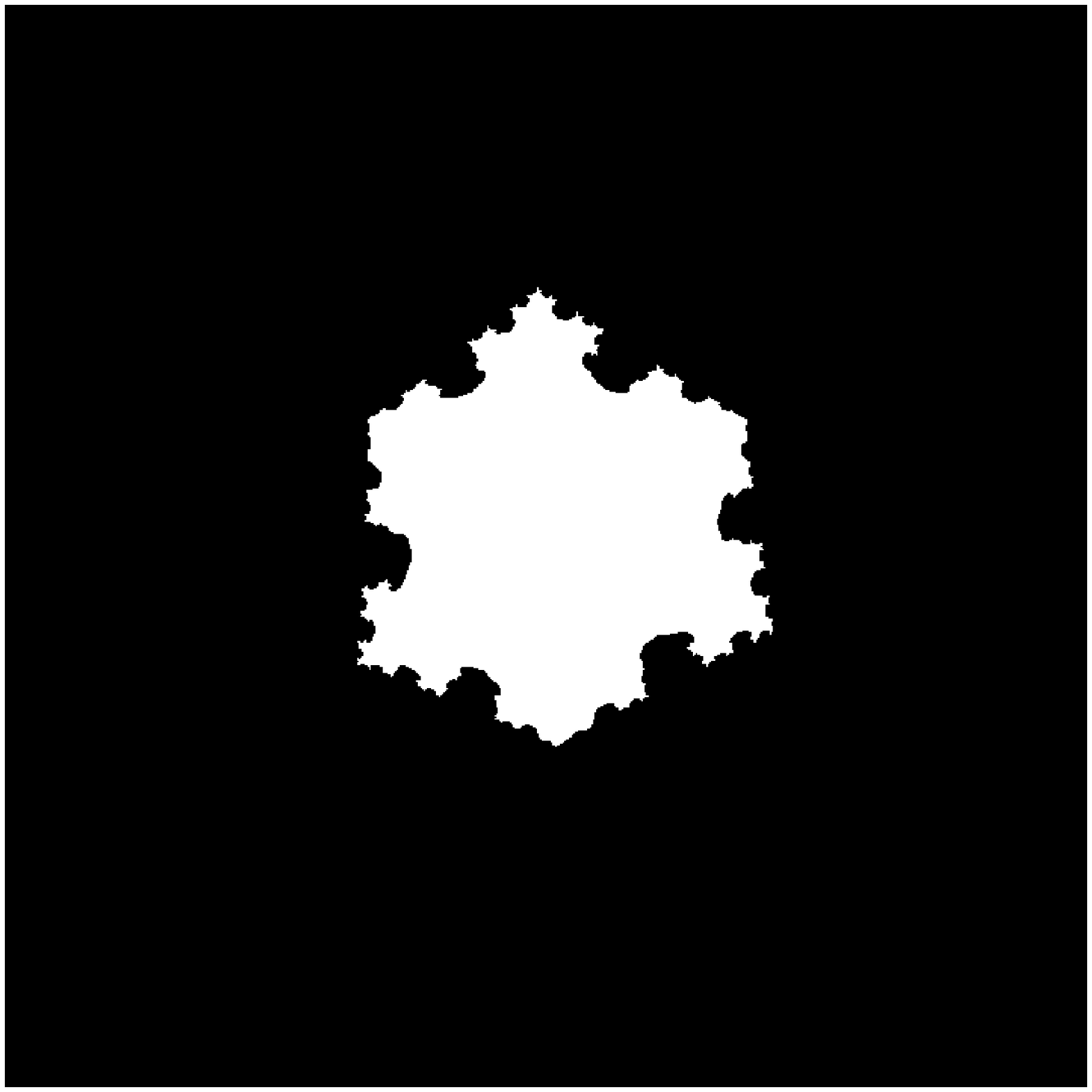}
\caption{The iterative evolution of our sandpile problem algorithm applied to the fractal region $S$ on the top left with $f = 1.2\chi_S$. The final state appears in the bottom right corner.}
\label{Sandpile:snow}
\end{figure}

\begin{figure}
\includegraphics[width = 2.25in]{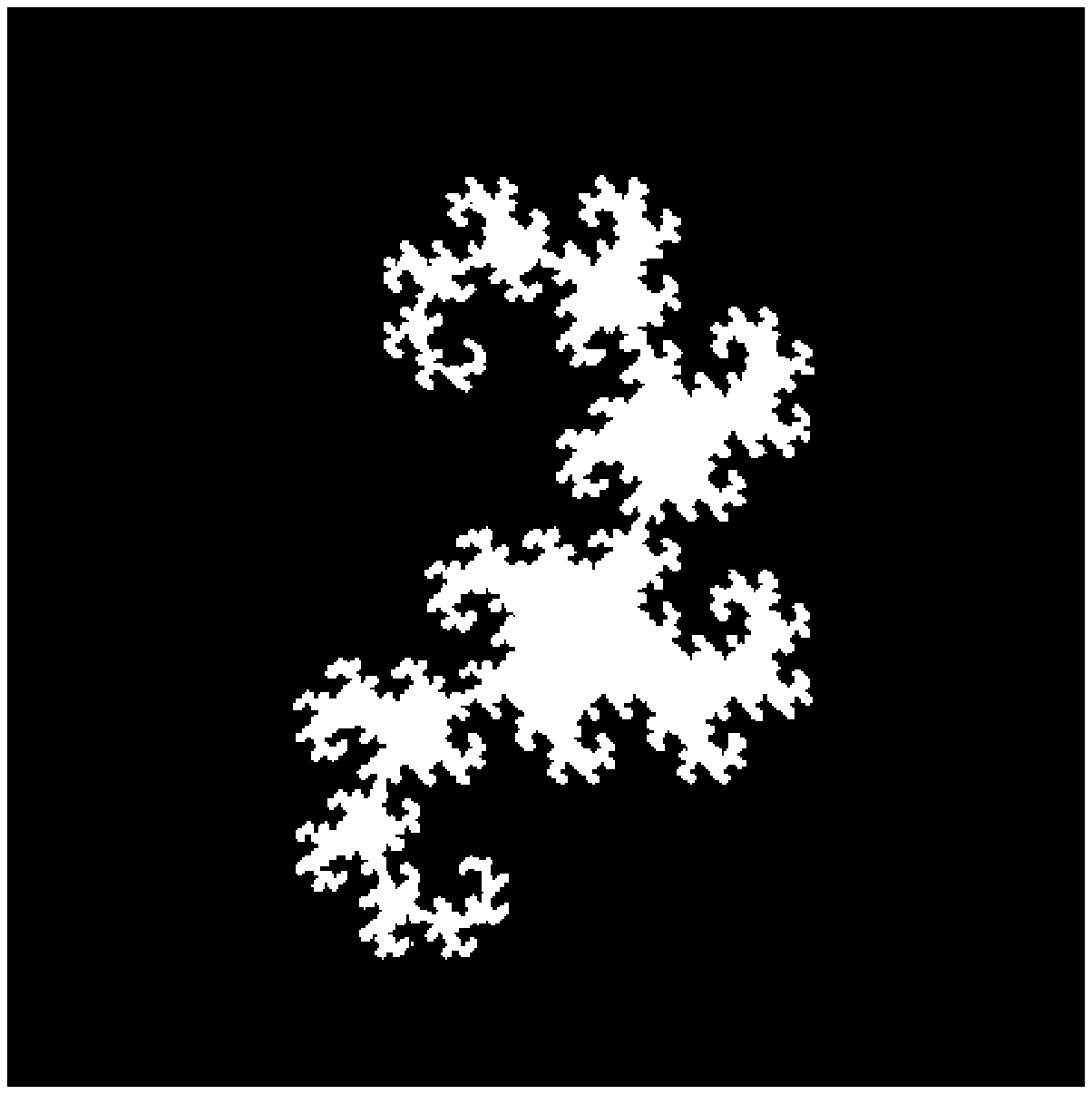}
\includegraphics[width = 2.25in]{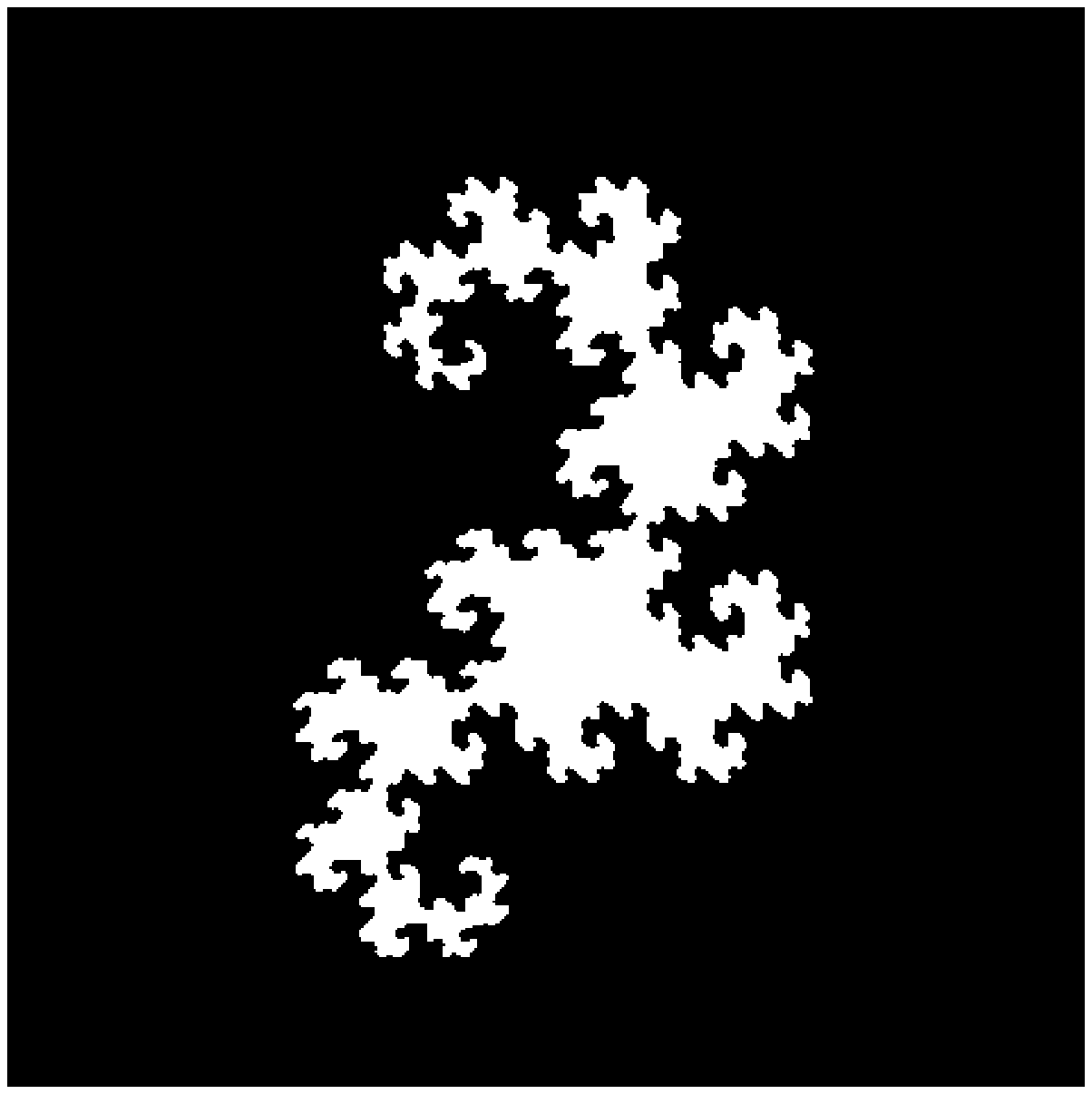}\\
\includegraphics[width = 2.25in]{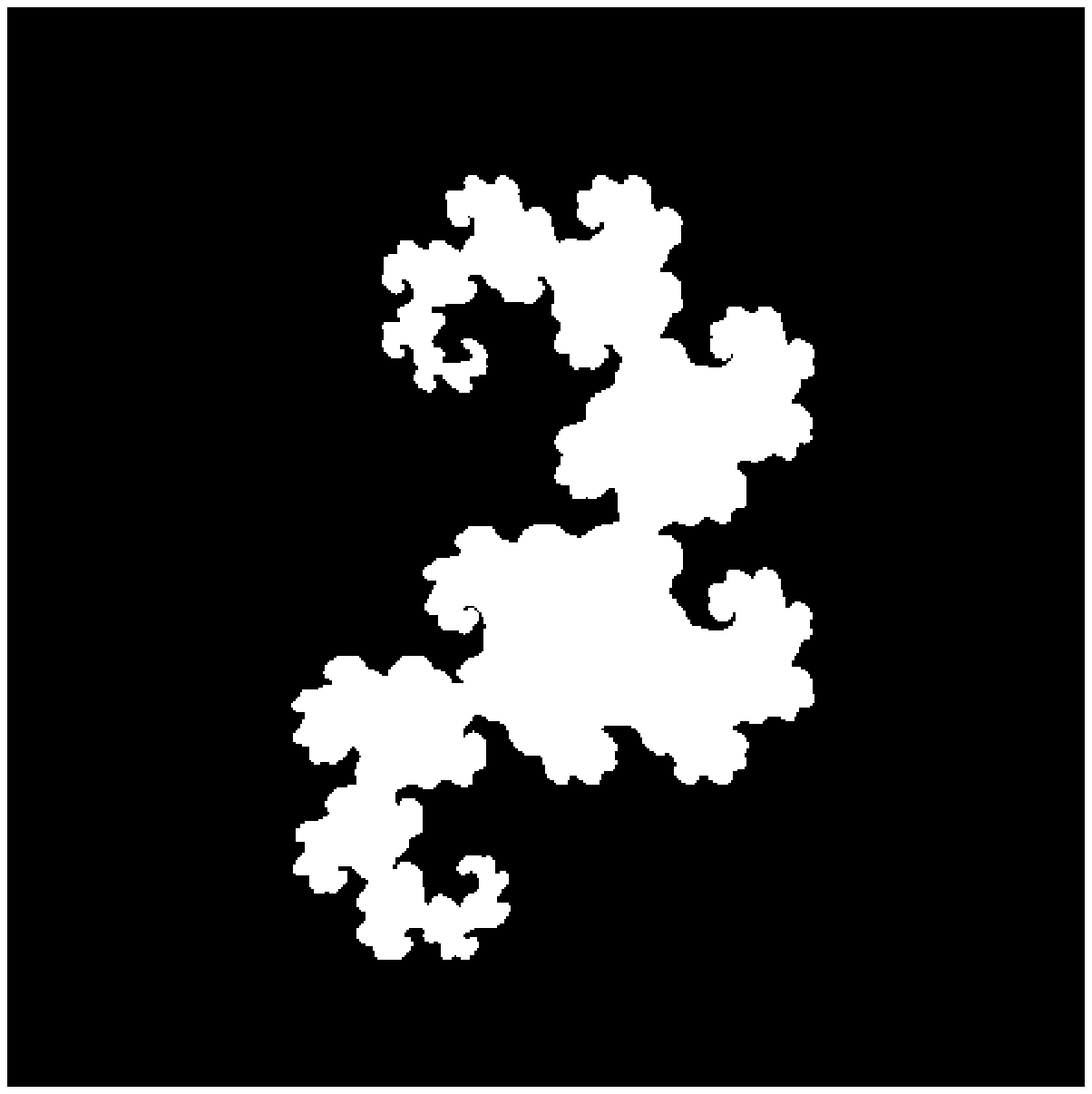}
\includegraphics[width = 2.25in]{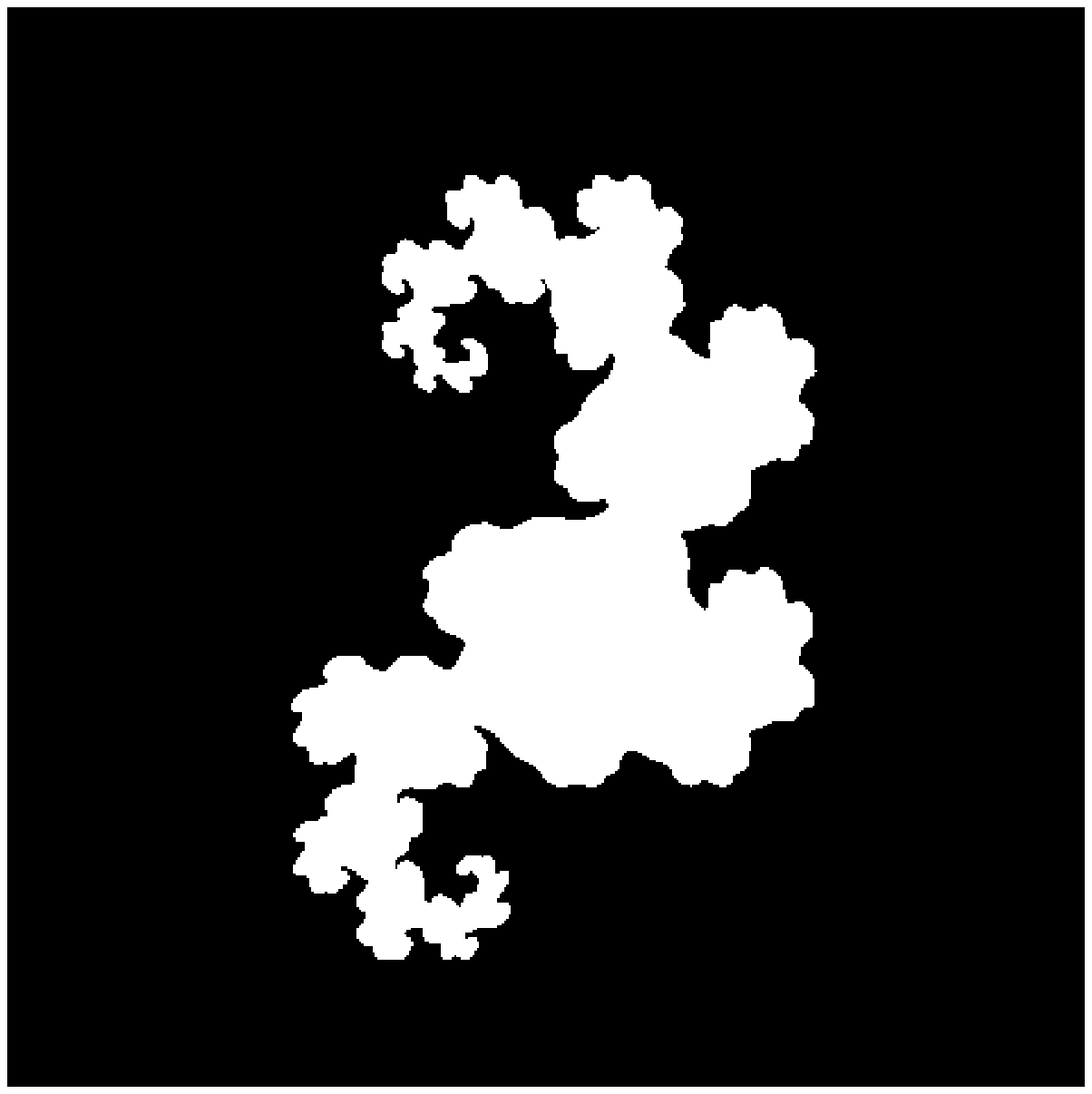}
\caption{The iterative evolution of our sandpile problem algorithm applied to the fractal region $S$ on the top left with $f = 1.5\chi_S$. The final state appears in the bottom right corner.}
\label{Sandpile:dragon}
\end{figure}

\begin{figure}
\includegraphics[width = 3.5in]{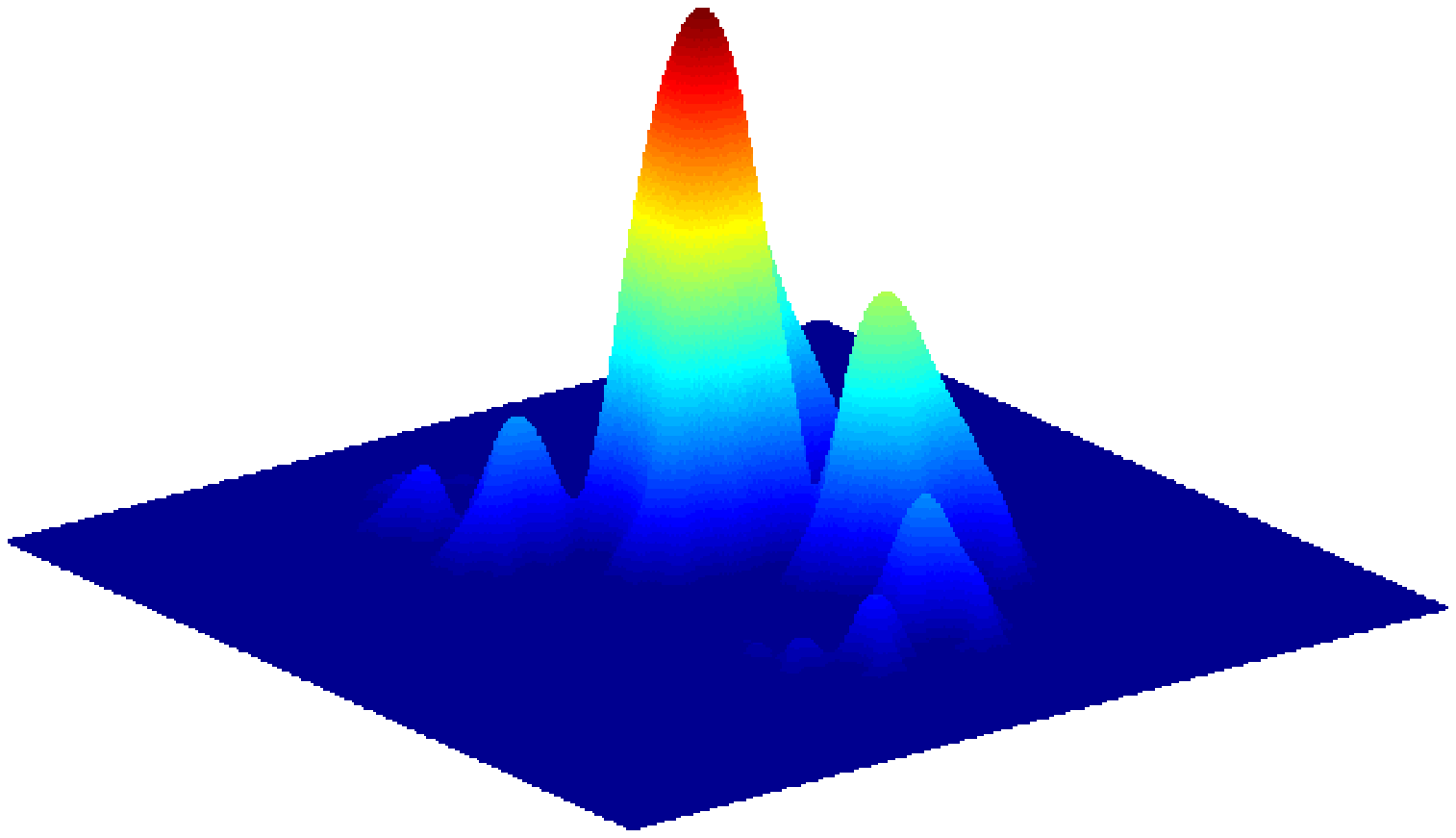}
\caption{The solution $u$ from Figure~\ref{Sandpile:dragon} (bottom right).}
\label{Sandpile:dragon2}
\end{figure}

As a model for self-assembly and internal diffusion limited aggregation, the sandpile problem has received attention recently \cite{pegden2013convergence,LevineSandpile, LevineSandpile2,LevineApollonian1,LevineApollonian2}. The problem is posed discretely, but has the following continuous formulation for the divisible sandpile problem \cite{LevineSandpile, LevineSandpile2}:
\begin{align} 
\Delta u = 1 - f, \text{ \ if } u \geq 0,
\label{Sandpile00}
\end{align}
where $f$ is some non-negative external force. By multiplying Equation~\eqref{Sandpile00} with $u$ and integrating over $\mathbb{R}^2$, the associated variational energy is: 
\begin{align} 
\min_u \int_{u \geq 0} \frac{1}{2} |\nabla u|^2+u - u f \text{\ dx}.
\label{sp}
\end{align}
There are several choices for relaxing the constraint $u\geq 0$, in particular, we use the following:
\begin{align} 
\min_u \int \frac{1}{2} |\nabla u|^2+|u| - u f \text{\ dx}.
\label{L1sp}
\end{align}
It can be shown (via maximum principle) that for $f \geq 0$ the solution of Equation~\eqref{L1sp} and Equation~\eqref{sp} are the same. The Euler-Lagrange equation for the $L^1$ sandpile problem is:
\begin{align} 
\Delta u = p(u) - f, 
\label{L1sp1}
\end{align} 
and is solved numerically via the Douglas-Rachford algorithm (see Equation~\eqref{Douglas_Rachford}). Note that if the external force  is a finite sum of characteristic functions
$f = \sum^N_{j=1} \alpha_j \chi_{S_j}$ where $S_j$ are compact sets and $\alpha_j\geq 0$, then by integrating Equation~\eqref{L1sp1} over $ \mathbb{R}^2$  we get:
\begin{equation}
|\text{supp}(u)|  = \sum^N_{j=0} \alpha_j |S_j|,
\label{conmass}
\end{equation}
since $u\geq 0$ and $\text{supp}(u)$ is compact. This refers to preservation of mass.

In Figure~\ref{Sandpile}, we take $f = \chi_{S_1} + \chi_{S_2}$, where $S_1$ and $S_2$ are the two overlapping square domains (on the left). The support set of $u$, given in Figure~\ref{Sandpile} (right), agrees with direct numerical simulation of the discrete sandpile problem. The direct simulation follows  a topping rule described in \cite{LevineSandpile}. 

In Figures~\ref{Sandpile:flower}-\ref{Sandpile:dragon}, we take $f = \alpha \chi_{S}$ where $S$ is the shape given in  Figures~\ref{Sandpile:flower}-\ref{Sandpile:dragon} (the top left), and $\alpha = 2.0, 1.2$ and $1.5$, respectively. The support set of $u$ is given in  Figures~\ref{Sandpile:flower}-\ref{Sandpile:dragon} (the bottom right) with intermediate calculation shown in  Figures~\ref{Sandpile:flower}-\ref{Sandpile:dragon} (the remaining plots). To verify that the solutions from our algorithm correspond to the correct solutions for the sandpile problem, we use the mass  conservation property, Equation~\eqref{conmass}. Unlike direct simulation, our method also calculates the function $u$ as shown in Figure~\ref{Sandpile:dragon2}. One of the benefits of our approach is that the solutions can be computed quickly, for example, our method is at least $8$ times faster than direct simulation (76 seconds vs. 652 seconds) at approximating the solution found in Figure~\ref{Sandpile:dragon2}. 

\section{Conclusion} \label{conc_sec}
 
By adding the subdifferential of $L^1$ to certain PDEs, we have shown (numerically and theoretically) various properties of the solutions.  These problems arise from physical models as well as exact relaxation of other PDEs, and could provide useful tools in computing fast approximations to nonlinear problems with a compactly supported free boundary. This is all in the spirit of borrowing the key idea from compressed sensing, that $L^1$ regularization implies sparsity of discrete systems \cite{donoho2006compressed}, and transferring it to classical problems in PDE. See \cite{Schaeffer:2013sparse,Ozolins:2013CMs} for earlier work in this direction.

\section*{Acknowledgments}
The authors would like to thank Farzin Barekat, Jerome Darbon, William Feldman, Inwon C. Kim and James H. von Brecht for their helpful discussions and comments.

R. Caflisch was supported by ONR N00014-14-1-0444. S. Osher was supported by ONR N00014-14-1-0444 and N000141110719. H. Schaeffer was supported by NSF 1303892 and University of California Presidents Postdoctoral Fellowship Program. G. Tran was supported by ONR N00014-14-1-0444 and N000141110719.
\appendix

\bibliographystyle{plain}

\section{appendix}

\subsection{Proof of a free boundary formula in Section \ref{freebndysection}}

We derive the short time asymptotic equation for the support set Equation~\eqref{parabolic_pde}. First, we provide a natural boundary condition for the problem.
\\

\noindent \textbf{Flux Condition.} \textit{
Let $u(x,t) \in C^0(C^1({\mathbb{R}}); (0,T))$ and $u_t\in L^\infty(C^1({\mathbb{R}}); (0,T))$ be a solution to 
\begin{equation}
u_t-u_{xx} =h(x,t,\gamma).
\label{eq:heat}
\end{equation}
Assume that there exists a positive valued function $a\in C^1(0,T)$ such that $h=0$ for $|x|>a(t)$, $g=0$ for $|x| >a(0)$, and the exterior mass,
$$m(t) = \int\limits_{a(t)}^{\infty} u(x,t)dx,$$
is conserved, then $u(a(t),t)=0$ and $u_x(a(t),t)=0$.
}

To derive this condition, consider the heat equation \eqref{eq:heat}. Differentiate the one sided mass in time yields:
\begin{equation*}
\begin{aligned}
\dfrac{dm}{dt} &= -u(a(t),t)a'(t) + \int\limits_{a(t)}^{\infty} u_t(x,t)dx\\
 &= -u(a(t),t)a'(t) + \int\limits_{a(t)}^{\infty} u_{xx}(x,t)dx\\
 &= -u(a(t),t)a'(t) - u_{x}(a(t),t)\\
& = - F(t),
\end{aligned}
\end{equation*}
in which $F$ is the flux across the moving boundary $x=a(t)$.

We now can see that if the flux across a moving boundary $x=a(t)$ is zero (\textit{i.e.} the mass is conserved), we have
\begin{equation}
 F(t) = u(a(t),t)a'(t) + u_{x}(a(t),t) =0.
\label{eq:BC}
\end{equation}
This is the natural boundary condition for this problem. In the time-dependent region ${\mathcal F}=\left\{(x,t): x>a(t)\right\}$, the initial data $g$, force $h$ and and incoming flux $F$ are all zero, so that the solution is identically zero. In particular, $u =u_x =0$ on $x=\pm a(t).$

Next, consider the following equation:
\begin{equation*}
	u_t - u_{xx} =
	\begin{cases}
       & f(x) -\gamma,\quad |x|<a(t)\\ 
       & 0,\hspace{1.5 cm}|x|>a(t) 
	\end{cases}\label{eq:freeBoundary}
\end{equation*}
\begin{equation*}
u(x,0) = 0.
\end{equation*}
For simplicity assume that $f(x) = f(|x|)$ and $f$ is a decreasing function with $f(|x|)\rightarrow 0$.  Denote $a_0\geq 0$ such that $f(a_0) =\gamma$ and w.l.o.g. $f_x(a_0) \neq 0$. By studying the exterior mass of Equation \eqref{eq:freeBoundary}, we want to show that in small time:
\begin{equation*}
a(t) = a_0+a_1\sqrt {t}, 
\end{equation*}
for some $a_1\geq 0$. 

We look for an increasing function $a(t)$ such that the exterior mass of Equation \eqref{eq:freeBoundary} is zero:
\begin{equation*}
m(t) = \int\limits_{a(t)}^{\infty} dx  \int\limits_0^t ds\int\limits_{-a(s)}^{a(s)} G(x-y,t-s)(f(y)-\gamma) dy.\label{eq:mass}
\end{equation*}
where we use the Greens formula to represent $u$. Since $a(t)$ is an increasing function, we have
\[y\leq a(s) \leq a(t)\leq x.\]
Therefore, for $t$ small, the Green's function $G(x-y,t-s)$ is sharply peaked near the point 
$$y=a(t), \ s =t , \ x = a(t).$$
So we can replace $(f(y) -\gamma)$ by the first few terms in its Taylor expansion
\begin{equation*}
f(y) -\gamma = (y-a_0)f_1 +\mathcal{O}((y-a_0)^2),
\end{equation*}
in which $f_1 = f_x(a_0)$. Also, since $G(x-y,t-s)$ decays exponentially as $y\rightarrow -\infty$, we replace the lower limit $y=-a(s)$ by $-\infty.$ Now the mass can be approximated by
\begin{equation*}
m(t) = f_1\int\limits_{a(t)}^{\infty} dx  \int\limits_0^t ds\int\limits_{-\infty}^{a(s)} (y-a_0) G(x-y,t-s) dy
\end{equation*}
Next we show the existence of $a_1$ satisfying the following approximations
\begin{equation*}
a(t) = a_0+a_1\sqrt{t},\quad\text{and}\quad m(t) = 0.
\end{equation*}
We change the variables to
\begin{equation*}
\begin{aligned}
x = & x_1\sqrt{t}  + a_0,& x_1\in [a_1,\infty),\\
y = & y_1 \sqrt{t}  + a_0,& y_1\in  (-\infty,a_1\sqrt{s_1} \,],\\
s = & s_1 t,& s_1\in[0,1],
\end{aligned}
\label{changevariable}
\end{equation*}
and 
\begin{equation*}
\begin{aligned}
x_1 = & x_2a_1,& x_2\in [1,\infty),\\
y_1 = & y_2 a_1,& y_2\in  (-\infty,\sqrt{s_1} \,],
\end{aligned}
\label{changevariable2}
\end{equation*}
and note that $G(x-y,t-s)=t^{-1/2} G(x_1-y_1,1-s_1)=t^{-1/2}G(a_1(x_2-y_2),1-s_1)$. Then
\begin{equation*}
\begin{aligned}
m(t) = & f_1 t^2\int\limits_{a_1}^{\infty} dx_1  \int\limits_0^1 ds_1\int\limits_{-\infty}^{a_1\sqrt{s_1}} y_1 G(x_1-y_1,1-s_1) dy_1
\\
= & a_1^3 f_1 t^2\int\limits_{1}^{\infty} dx_2  \int\limits_0^1 ds_1\int\limits_{-\infty}^{\sqrt{s_1}} y_2 G(a_1(x_2-y_2),1-s_1) dy_2.
\end{aligned}
\end{equation*}
Consider the rescaled masses $\widetilde{m_1}(a_1)=m(t)/(f_1 t^2)$ and $\widetilde{m_2}(a_1)=m(t)/(a_1^2 f_1 t^2)$; \textit{i.e.},
\begin{equation*}
\begin{aligned}
\widetilde{m_1}(a_1) = &\int\limits_{a_1}^{\infty} dx_1  \int\limits_0^1 ds_1\int\limits_{-\infty}^{a_1\sqrt{s_1}} y_1 G(x_1-y_1,1-s_1) dy_1,
\\
\widetilde{m_2}(a_1) =a_1 &\int\limits_{1}^{\infty} dx_2  \int\limits_0^1 ds_1\int\limits_{-\infty}^{\sqrt{s_1}} y_2 G(a_1(x_2-y_2),1-s_1) dy_2.
\end{aligned}
\end{equation*}

As $a_1\rightarrow 0$, $\widetilde{m_1}(a_1)$ goes to
\begin{equation*}
\widetilde{m_1}(0) = \int\limits_{0}^{\infty} dx_1  \int\limits_0^1 ds_1\int\limits_{-\infty}^{0} y_1 G(x_1-y_1,1-s_1) dy_1
\end{equation*}
with 
$\widetilde{m_1}(0)<0$. This shows that $m(t)<0$ for $a_1=0$.

On the other hand, for $a_1 \gg 1$, $a_1 G(a_1(x_2-y_2),1-s_1) $ is approximately the Dirac delta function at $x_2=y_2$, $s_1 = 1$.  At this point, we have $y_2>0$, therefore $\widetilde{m_2}(a_1) >0$. This shows that $m(t)>0$ for large values of $a_1$. Thus there exists a positive value $a_1$ so that $m(t) =0$.

\subsection{Proof of support size estimate in Section \ref{sec_compsupp}}
\begin{proof}
First, observe that if $\gamma \geq \max |f|,$ then the unique solution of Equation~\eqref{general_elliptic_pde}. is $u\equiv 0.$ Indeed, if $u=0$, since $\dfrac{f}{\gamma}\in [-1,1]$, we can choose $p(u) = \dfrac{f}{\gamma}$ and Equation~\eqref{general_elliptic_pde} is satisfied.

Now, take $\mathcal{S}=\text{supp}(u)$ and integrating both sides of Equation~\eqref{general_elliptic_pde} gives us
\begin{equation*}
\int_{\partial \mathcal{S}} M \nabla u \cdot N ds = -\int_{\mathcal{S}} f dx+ \gamma \, \text{sign}(u) |\mathcal{S}|.
\end{equation*}
Since the left hand side is nonpositive, we have
\begin{equation*}
 |\text{supp}(u)| \leq \gamma^{-1} {\int_{\supp(u)} |f| dx}.
 \end{equation*}

For the parabolic case, define the time dependent support set $\mathcal{S}(t):=\text{supp}(u(x,t))$. Differentiating the integral of $u$ over $\mathcal{S}(t)$ and using the boundary conditions (\textit{i.e.}, $u=0$ on $\partial \mathcal{S}(t)$) yields:
\begin{align*}
\frac{d}{dt} \int_{\mathcal{S}(t)} u(x,t) dx &= \int_{\mathcal{S}(t)} u_t dx =  \int_{\mathcal{S}(t)} \nabla \cdot M \nabla u + f - \gamma p(u)  \ dx.
\end{align*}
Because of the divergence theorem and the fact that $M$ is positive definite, we have
\begin{align*}
\frac{d}{dt} \int_{\mathcal{S}(t)} |u(x,t)| dx & \leq\int_{\mathcal{ S}(t)}  |f| dx- \gamma |\mathcal{S}(t)|.
\end{align*}
Integrating the expression in time yields the following bound on the support size:
\begin{align*}
| \supp_{(x,t)} u(x,t)|  \leq \int_{\mathcal{ S}(t)}  |g| dx
+ \iint_{\mathcal{ S}(t)}  |f| dx \, dt.
\end{align*}
 \end{proof}

\end{document}